\documentclass[12pt]{article}
\usepackage{amssymb,latexsym}
\usepackage{amsfonts}

\input epsf

\newcommand{\gl}{{\ensuremath{{\rm gl}}}}
\newcommand{\gloo}{{\ensuremath{{\rm gl}(1\vert 1)}}}
\newcommand{\A}{{\cal A}}
\newcommand{\Ac}{{\widehat{\cal A}}}

\newcommand{\coeff}{{\bf k}}
\newcommand{\sdim}{\ensuremath{\rm sdim}}
\newcommand{\str}{\ensuremath{\rm str}}
\newcommand{\tr}{\ensuremath{\rm tr}}
\newcommand{\id}{\ensuremath{\rm id}}
\newcommand{\Hom}{{\ensuremath{\rm Hom}}}
\newcommand{\End}{{\ensuremath{\rm End}}}
\newcommand{\Ker}{{\ensuremath{\rm Ker}}}
\newcommand{\mod}{{\ensuremath{\rm \, mod \,}}}
\newcommand{\Gal}{{\ensuremath{\rm Gal}}}
\newcommand{\ad}{{\ensuremath{\rm ad}}}

\newcommand{\Sym}{{\ensuremath{\rm Sym}}}

\newcommand{\rank}{{\ensuremath{\rm rank}}}
\newcommand{\op}{{\ensuremath{\rm op}}}

\newcommand{\Anab}{{\cal C}}
\newcommand{\Anabo}{{\ensuremath{\Anab_0}}´}
\newcommand{\Anabc}{{\ensuremath{\widehat{\Anab}}}}

\newcommand{\Bnab}{{\cal D}}

\newcommand{\Wnab}{\ensuremath{W}}
\newcommand{\Wnabt}{\ensuremath{W^\circ\!}}
\newcommand{\Ws}{{\widehat{W}}}
\newcommand{\Wst}{\widehat{\Wnabt}}
\newcommand{\Zb}{\ensuremath{\overline{Z}}}
\newcommand{\G}{{\cal G}}
\newcommand{\V}{{\ensuremath{V_t}}}
\newcommand{\Vo}{{\ensuremath{V_{t_1}}}}
\newcommand{\Vp}{{\ensuremath{V_{t_2}}}}
\newcommand{\Vpp}{{\ensuremath{V_{t_3}}}}
\newcommand{\vv}{{\ensuremath{v_t}}}
\newcommand{\vvo}{{\ensuremath{v_{t_1}}}}
\newcommand{\vvp}{{\ensuremath{v_{t_2}}}}

\newcommand{\ww}{{\ensuremath{w_t}}}
\newcommand{\wwo}{{\ensuremath{w_{t_1}}}}
\newcommand{\wwp}{{\ensuremath{w_{t_2}}}}

\newcommand{\gcap}{{\setbox1=\hbox{\begin{picture}(0,0)%
\includegraphics{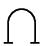}%
\end{picture}%
\setlength{\unitlength}{3947sp}%
\begingroup\makeatletter\ifx\SetFigFont\undefined%
\gdef\SetFigFont#1#2#3#4#5{%
  \reset@font\fontsize{#1}{#2pt}%
  \fontfamily{#3}\fontseries{#4}\fontshape{#5}%
  \selectfont}%
\fi\endgroup%
\begin{picture}(195,184)(2174,-1991)
\end{picture}
}\vcenter{\box1}\!\!}}

\newcommand{\gcup}{{\setbox1=\hbox{\begin{picture}(0,0)%
\includegraphics{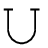}%
\end{picture}%
\setlength{\unitlength}{3947sp}%
\begingroup\makeatletter\ifx\SetFigFont\undefined%
\gdef\SetFigFont#1#2#3#4#5{%
  \reset@font\fontsize{#1}{#2pt}%
  \fontfamily{#3}\fontseries{#4}\fontshape{#5}%
  \selectfont}%
\fi\endgroup%
\begin{picture}(195,183)(2174,-1764)
\end{picture}
}\vcenter{\box1}\!\!\!}}
\newcommand{\Ystd}{{\setbox1=\hbox{\begin{picture}(0,0)%
\includegraphics{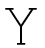}%
\end{picture}%
\setlength{\unitlength}{3947sp}%
\begingroup\makeatletter\ifx\SetFigFont\undefined%
\gdef\SetFigFont#1#2#3#4#5{%
  \reset@font\fontsize{#1}{#2pt}%
  \fontfamily{#3}\fontseries{#4}\fontshape{#5}%
  \selectfont}%
\fi\endgroup%
\begin{picture}(174,192)(2248,-1620)
\end{picture}
}\vcenter{\box1}\!\!\!}}
\newcommand{\Ydot}{{\setbox1=\hbox{\begin{picture}(0,0)%
\includegraphics{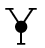}%
\end{picture}%
\setlength{\unitlength}{3947sp}%
\begingroup\makeatletter\ifx\SetFigFont\undefined%
\gdef\SetFigFont#1#2#3#4#5{%
  \reset@font\fontsize{#1}{#2pt}%
  \fontfamily{#3}\fontseries{#4}\fontshape{#5}%
  \selectfont}%
\fi\endgroup%
\begin{picture}(174,192)(2248,-2070)
\end{picture}
}\vcenter{\box1}\!\!\!}}
\newcommand{\Astd}{{\setbox1=\hbox{\begin{picture}(0,0)%
\includegraphics{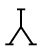}%
\end{picture}%
\setlength{\unitlength}{3947sp}%
\begingroup\makeatletter\ifx\SetFigFont\undefined%
\gdef\SetFigFont#1#2#3#4#5{%
  \reset@font\fontsize{#1}{#2pt}%
  \fontfamily{#3}\fontseries{#4}\fontshape{#5}%
  \selectfont}%
\fi\endgroup%
\begin{picture}(174,192)(2248,-1844)
\end{picture}
}\vcenter{\box1}\!\!}}
\newcommand{\Adot}{{\setbox1=\hbox{
\begin{picture}(0,0)%
\includegraphics{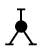}%
\end{picture}%
\setlength{\unitlength}{3947sp}%
\begingroup\makeatletter\ifx\SetFigFont\undefined%
\gdef\SetFigFont#1#2#3#4#5{%
  \reset@font\fontsize{#1}{#2pt}%
  \fontfamily{#3}\fontseries{#4}\fontshape{#5}%
  \selectfont}%
\fi\endgroup%
\begin{picture}(174,192)(2248,-1844)
\end{picture}
}\vcenter{\box1}\!\!}}
\newcommand{\I}{{\setbox1=\hbox{\begin{picture}(0,0)%
\includegraphics{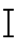}%
\end{picture}%
\setlength{\unitlength}{3947sp}%
\begingroup\makeatletter\ifx\SetFigFont\undefined%
\gdef\SetFigFont#1#2#3#4#5{%
  \reset@font\fontsize{#1}{#2pt}%
  \fontfamily{#3}\fontseries{#4}\fontshape{#5}%
  \selectfont}%
\fi\endgroup%
\begin{picture}(73,152)(2089,-2001)
\end{picture}
}\vcenter{\box1}\!\!}}
\newcommand{\Xt}{{\setbox1=\hbox{\begin{picture}(0,0)%
\includegraphics{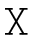}%
\end{picture}%
\setlength{\unitlength}{3947sp}%
\begingroup\makeatletter\ifx\SetFigFont\undefined%
\gdef\SetFigFont#1#2#3#4#5{%
  \reset@font\fontsize{#1}{#2pt}%
  \fontfamily{#3}\fontseries{#4}\fontshape{#5}%
  \selectfont}%
\fi\endgroup%
\begin{picture}(131,144)(2313,-2008)
\end{picture}
}\vcenter{\box1}\!}}

\newcommand{\na}{{\rm na}}
\newcommand{\mercedes}{{\setbox1=\hbox{\begin{picture}(0,0)%
\includegraphics{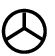}%
\end{picture}%
\setlength{\unitlength}{3947sp}%
\begingroup\makeatletter\ifx\SetFigFont\undefined%
\gdef\SetFigFont#1#2#3#4#5{%
  \reset@font\fontsize{#1}{#2pt}%
  \fontfamily{#3}\fontseries{#4}\fontshape{#5}%
  \selectfont}%
\fi\endgroup%
\begin{picture}(230,228)(486,-790)
\end{picture}
}\vcenter{\box1}\!}}
\newcommand{\nabh}{{\ensuremath{\widehat{\nabla}}}}
\def\Aa[#1,#2,#3]{{A(#1,#2,#3)}}
\newcommand{\indx}{{k}}

\begingroup\makeatletter\ifx\SetFigFont\undefined%
\gdef\SetFigFont#1#2#3#4#5{%
  \reset@font\fontsize{#1}{#2pt}%
  \fontfamily{#3}\fontseries{#4}\fontshape{#5}%
  \selectfont}%
\fi\endgroup%

\def\Y[#1,#2,#3]{{
\setbox1=\hbox{
\begin{picture}(0,0)%
\includegraphics{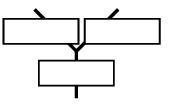}%
\end{picture}%
\setlength{\unitlength}{3947sp}%
\begin{picture}(394,419)(879,67)
\put(1326,379){\makebox(280,0)[cb]{\smash{\SetFigFont{2}{2.0}{\familydefault}{\mddefault}{\updefault}{\color[rgb]{0,0,0}$\nu^{#2}$}%
}}}
\put(1106,179){\makebox(280,0)[cb]{\smash{\SetFigFont{2}{2.0}{\familydefault}{\mddefault}{\updefault}{\color[rgb]{0,0,0}$\nu^{#3}$}%
}}}
\put(931,379){\makebox(280,0)[cb]{\smash{\SetFigFont{2}{2.0}{\familydefault}{\mddefault}{\updefault}{\color[rgb]{0,0,0}$\nu^{#1}$}%
}}}
\end{picture}\qquad
} \vcenter{\box1} } }

\def\Aa[#1,#2,#3]{{
\setbox1=\hbox{
\begin{picture}(0,0)%
\includegraphics{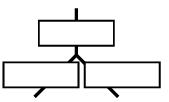}%
\end{picture}%
\setlength{\unitlength}{3947sp}%
\begin{picture}(774,469)(689,-508)
\put(1056,-225){\makebox(0,0)[cb]{\smash{\SetFigFont{2}{2.0}{\familydefault}{\mddefault}{\updefault}{\color[rgb]{0,0,0}$\nu^{#1}$}%
}}}
\put(881,-425){\makebox(0,0)[cb]{\smash{\SetFigFont{2}{2.0}{\familydefault}{\mddefault}{\updefault}{\color[rgb]{0,0,0}$\nu^{#2}$}%
}}}
\put(1276,-425){\makebox(0,0)[cb]{\smash{\SetFigFont{2}{2.0}{\familydefault}{\mddefault}{\updefault}{\color[rgb]{0,0,0}$\nu^{#3}$}%
}}}
\end{picture}\qquad
} \vcenter{\box1} } }

\def\Ip[#1,#2]{{
\setbox1=\hbox{
\begin{picture}(0,0)%
\includegraphics{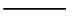}%
\end{picture}%
\setlength{\unitlength}{3947sp}%
\begin{picture}(324,239)(-1061,1097)
\put(-1140,1180){\makebox(0,0)[lb]{\smash{\SetFigFont{7}{8.4}{\familydefault}{\mddefault}{\updefault}{\color[rgb]{0,0,0}$#1$}%
}}}
\put(-715,1180){\makebox(0,0)[lb]{\smash{\SetFigFont{7}{8.4}{\familydefault}{\mddefault}{\updefault}{\color[rgb]{0,0,0}$#2$}%
}}}
\end{picture}
} \vcenter{\box1} } }

\def\Ips[#1,#2]{{
\setbox1=\hbox{
\begin{picture}(0,0)%
\includegraphics{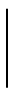}%
\end{picture}%
\setlength{\unitlength}{4144sp}%
\begin{picture}(32,555)(251,164)
\put(251,164){\makebox(0,0)[lb]{\smash{\SetFigFont{8}{9.6}{\familydefault}{\mddefault}{\updefault}{\color[rgb]{0,0,0}$#2$}%
}}}
\put(251,634){\makebox(0,0)[lb]{\smash{\SetFigFont{8}{9.6}{\familydefault}{\mddefault}{\updefault}{\color[rgb]{0,0,0}$#1$}%
}}}

\end{picture}
} \vcenter{\box1} } }

\def\Yp[#1,#2,#3]{{
\setbox1=\hbox{
\begin{picture}(0,0)%
\includegraphics{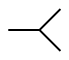}%
\end{picture}%
\setlength{\unitlength}{3947sp}%
\begin{picture}(380,321)(951,-136)
\put(1331,-136){\makebox(0,0)[lb]{\smash{\SetFigFont{6}{9.6}{\familydefault}{\mddefault}{\updefault}{\color[rgb]{0,0,0}$#2$}%
}}}
\put(1331, 89){\makebox(0,0)[lb]{\smash{\SetFigFont{8}{9.6}{\familydefault}{\mddefault}{\updefault}{\color[rgb]{0,0,0}$#3$}%
}}}
\put(996,-23){\makebox(0,0)[lb]{\smash{\SetFigFont{8}{9.6}{\familydefault}{\mddefault}{\updefault}{\color[rgb]{0,0,0}$#1$}%
}}}
\end{picture}
} \vcenter{\box1} } }

\def\H[#1,#2,#3,#4]{{
\setbox1=\hbox{
\begin{picture}(0,0)%
\includegraphics{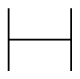}%
\end{picture}%
\setlength{\unitlength}{3947sp}%
\begin{picture}(337,550)(501,-246)
\put(801,-246){\makebox(0,0)[lb]{\smash{\SetFigFont{8}{9.6}{\familydefault}{\mddefault}{\updefault}{\color[rgb]{0,0,0}$#4$}%
}}}
\put(501,199){\makebox(0,0)[lb]{\smash{\SetFigFont{8}{9.6}{\familydefault}{\mddefault}{\updefault}{\color[rgb]{0,0,0}$#1$}%
}}}
\put(801,199){\makebox(0,0)[lb]{\smash{\SetFigFont{8}{9.6}{\familydefault}{\mddefault}{\updefault}{\color[rgb]{0,0,0}$#2$}%
}}}
\put(501,-246){\makebox(0,0)[lb]{\smash{\SetFigFont{8}{9.6}{\familydefault}{\mddefault}{\updefault}{\color[rgb]{0,0,0}$#3$}%
}}}
\end{picture}}
\vcenter{\box1} }}

\def\Ha[#1,#2,#3,#4,#5]{{
\setbox1=\hbox{
\begin{picture}(0,0)%
\includegraphics{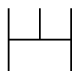}%
\end{picture}%
\setlength{\unitlength}{3947sp}%
\begin{picture}(337,550)(501,-246)
\put(501,-246){\makebox(0,0)[lb]{\smash{\SetFigFont{8}{9.6}{\familydefault}{\mddefault}{\updefault}{\color[rgb]{0,0,0}$#4$}%
}}}
\put(801,-246){\makebox(0,0)[lb]{\smash{\SetFigFont{8}{9.6}{\familydefault}{\mddefault}{\updefault}{\color[rgb]{0,0,0}$#5$}%
}}}
\put(801,199){\makebox(0,0)[lb]{\smash{\SetFigFont{8}{9.6}{\familydefault}{\mddefault}{\updefault}{\color[rgb]{0,0,0}$#3$}%
}}}
\put(501,199){\makebox(0,0)[lb]{\smash{\SetFigFont{8}{9.6}{\familydefault}{\mddefault}{\updefault}{\color[rgb]{0,0,0}$#1$}%
}}}
\put(651,199){\makebox(0,0)[lb]{\smash{\SetFigFont{8}{9.6}{\familydefault}{\mddefault}{\updefault}{\color[rgb]{0,0,0}$#2$}%
}}}
\end{picture}}
\vcenter{\box1} }}

\newtheorem{lemma}{Lemma}
\newtheorem{prop}[lemma]{Proposition}
\newtheorem{remark}[lemma]{Remark}
\newtheorem{theorem}[lemma]{Theorem}

\newtheorem{coro}[lemma]{Corollary}
\newtheorem{defi}[lemma]{Definition}
\newtheorem{fact}[lemma]{Fact}

\newenvironment{proof}{{\bf Proof:}}{$\Box$\medskip}
\newenvironment{skoproof}{{\bf Sketch of proof:}}{$\Box$\medskip}

\newcommand{\N}{{\ensuremath{\mathbb N}}}
\newcommand{\Z}{{\ensuremath{\mathbb Z}}}
\newcommand{\Q}{{\ensuremath{\mathbb Q}}}
\newcommand{\R}{{\ensuremath{\mathbb R}}}
\newcommand{\C}{{\ensuremath{\mathbb C}}}

\def\color[#1]#2{}
\def\epsfigbracket[#1=#2]{\epsffile{#2}}

\def\includegraphics#1{\epsffile{#1}}

\begin{document}
\parindent0cm

\title{
\bf The Drinfeld associator of $\gloo$}
\author{Jens Lieberum\\
\footnotesize Mathematisches Institut, Rheinsprung 21, CH-4051
Basel, lieberum@math.unibas.ch}

\date{}
\maketitle


\begin{abstract}
We determine explicitly a rational even
Drinfeld associator~$\overline{\Phi}$ 
in a completion of the universal enveloping
algebra of the Lie superalgebra~$\gl(1\vert 1)^{\oplus 3}$. 
More generally,
we define a new algebra of trivalent diagrams that has a
unique even horizontal group-like  
Drinfeld associator~$\Phi$. The associator~$\Phi$ 
is mapped to~$\overline{\Phi}$ by a weight system.
As a related result of independent interest, we show how
O.\ Viro's generalization~$\underline{\Delta}^1$  
of the multi-variable Alexander polynomial 
can be obtained from the universal Vassiliev invariant of trivalent graphs.
We determine~$\Phi$ by using the invariant
$\underline{\Delta}^1(\mercedes)$ of a planar tetrahedron~$\mercedes$.
\end{abstract}

\bigskip



\bigskip

\section*{Introduction} { 
Building on concepts of Mac Lane and Kohno,
Drinfeld introduced associators in~1989/90 in~\cite{Dr1} 
in context with a weakened version of the 
coassociativity axiom of Hopf algebras and quasitriangular Hopf 
algebras. 
Around the same time, the development of a systematic 
approach to knot theory
started with the concept 
of Vassiliev invariants (\cite{Vas}). The relation between 
Vassiliev invariants and Drinfelds work was established by 
Kontsevichs analytic construction of a universal Vassiliev 
invariant~$Z$ and its algebraic description that requires the 
existence of a Drinfeld associator. 
The most useful and most convenient associators for topological applications
are horizontal even group-like Drinfeld associators.

In the main result of this paper 
we determine explicitly an even horizontal group-like 
associator~$\overline{\Phi}$ in a completion of
$U(\gl(1\vert 1))^{\otimes 3}$. 
More generally,
we define an algebra~$\Anab(n)$ over 
$\Lambda(n)=\Q[d_1^{\pm 1}, \ldots,
d_n^{\pm 1}]$ that is generated by
trivalent diagrams 
on~$n$ strings modulo some graphical relations.
In a completion~$\Anabc^0(3)$ of a~$\Q$-subalgebra of~$\Anab(3)$ 
there exists a unique Drinfeld
associator of the form

\begin{equation}\label{e:Phiprime}
\Phi\ =\ \exp\left(F\cdot [t^{12}, t^{23}]\right)\in\Anabc^0(3),
\end{equation}

where $F\in\Q[[C,D,E]]\subset\Q[[d_1,d_2,d_3]]$ is a formal power series that starts with

\begin{eqnarray*}
F & = & \frac{1}{24}-\frac{C+4D}{5760}+ \frac{4 C^2+36C D+48
D^2-31 E}{2903040}\samepage\\\samepage
                    & & -\frac{6C^3+96C^2 D+240 C D^2+192 D^3+13 C E-184
                    DE}{464486400}+\ldots
\end{eqnarray*}

and $C=d_1 d_3-d_2(d_1+d_2+d_3)$, 
$D=(d_1+d_2)^2+(d_2+d_3)^2$, $E=(d_1+d_2)^2(d_2+d_3)^2$. 
The main result of this paper is the following theorem.

\begin{theorem}\label{t:F}
The 
series~$F$ is the unique solution of the equation 

$$
\cosh(Fk)+d_1d_3\frac{\sinh(Fk)}{k}=\sqrt{\frac{\varphi(d_2)\varphi(d_1+d_2+d_3)}
{\varphi(d_1+d_2)\varphi(d_2+d_3)}}
$$

where $\varphi(x)=2\sinh(x/2)/x$ and $k^2=-d_1d_2d_3(d_1+d_2+d_3)$.
\end{theorem}

Due to a computation of P.\ Vogel the solution~$F$ 
is given
explicitly by

\begin{eqnarray}
F & = & X\, \Psi(d_1d_3X, -d_2(d_1+d_2+d_3)X),\quad \mbox{where}\label{e:fF}\\
\Psi(u,v) & = & \sum_{p,q=0}^\infty \left({-1/2\atop p}\right)\left({-1/2\atop q}\right)
\frac{u^pv^q(u+2)^{p+q+1}}{2(p+q+1)}\quad \mbox{and}\label{e:defPsi}\\
X & = & d_1^{-1}d_3^{-1}\left(\sqrt{\frac{\varphi(d_2)\varphi(d_1+d_2+d_3)}
{\varphi(d_1+d_2)\varphi(d_2+d_3)}}-1\right).\label{e:defX}
\end{eqnarray}

There is a map from~$\Anabc^0(3)$ to 
a completion of $U(\gl(1\vert 1))^{\otimes 3}$ that maps~$\Phi$ to
the associator~$\overline{\Phi}$ mentioned before.
The associator~$\Phi\in\Anabc^0(3)$ is the unique
even horizontal group-like Drinfeld associator inside of
$\Anabc_\coeff^0(3)$. The image of~$\Phi$ in a quotient~$\Anabc^1(3)$ is the
unique even group-like Drinfeld associator inside of this quotient 
(see Theorem~\ref{t:existF2}).

The quantum supergroup $U_q(\gl(1\vert 1))$ has been used by O.\ Viro to extend
the 
multi-variable Alexander polynomial of links to an invariant~$\underline{\Delta}^1$ of 
colored embedded trivalent graphs (\cite{Vir}). 
In general, the invariant
of a trivially embedded tetrahedron colored by representations of a quantum
group provides the relation 
between $R$-matrices and $6j$-symbols (\cite{Tu2}). 
Although this construction 
neither extends in its full generality to quantum supergroups nor to 
versions of Vassiliev invariants for $3$-manifolds
the value~$\underline{\Delta}^1(\mercedes)$ 
and its relation to the Kontsevich integral~$Z$
turned out to be useful to determine~$\overline{\Phi}$ explicitly.
The translation between~$\underline{\Delta}^1(\mercedes)$ 
and~$\overline{\Phi}$ is in the spirit of a general
connection between well-behaved invariants of trivalent graphs and associators 
that is investigated by D.\ Bar-Natan and D.\ Thurston. 
We extend the relation between~$\underline{\Delta}^1$ and the Kontsevich integral~$Z$
from~$\mercedes$ to arbitrary trivalent graphs (Theorem~\ref{t:nabvir}). 
This is a result of independent interest. It
generalizes an unpublished proof of A.\ Vaintrob who related the multi-variable Alexander
polynomial of links to the Kontsevich integral.
With the standard definition the Kontsevich integral of trivalent graphs (\cite{MuO})
Theorem~\ref{t:nabvir} 
would only hold up to a factor that depends on the 
colored graph but not on its embedding. 
In order to avoid this factor, to simplify
computations, and to emphasize the roles played by cyclic orientations of vertices and 
half-framings
we introduce a different normalization of the Kontsevich integral of oriented trivalent
graphs (Theorem~\ref{t:Ztriv}).


 
The paper is organized as follows. In Section~\ref{s:assoA} we recall definitions
and properties of a module~$\A(\Gamma)$ of 
trivalent diagrams on an oriented unitrivalent graph~$\Gamma$ 
and of Drinfeld associators in a completion~$\Ac(3)$ of~$\A(\Gamma_3)$ where
$\Gamma_n$ consists of~$n$ intervals. 
In Section~\ref{s:basAnab} we introduce the module~$\Anab(\Gamma)$
and investigate its structure for~$\Gamma=\Gamma_n$
(Theorem~\ref{t:Bnabgen} and Corollary~\ref{c:Anabfree}). 
Section~\ref{s:assoC} contains 
results about Drinfeld associators in~$\Anabc^0(3)$ and in a quotient~$\Anabc^1(3)$ 
of~$\Anabc^0(3)$ 
(Theorem~\ref{t:existF2}). In particular, we establish 
the existence and uniqueness of the series~$F$ and deduce equations~(\ref{e:fF})
to~(\ref{e:defX}) from Theorem~\ref{t:F}.
In Sections~\ref{s:rep1} to~\ref{s:lemmas} we prepare the proof 
of Theorem~\ref{t:F} that will be given in Section~\ref{s:tetra2}. 
Sections~\ref{s:rep1} to~\ref{s:tetra1} are also used in 
Section~\ref{s:MultAlexViro} where we relate the Kontsevich integral~$Z$ to
Viro's Alexander invariant~$\underline{\Delta}^1$ (Theorem~\ref{t:nabvir}).

\subsection*{Acknowledgements}

I would like to thank D.\ Bar-Natan, C.\ Kassel, G.\ Masbaum, 
D.\ Thurston, and P.\ Vogel for helpful
discussions, and O.\ Viro for writing~\cite{Vir} and for sending me 
a preliminary version of that paper.

}

\section{Drinfeld associators in~$\Ac_\coeff(3)$}\label{s:assoA}

A graph is called unitrivalent if all of
its vertices have valency one or three.
Let~$\Gamma$ be a unitrivalent graph
with oriented edges and cyclically oriented vertices. 
As an exception, we allow circles as connected components of~$\Gamma$ 
that we consider as a single oriented edge without vertex. 
The graph~$\Gamma$ may have multiple edges between vertices. When~$\Gamma$ has no
univalent vertex, we call it a trivalent graph.
Let~$V(G)$
(resp.\ $E(G)$) be the set of vertices (resp.\ edges)
of a unitrivalent graph~$G$.
A {\em trivalent diagram}~$D$ with {\em skeleton}~$\Gamma$ is a unitrivalent graph~$G$ whose
univalent vertices are glued to~$\Gamma\setminus V(\Gamma)$ by an injective gluing map.
The unitrivalent graph~$G$ has the following properties:
trivalent vertices of~$G$ are cyclically oriented, but in contrast to~$\Gamma$, 
edges of~$G$ are not oriented. In addition, we require that
each connected component of~$G$ has at least one univalent vertex.

We represent a trivalent diagram graphically by a generic picture 
of~$\Gamma\cup G$ in the plane. We use thicker lines
to draw~$\Gamma$ than we use for~$G$. We assume that cyclic
orientations of oriented vertices are always counterclockwise.
When it is of importance we indicate orientations of edges
of~$\Gamma$ by arrows, and
we include the
names~$i,j,\ldots$ of edges (resp.\ vertices) in our graphical representation of
trivalent diagrams by writing them close to the corresponding
edges (resp.\ vertices). 
Homeomorphisms $h:\Gamma\cup G\longrightarrow \Gamma\cup G'$
between trivalent diagrams on~$\Gamma$ have to respect
orientations of edges and vertices and induce a homeomorphism of~$\Gamma$ that is 
homotopic to the identity. By abuse of language we call the homeomorphism 
class of a
trivalent diagram simply trivalent diagram. 
For a vertex~$v\in V(\Gamma)$ and an edge $i\in E(\Gamma)$ that
is incident to~$v$, we define~$s_{v,i}\in\{\pm 1\}$ by~$s_{v,i}=1$ if the edge~$i$ is
oriented towards~$v$ and~$s_{v,i}=-1$ otherwise.
We specify relations
between trivalent diagrams by using graphical representations of
the part where these diagrams differ.

\begin{defi}\label{d:defA}
Let $\A(\Gamma)$ be the $\Q$-vector space generated by trivalent
diagrams on~$\Gamma$ modulo the relations~$(STU)$ and $(InvV)$.
\end{defi}
\begin{eqnarray*}
(STU) & & \setbox1=\hbox{\input{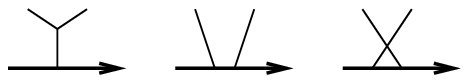}}\vcenter{\box1}\!,\\ (InvV) & & s_{v,i}\setbox1=\hbox{\input{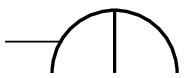}}\vcenter{\box1}\!+s_{v,j}\setbox1=\hbox{\input{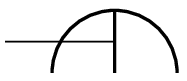}}\vcenter{\box1}\!+
s_{v,k}\setbox1=\hbox{\input{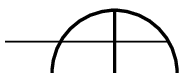}}\vcenter{\box1}\!\ = \ 0.
\end{eqnarray*}
 
The signs on the right side of relation~$(STU)$ depend on the cyclic order of the trivalent
vertex and on the orientation of the shown part of~$\Gamma$ in this relation.
The degree of a trivalent diagram~$D=\Gamma\cup G$ is defined by~$\deg(D)=(1/2)\#V(G)$.
This definition induces
a grading on~$\A(\Gamma)$.
It is well-known (see \cite{BN1}) that the relations~$(IHX)$
and~$(AS)$ below are consequences of the $(STU)$-relation.

$$ (IHX)\quad \setbox1=\hbox{\input{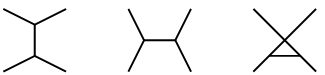}}\vcenter{\box1}\!,\qquad (AS)\quad \setbox1=\hbox{\input{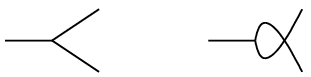}}\vcenter{\box1}\!.$$

A basis of $\A(\Gamma)$ is not known explicitly. 
Let~$\Gamma_n=\coprod_{i=1}^n I_i$ be the disjoint union of
oriented intervals~$I_i$. The bijection between~$\{1,\ldots,n\}$ and
the intervals of~$\Gamma_n$ is a part
of the definition of~$\Gamma_n$.
The vector spaces~$\A(n)=\A(\Gamma_n)$ are particularly interesting for several
reasons. One reason is that for connected trivalent graphs~$\Gamma$
there exist isomorphisms~$\A(\Gamma)\cong\A(b_1(\Gamma))$ where~$b_1$ is the first
Betti number. A second reason is that~$\A(n)$ is an algebra in the following way.
We represent trivalent diagrams
on~$\Gamma_n$ graphically in the strip~$\R\times [0,1]$ 
such that~$\partial I_i=\{i\}\times\{0,1\}$ and these intervals
direct from $(i,1)$ to~$(i,0)$. We say that~$\{1\}\times
\{1,\ldots,n\}$ (resp.\ $\{0\}\times \{1,\ldots,n\}$) is the upper
(resp.\ lower) boundary of a trivalent diagram on~$\Gamma_n$. The 
product~$ab$ of trivalent diagrams~$a$, $b$ is induced by gluing
the lower boundary points of~$a$ to the corresponding upper
boundary points of~$b$. 
The $1$-element of~$\A(n)$ is the diagram~$\Gamma_n$.

For a commutative~$\Q$-algebra~$\coeff$ 
we define the $\N$-graded $\coeff$-module
$$\A_\coeff(\Gamma):=\coeff\otimes_\Q\A(\Gamma)$$
and denote its homogeneous 
components~$\coeff\otimes_\Q\A(\Gamma)_i$ by~$\A_\coeff(\Gamma)_i$. 
The case~$\coeff\not=\Q$ will only 
be important in this section and in 
Section~\ref{s:assoC}.
When~$\coeff=\Q$ we omit the symbol~$\coeff$.
We define the
completion~$\Ac_\coeff(\Gamma)$ by

$$ \Ac_\coeff(\Gamma)=\prod_{i=0}^\infty\Ac_\coeff(\Gamma)_i. $$

We represent elements of completions of graded vector spaces by formal
power series~$\sum_{i=k}^\infty a_i$ with homogeneous
elements~$a_i$ of degree~$i$. We consider completions 
as metric spaces (and in
particular as topological spaces) by

$$d\left(\sum_{i=k}^\infty a_i,\sum_{i=k}^\infty b_j\right)
=\sum_{i=k}^\infty \delta_{a_i,b_i}2^{-i}$$ 

where~$\delta_{a_i,b_i}\in\{0,1\}$ is~$0$ iff
$a_i=b_i$.
The algebra structure on~$\A_\coeff(n):=\A_\coeff(\Gamma_n)$
extends in a unique way to a topological algebra structure
on~$\Ac_\coeff(n):=\Ac_\coeff(\Gamma_n)$.

Define continous 
$\coeff$-linear maps~$\Delta_i:\Ac_\coeff(n)\longrightarrow\Ac_\coeff(n+1)$ ($i=0,\ldots, n+1$), 
where~$\Delta_i(D)$ ($i=1,\ldots,n$)
is obtained from the trivalent diagram~$D$ by
replacing the $i$-th interval of~$D$ by two copies and by summing over all ways 
of lifting the univalent vertices that are glued to the $i$-th interval 
to the copies of that interval. This sum has~$2^\ell$ terms when 
$\ell$ univalent vertices are glued
to the $i$-th interval of~$D$. 
We label the new intervals in~$\Delta_i(D)$ by~$i$ and $i+1$ and replace labels
$j$ with $j>i$ by $j+1$. Define~$\Delta_{n+1}(D)$ as 
the union of~$D$ with a new skeleton component labeled~$n+1$, and 
define 
$\Delta_0(D)$
by first adding to~$D$ a new skeleton component labeled~$0$ followed 
by replacing
all labels~$i$ of the skeleton components by~$i+1$.

Define $\epsilon_i:\Ac_\coeff(n)\longrightarrow\Ac_\coeff(n-1)$ by $\epsilon_i(D)=0$ if
the $i$-th interval of~$D$ contains a univalent vertex of~$D\setminus\Gamma$, 
and by deleting the $i$-th
interval and by replacing labels~$j>i$ by~$j-1$ otherwise.

Let~$f:\{1,\ldots,m\}\longrightarrow\{1,\ldots n\}$ be an injective map. 
For $a\in\Ac(m)$ we
denote by $a^{f(1)\ldots f(m)}\in\Ac(n)$ the element 
obtained by applying~$f$ to the labels of the skeleton components 
of trivalent diagrams in~$a$ and by adding intervals with labels 
$i\in\{1,\ldots,n\}\setminus f(\{1,\ldots,m\})$ to the skeleton of these diagrams.
Let $t^{ij}$ be the unique trivalent
diagram on~$\Gamma_n$ of degree~$1$ where the $i$-th edge
of~$\Gamma_n$ is connected to the $j$-th edge of~$\Gamma_n$ by a
single edge.
 Define~$R\in\Ac_\coeff(2)$ by
$R=\exp(t^{12}/2)$. Now we are ready to define a Drinfeld associator.

\begin{defi}\label{d:Phi}
A Drinfeld associator~$\Phi$ in~$\Ac_\coeff(3)$ is a solution of 
equations (DA1)-(DA4) in $\Ac_\coeff(n)$ ($n=4,3,3,2$).

\begin{eqnarray*}
(DA1) & & \Delta_0(\Phi)\circ\Delta_2(\Phi)\circ\Delta_4(\Phi)=\Delta_3(\Phi)\circ\Delta_1(\Phi),\\
(DA2) & & \Delta_1(R^{12})=\Phi^{312}\circ R^{13}\circ({\Phi^{132}}
)^{-1}
\circ R^{23}\circ\Phi\\
(DA3) & & \Phi\circ\Phi^{321}=1\\
(DA4) & & \epsilon_1(\Phi)=\epsilon_2(\Phi)=\epsilon_3(\Phi)=1\\
\end{eqnarray*}
%
\end{defi}

Equation~$(DA4)$ implies that there exists a unique~$P\in\Ac_\coeff(3)$ such that
$\Phi=\exp(P)$.
Let ${\cal P}(\Ac_\coeff(n))$ be the closed 
$\coeff$-submodule of~$\Ac_\coeff(n)$ consisting of series that
involve only trivalent diagrams~$D$ 
with the property that $D\setminus\Gamma_n$ is non-empty and connected.
When~$P\in {\cal P}(\Ac_\coeff(n))$ we say that~$\exp(P)$ is {\em group-like}.
A Drinfeld associator is called {\em horizontal}, if it lies in the closed 
subalgebra of~$\Ac_\coeff(3)$ generated by~$t^{12}$ and~$t^{23}$.
A horizontal Drinfeld associator uniquely determines a formal 
series~$S\in \coeff\langle\langle A,B\rangle\rangle$ in non-commutative, associative
indeterminates~$A$ and~$B$ such that
$\Phi=S(t^{12}, t^{23})$ (see~\cite{BN3}, Fact~9 and Corollary~4.4). 
When~$\Phi$ is group-like
the series~$S$ 
satisfies

\begin{eqnarray*}
(DA5) & & \mbox{$S=\exp(p)$ for a Lie series~$p$ in~$A$ and $B$ over~$\coeff$.}
\end{eqnarray*}

We say that a Drinfeld
associator~$\Phi$ is {\em even}, if~$\Phi=\sum_{i=0}^\infty a_{2i}$ 
with~$a_{2i}\in\A_\coeff(3)_{2i}$. For the following fact 
see~\cite{Dri}, Theorem A'' and~\cite{BN4}, Corollary~4.2.

\begin{fact}\label{f:existPhi}
There exists an even horizontal group-like Drinfeld associator in~$\Ac(3)$.
\end{fact}

Associators in~$\Ac_\coeff(3)$ are not unique, but
for two associators~$\Phi_1$, $\Phi_2$ there exists an 
element~$T\in\Ac_\coeff(2)$ satisfying~$\epsilon_1(T)=\epsilon_2(T)=1$ and~$T^{21}=T$ such
that

\begin{equation}\label{e:twist}
\Phi_1=\Delta_0(T)\Delta_2(T)\Phi_2\Delta_1(T^{-1})\Delta_3(T^{-1})
\end{equation}

(see Theorem~8 of~\cite{LeM}). 
We say that $\Phi_1$ and $\Phi_2$ are related by a {\em twist}~$T$.
Group-like (resp.\ even) associators are related by group-like 
(resp.\ even) twists. The relation between horizontal
associators is more involved~:
the proofs of Fact~\ref{f:existPhi} rely on an action 
of the formal Grothendieck-Teichm\"uller group 
on horizontal group-like Drinfeld associators and on the existence of
a horizontal group-like Drinfeld associator~$\Phi_{KZ}\in\Ac_\C(3)$.

\section{The~$\Lambda(n)$-algebras $\Anab(n)$ and $\Bnab(n)$} 
\label{s:basAnab}


Let $\Lambda(\Gamma)$ be the commutative $\Q$-algebra generated by
elements~$d_i, d_i^{-1}$ ($i\in E(\Gamma)$) modulo
relations~$d_id_i^{-1}=1$ for each $i\in E(\Gamma)$ and
$s_{v,i}d_i+s_{v,j}d_j+s_{v,k}d_k=0$ for each trivalent vertex~$v$ of~$\Gamma$
that is incident to the three edges $i, j, k\in E(\Gamma)$.
Let us investigate the structure of~$\Lambda(\Gamma)$.
For $e\in E(\Gamma)$
let~$\alpha_e\in H^1(\Gamma,\partial\Gamma,\Z)=:H$ 
be given by the~$1$-cocycle that evaluates on a $1$-chain~$c$
to the coefficient of~$e$ in~$c$. Then there exists a homomorphism 
from~$H$ into the additive group of~$\Lambda(\Gamma)$ that sends~$\alpha_e$ 
to~$d_e$. It is easy to see that~$\Lambda(\Gamma)=\{0\}$ iff~$\alpha_e=0$ for
some edge~$e$ of~$\Gamma$ (or more explicitly, $\Lambda(\Gamma)=\{0\}$ iff~$\Gamma$ has an
edge~$e$ such that
the number of connected components of~$\Gamma$
increases when we cut $\Gamma$ at a point~$p$ in the interior of~$e$, 
and at least one of the new connected components 
contains no univalent
vertex besides~$p$). When~$\Lambda(\Gamma)\not=\{0\}$ then
$\Lambda(\Gamma)$ is isomorphic to a localization of 
a~$\Q$-algebra of polynomials in~$\rank(H)$ 
indeterminates.

\begin{defi}
Let $\Anab(\Gamma)$ be the quotient of $\Lambda(\Gamma)\otimes_\Q\A(\Gamma)$ by
the relations $(CL1A)$ and $(CL2A)$.
\begin{eqnarray*}
(CL1A) & & \setbox1=\hbox{\input{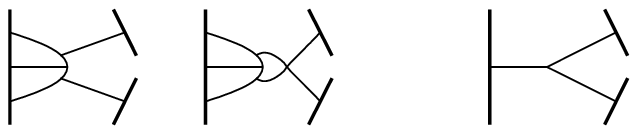}}\vcenter{\box1}\!,\\[3mm] (CL2A) & & \setbox1=\hbox{\input{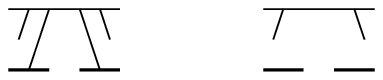}}\vcenter{\box1}\!.
\end{eqnarray*}
\end{defi}

Both relations in the definition of~$\Anab(\Gamma)$ relate 'coefficients' to 'legs' 
of trivalent diagrams what explains the letters~$C$ and~$L$ in the names of the 
relations. The following relations are consequences of the definition of~$\Anab(\Gamma)$.

\begin{eqnarray*}
(LS) & & \setbox1=\hbox{\input{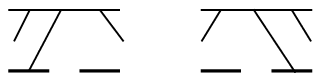}}\vcenter{\box1}\!\\
(IntV) & & \setbox1=\hbox{\begin{picture}(0,0)%
\includegraphics{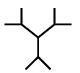}%
\end{picture}%
\setlength{\unitlength}{4144sp}%
\begingroup\makeatletter\ifx\SetFigFont\undefined%
\gdef\SetFigFont#1#2#3#4#5{%
  \reset@font\fontsize{#1}{#2pt}%
  \fontfamily{#3}\fontseries{#4}\fontshape{#5}%
  \selectfont}%
\fi\endgroup%
\begin{picture}(329,303)(439,98)
\end{picture}
}\vcenter{\box1}\!\ =\ \setbox1=\hbox{\begin{picture}(0,0)%
\includegraphics{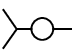}%
\end{picture}%
\setlength{\unitlength}{4144sp}%
\begingroup\makeatletter\ifx\SetFigFont\undefined%
\gdef\SetFigFont#1#2#3#4#5{%
  \reset@font\fontsize{#1}{#2pt}%
  \fontfamily{#3}\fontseries{#4}\fontshape{#5}%
  \selectfont}%
\fi\endgroup%
\begin{picture}(339,204)(574,-253)
\end{picture}
}\vcenter{\box1}\!\ =\ 0
\end{eqnarray*}

Relation~$(LS)$ ('leg slide') is implied by the invertibility of the elements~$d_i$ 
and by relation~$(CL2A)$. In context with relation~$(STU)$, relation~$(LS)$
is equivalent to relations~$(IntV)$ ('internal vertex').

The definition of~$\deg(D)$ for a trivalent diagram~$D$ 
and~$\deg(d_e)=1$ ($e\in E(\Gamma)$)
induce a $\Z$-grading of~$\Anab(\Gamma)$. 
Consider the $\N$-graded $\Q$-subalgebra~$\Lambda^+(\Gamma)$ of~$\Lambda(\Gamma)$ generated by
elements~$d_e$ ($e\in E(\Gamma)$). Define

\begin{equation}\label{e:subalg}
\Anab^+(\Gamma)=\tau(\Lambda^+(\Gamma)\otimes\A(\Gamma))\quad\mbox{and}\quad
\Anab^0(\Gamma)=\tau(\A(\Gamma))
\end{equation}

where $\tau:\Lambda(\Gamma)\otimes\A(\Gamma)\longrightarrow\Anab(\Gamma)$ denotes the
canonical projection. 
Let~$\Anab(n)=\Anab(\Gamma_n)$ and~$\Lambda(n)=\Lambda(\Gamma_n)$.  
We have~$\Lambda(n)=\Q[d_1^{\pm 1},\ldots, d_n^{\pm 1}]$. 
There exists a unique structure of a~$\Lambda(n)$-algebra on~$\Anab(n)$ 
such that the map~$\tau$ is a homomorphism of rings. The space~$\Anab^+(n):=
\Anab^+(\Gamma_n)$ (resp.\ $\Anab^0(n):=\Anab^0(\Gamma_n)$) is 
a~$\Lambda^+(n):=\Lambda^+(\Gamma_n)$-subalgebra (resp.\ $\Q$-subalgebra).

A unitrivalent diagram on a set~${\cal M}$ is a
unitrivalent graph with oriented trivalent vertices and at least one univalent vertex
on each connected component together with an assignment of a label in~${\cal M}$ to
each univalent vertex.
Examples of unitrivalent diagrams can be obtained from trivalent diagrams~$D$
on~$\Gamma_n$ by labeling the univalent vertices
of~$D\setminus\Gamma_n$ according to the connected components
of~$\Gamma_n$. 

\begin{defi}
Let $\Bnab(n)$ be the $\Lambda(n)$-module generated by
unitrivalent diagrams on the set~$\{1,\ldots, n\}$ modulo the relations
$(IHX)$, $(AS)$, $(IntV)$, $(CL1B)$, $(CL2B)$.

\begin{eqnarray*} 
(CL1B) & &
\Ha[i,i,i,j,k]\ = \ d_i^2\Yp[i,j,k] \\[4mm] (CL2B) & & \setbox1=\hbox{\input{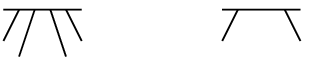}}\vcenter{\box1}\!
\end{eqnarray*}
\end{defi}

Relation~$(CL1B)$ concerns connected components of unitrivalent diagrams 
because all univalent vertices of the diagrams in the relation have labels.
Relation~$(CL2B)$ can be applied to parts of connected components of unitrivalent diagrams.

Let~$D$ be a unitrivalent diagram
on~$\{1,\ldots,n\}$. Let $W_i$ be the set of all linear orders on
the set of univalent vertices of~$D$ labeled~$i$ with the following property~:
for all connected components~$C$ of~$D$ and all univalent vertices $l_1<l_2$ of $C$
labeled~$i$ a univalent vertex~$l_3$ labeled~$i$ with $l_1<l_3<l_2$ also
belongs to~$C$. We define $W=W_1\times\ldots \times W_n$ and

$$\chi (D)=\frac{1}{\# W}\sum_{w \in W} D_w,$$

where $D_w$ is the trivalent diagram on~$\Gamma_n$ obtained by
gluing the univalent vertices labeled~$i$ of $D$ to the $i$-th interval
of~$\Gamma_n$ according to the order~$w$.

\begin{prop}\label{p:PBW}
The definition of~$\chi$ induces an isomorphism
of $\Z$-graded $\Lambda(n)$-mo\-dules $\Bnab(n)\longrightarrow\Anab(n)$.
\end{prop}
\begin{proof}
First verify that~$\chi$ is well-defined: it is clear 
that~$\chi$ is compatible with relations~$(AS)$, $(IHX)$, $(IntV)$ and
it is easy to see that~$\chi$ is compatible with relation~$(CL2B)$.
By relation~$(LS)$ the second vertex\footnote{Notice that our definition
of~$\chi$ is slightly different from the standard definition of~$\chi$ because the standard
definition is not compatible with relation~$(CL1B)$.} on the interval~$i$ in 
relation~$(CL1A)$ can be moved freely on that interval. We use this to
see that~$\chi$ is compatible with
relation~$(CL1B)$ in the case~$j\not=i\not=k$. 
The compatibility of~$\chi$ with relation~$(AS)$ implies compatibilby with relation~$(CL1B)$ 
in the cases~$i=j$ and~$i=k$.

When $D$ and~$D'$ are trivalent diagrams on~$\Gamma_n$ that differ only by the order
of their univalent vertices on~$\Gamma_n$, then by relation~$(STU)$ the element
$D-D'\in\Anab(n)$ can be expressed in terms of diagrams~$D''$ such that~$D''\setminus\Gamma_n$
has more trivalent vertices than
$D\setminus\Gamma_n$. This implies that~$\chi$ is surjective.
The proof of the injectivity
of~$\chi$ is more difficult and similar to the proof of Theorem~8 in~\cite{BN1}.
\end{proof}

Let

$$
G(n)=\left\{\left.\Ips[a,b]\ , \Yp[1,c,d]\ , \H[1,1,e,f]\ ,\, \setbox1=\hbox{\input{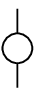}}\vcenter{\box1}\! \ \right\vert\,
1\leq a\leq b\leq n,
2\leq c<d\leq n,
2\leq e\leq f\leq n\right\}.
$$

The disjoint union of diagrams turns~$\Bnab(n)$ into a commutative
$\Lambda(n)$-algebra whose~$1$-element is the empty diagram. 

\begin{theorem}\label{t:Bnabgen}
The commutative $\Lambda(n)$-algebra~$\Bnab(n)$ is freely generated by~$G(n)$.
\end{theorem}

The proof of Theorem~\ref{t:Bnabgen} will occupy the rest of this section. 
We obtain the following corollary of Theorem~\ref{t:Bnabgen} by using
Proposition~\ref{p:PBW} and the ascending filtration of~$\Anab(n)$ 
defined for trivalent diagrams~$D$ by 
the number of connected components of~$D\setminus\Gamma_n$.

\begin{coro}\label{c:Anabfree}
The $\Lambda(n)$-module $\Anab(n)$ is free. For any order on $G(n)$ a basis
of~$\Anab(n)$ is given by ordered monomials in~$\chi(G(n))$. 
\end{coro}

By Corollary~\ref{c:Anabfree} 
the~$\Lambda^+(n)$-module~$\Anab^+(n)$ is torsion-free. 
This can be seen directly by using that the $\Q$-vector space~$\Anab^+(n)$ 
admits an~$\N^n$-grading where the elements~$d_i$ act by isomorphisms 
of degree~$(0,\ldots,0,1,0,\ldots,0)$ where the~$1$ is at the~$i$-th place. 
Using the same grading and Corollary~\ref{c:Anabfree} 
one sees that for~$n>1$ the $\Lambda^+(n)$-module $\Anab^+(n)$ is not free.

We call a unitrivalent diagram~$D$ 
a {\em comb}, if it is a tree whose unique spanning
tree for the set of trivalent vertices is empty or a point or homeomorhic to an interval. 
We say that a unitrivalent diagram~$D$ is a {\em wheel} when it contains a 
circle whose complement is a disjoint union of intervals and each of these intervals has
exactly one labeled univalent vertex.
The first step in the proof of Theorem~\ref{t:Bnabgen} is the following lemma.

\begin{lemma}\label{l:Bgen2}
The $\Lambda(n)$-algebra~$\Bnab(n)$ is generated by~$G(n)$.
\end{lemma}
\begin{proof}
By relation~$(IntV)$ the algebra~$\Bnab(n)$ is generated by trees and wheels. By 
relation~$(AS)$ wheels of odd degree whose univalent vertices have the same label are
trivial in~$\Bnab(n)$. Relation~$(CL2B)$ then implies
that wheels of odd degree are trivial in~$\Bnab(n)$ and wheels of even
degree are related to~$\setbox1=\hbox{\input{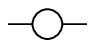}}\vcenter{\box1}\!$ by a monomial 
in~$d_i^{\pm 1}$ ($i=1,\ldots,n$).

Now we consider trees. 
By relation~$(IHX)$ it is sufficient to consider combs.
By relation~$(CL2B)$ we only have to
consider combs of degrees $1$, $2$, $3$, and $4$. All combs of
degree~$1$ appear in our list of generators, so there is nothing
to do. Now we treat combs of degree~$3$. In the computation below
we use relations~$(CL2B)$, $(IHX)$ and~$(AS)$, and~$(CL2B)$ and~$(AS)$ 
to write
such a comb as a linear combination of combs that have two
univalent vertices labeled by~$1$.

\begin{eqnarray}
d_1^2\H[i,j,k,\ell] & = & \setbox1=\hbox{\input{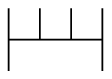}}\vcenter{\box1}\!\ = \
\setbox1=\hbox{\input{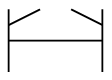}}\vcenter{\box1}\!-\setbox1=\hbox{\input{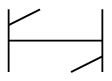}}\vcenter{\box1}\!-\setbox1=\hbox{\input{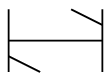}}\vcenter{\box1}\!+\setbox1=\hbox{\input{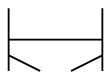}}\vcenter{\box1}\!\nonumber\\
& = &
d_kd_\ell\H[1,1,i,j]-d_jd_k\H[1,1,i,\ell]-d_id_\ell\H[1,1,k,j]+d_id_j\H[1,1,k,\ell]\label{e:combdeg3}
\end{eqnarray}

By relation~$(AS)$, the symmetry of the comb of degree~$3$, and by
invertibility of~$d_1$ we see that we can write a comb of
degree~$3$ as a linear combination of elements of~$G(n)$.

We continue with combs of degree~$4$. The proof of the first
equality below is similar to the computation in
equation~(\ref{e:combdeg3}), and the second equality follows from
relations~$(CL1B)$ and~$(AS)$.

\begin{eqnarray}
d_1^3d_m^{-1}\Ha[i,m,j,k,\ell] & = &
d_kd_\ell\Ha[1,1,1,i,j]-d_jd_k\Ha[1,1,1,i,\ell]-d_id_\ell\Ha[1,1,1,k,j]+d_id_j\Ha[1,1,1,k,\ell]\nonumber\\
& = &
d_1^2\left(d_kd_\ell\Yp[1,i,j]+d_jd_k\Yp[1,\ell,i]+d_id_\ell\Yp[1,j,k]+d_id_j\Yp[1,k,\ell]\
\right)\label{e:combdeg4}
\end{eqnarray}

By equation~(\ref{e:combdeg4}) and relation~$(AS)$ we can write
combs of degree~$4$ as linear combinations of elements of~$G(n)$. We 
apply equation~(\ref{e:combdeg4}) to combs of
degree~$2$ as follows~:

\begin{equation}
d_1\Yp[j,k,\ell]\ =\ d_1d_j^{-2}\Ha[j,j,j,k,\ell]\ =\
d_k\Yp[1,\ell,j]+d_\ell\Yp[1,j,k]+d_j\Yp[1,k,\ell].
\end{equation}

This completes the proof.
\end{proof}

In the proof of the following lemma we use 
the Lie superalgebra~$\gl(1\vert 1)$ that we will treat later in more detail. 

\begin{lemma}\label{l:Gnlinind2}
The elements of~$G(n)$ are $\Lambda(n)$-linearly
independent in~$\Bnab(n)$.
\end{lemma}
\begin{proof} The algebra $\Sym(\gl(1\vert 1)^{\oplus n})$ is generated by $H_i,
D_i, E_i, F_i$ ($i=1,\ldots, n$), where $X_i$ denotes the
element~$X$ of the $i$-th copy of~$\gl(1\vert 1)$ (see equation~(\ref{e:HDEF})). 
We
consider~$\Sym(\gl(1\vert 1)^{\oplus n})$ as a module
over~$\Q[d_1,\ldots, d_n]$ where~$d_i$ acts by multiplication
with~$D_i$. 
A well-known construction using the element~$\omega$ (see equation~(\ref{e:defomega}))
and equation~(\ref{e:FKVrel}) shows that there
exist morphisms of~$\Lambda(n)$-algebras

$$ u_n :
\Bnab(n)\longrightarrow\Lambda(n)\otimes_{\Q[d_1,\ldots,d_n]}\Sym(\gl(1\vert
1)^{\oplus n})^{\gl(1\vert 1)}
$$

satisfying

\begin{eqnarray}
u_n\left(\Ip[i,j]\right) & = & (1/2)(d_jH_i+d_iH_j) +F_iE_j-E_iF_j,\\
u_n\left(\!\Yp[1,i,j]\,\right) & = &
d_1(E_jF_i+F_jE_i)+d_i(E_1F_j+F_1E_j)+d_j(E_iF_1+F_iE_1),\label{e:unY}\\
u_n\left(\!\H[1,1,i,j]\,\right) & = &
d_1^2(F_iE_j+F_jE_i)+d_1d_i(E_1F_j+E_jF_1)+d_1d_j(E_iF_1+E_1F_i)\nonumber\\
& & +2d_id_jF_1E_1,\label{e:unH}\\
u_n\left(\setbox1=\hbox{\input{bub11b}}\vcenter{\box1}\!\right) & = & -2d_1^2.
\end{eqnarray}

The $\Lambda(n)$-algebra $\Bnab(n)$ has an $\N^3$-grading given
for unitrivalent diagrams~$D$ by $\partial(D)=(\partial_1(D),\partial_2(D),\partial_3(D))$,
where $\partial_1(D)$ is the number of connected components of~$D$ of
degree~$1$, $\partial_2(D)$ is the number of connected components in~$D$ of even degree, 
and $\partial_3(D)$ is the number of connected components in~$D$ of 
odd degree~$\geq 3$.
The
elements of~$G(n)$ are homogeneous with respect to~$\partial$. The
formulas for~$u_n(D)$ ($D\in G(n)$) from above 
imply that elements of~$G(n)$ of the
same degree are~$\Lambda(n)$-linearly 
independent.
\end{proof}

It follows from computations of~\cite{FKV} that the maps~$u_n$ from the proof
of Lemma~\ref{l:Gnlinind2} are compatible with relations~$(IntV)$ and satisfy

\begin{equation}\label{e:FKVrel}
u_n\left(\setbox1=\hbox{\input{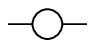}}\vcenter{\box1}\!\right)=  
-2d_id_j,\ \mbox{and}\ 
u_n\left(\setbox1=\hbox{\begin{picture}(0,0)%
\includegraphics{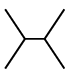}%
\end{picture}%
\setlength{\unitlength}{4144sp}%
\begingroup\makeatletter\ifx\SetFigFont\undefined%
\gdef\SetFigFont#1#2#3#4#5{%
  \reset@font\fontsize{#1}{#2pt}%
  \fontfamily{#3}\fontseries{#4}\fontshape{#5}%
  \selectfont}%
\fi\endgroup%
\begin{picture}(294,294)(394,197)
\end{picture}
}\vcenter{\box1}\!\right)=\frac{1}{2}
u_n\left(\setbox1=\hbox{\begin{picture}(0,0)%
\includegraphics{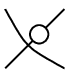}%
\end{picture}%
\setlength{\unitlength}{4144sp}%
\begingroup\makeatletter\ifx\SetFigFont\undefined%
\gdef\SetFigFont#1#2#3#4#5{%
  \reset@font\fontsize{#1}{#2pt}%
  \fontfamily{#3}\fontseries{#4}\fontshape{#5}%
  \selectfont}%
\fi\endgroup%
\begin{picture}(294,294)(664,-208)
\end{picture}
}\vcenter{\box1}\!+\setbox1=\hbox{\begin{picture}(0,0)%
\includegraphics{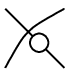}%
\end{picture}%
\setlength{\unitlength}{4144sp}%
\begingroup\makeatletter\ifx\SetFigFont\undefined%
\gdef\SetFigFont#1#2#3#4#5{%
  \reset@font\fontsize{#1}{#2pt}%
  \fontfamily{#3}\fontseries{#4}\fontshape{#5}%
  \selectfont}%
\fi\endgroup%
\begin{picture}(294,294)(664,-388)
\end{picture}
}\vcenter{\box1}\!-\setbox1=\hbox{\begin{picture}(0,0)%
\includegraphics{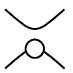}%
\end{picture}%
\setlength{\unitlength}{4144sp}%
\begingroup\makeatletter\ifx\SetFigFont\undefined%
\gdef\SetFigFont#1#2#3#4#5{%
  \reset@font\fontsize{#1}{#2pt}%
  \fontfamily{#3}\fontseries{#4}\fontshape{#5}%
  \selectfont}%
\fi\endgroup%
\begin{picture}(294,294)(664,197)
\end{picture}
}\vcenter{\box1}\!-\setbox1=\hbox{\begin{picture}(0,0)%
\includegraphics{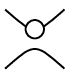}%
\end{picture}%
\setlength{\unitlength}{4144sp}%
\begingroup\makeatletter\ifx\SetFigFont\undefined%
\gdef\SetFigFont#1#2#3#4#5{%
  \reset@font\fontsize{#1}{#2pt}%
  \fontfamily{#3}\fontseries{#4}\fontshape{#5}%
  \selectfont}%
\fi\endgroup%
\begin{picture}(294,294)(664,17)
\end{picture}
}\vcenter{\box1}\!\right).
\end{equation}

Equation~(\ref{e:FKVrel}) or direct computations imply
that~$u_n$ is compatible with relations $(CL1B)$ and~$(CL2B)$. 
The maps~$u_n$ are not injective because
in contrast to Theorem~\ref{t:Bnabgen}  
the elements in equations~(\ref{e:unY}) and~(\ref{e:unH}) are nilpotent.

\smallskip

\begin{proof}[of Theorem~\ref{t:Bnabgen}]
Let $Q_n$ be the quotient field of~$\Lambda(n)$. Consider the
commutative $Q_n$-algebra $D_n=Q_n\otimes_{\Lambda(n)}\Bnab(n)$.
Let the degree of a unitrivalent diagram be the number of its components.
Then~$D_n$ admits the structure of a connected 
primitively generated graded Hopf algebra of finite type over~$Q_n$ whose
space of primitive elements~$P_n$ is the homogeneous part of degree~$1$ of~$D_n$ and
whose counit is the augmentation map. By~\cite{MiM} the algebra~$Q_n$ is freely
generated by~$P_n$. By Lemmas~\ref{l:Bgen2} and \ref{l:Gnlinind2} the 
map~$\Bnab(n)\longrightarrow D_n$ induced by~$\Lambda(n)\subset Q_n$ is 
injective and maps~$G(n)$ to a basis of~$P_n$. This completes the proof.
\end{proof}

\section{Drinfeld associators in $\Anabc_\coeff^m(3)$}\label{s:assoC}

Let $\Anab^1(n)$ be the quotient of $\Anab^0(n)$ by the ideal~$I_1$
generated by trivalent diagrams~$D$ with~$b_1(D\setminus\Gamma_n)>0$. 
Define $\widehat{\Lambda}_\coeff^+(\Gamma)$ (resp.\ 
$\Anabc_\coeff^0(\Gamma)$, $\widehat{\Lambda}_\coeff^+(n)$, $\Anabc_\coeff^+(n)$, 
$\Anabc_\coeff^0(n)$, $\Anabc_\coeff^1(n)$) 
by extending coefficients of
$\Lambda^+(\Gamma)$ (resp.\ $\Anab^0(\Gamma)$, 
$\Lambda^+(n)$, $\Anab^+(n)$, 
$\Anab^0(n)$, $\Anab^1(n)$)
from~$\Q$ to the commutative $\Q$-algebra~$\coeff$ followed by
completion.
Define continous 
morphisms of $\coeff$-algebras

\begin{eqnarray*} &
\widetilde{\Delta}_i:\widehat{\Lambda}_\coeff^+(n)
\longrightarrow\widehat{\Lambda}_\coeff^+(n+1) \quad
(i=0,\ldots,n+1),\\
& \widetilde{\Delta}_i(d_j)=d_j\ \mbox{for $j<i$},\quad 
\widetilde{\Delta}_i(d_j)=d_{j+1}\ \mbox{for $j>i$,}\quad
\widetilde{\Delta}_i(d_i)=d_i+d_{i+1},\quad \mbox{and} &\\
& \widetilde{\epsilon}_i:
\widehat{\Lambda}_\coeff^+(n)\longrightarrow\widehat{\Lambda}_\coeff^+(n-1) 
\quad (i=1,\ldots,n),&\\
& \widetilde{\epsilon}_i(d_j)=d_j\ \mbox{for $j<i$},\quad 
\widetilde{\epsilon}_i(d_j)=d_{j-1}\ \mbox{for $j>i$,}\quad
\widetilde{\epsilon}_i(d_i)=0. &
\end{eqnarray*}

Then the maps~$\widetilde{\Delta_i}\otimes\Delta_i:
\widehat{\Lambda}_\coeff^+(n)\otimes\Ac(n)\longrightarrow\widehat{\Lambda}_\coeff^+(n+1)
\otimes\Ac_\coeff(n+1)$ 
and $\widetilde{\epsilon}_i
\otimes\epsilon_i:
\widehat{\Lambda}_\coeff^+(n)\otimes\Ac(n)\longrightarrow
\widehat{\Lambda}_\coeff^+(n-1)\otimes\Ac(n-1)$
induce continous linear maps

\begin{equation}\label{e:Anabcoalg}
\Delta_i:\Anabc_\coeff^m(n)\longrightarrow\Anabc_\coeff^m(n+1) \quad \mbox{and} \quad 
\epsilon_i:\Anabc_\coeff^m(n)
\longrightarrow\Anabc_\coeff^m(n-1) \quad (m='+',0,1).
\end{equation}

The proof that $\epsilon_i$ in equation~(\ref{e:Anabcoalg}) 
is well-defined uses an $\N^n$-grading of~$\Anab^m(n)$.
A direct proof that~$\Delta_i$ is well-defined 
uses that~$\Anab^+(n)$ is a 
torsion-free $\Lambda^+(n)$-module and requires a
computation to ensure compatibility with relation~$(CL1A)$. 

\begin{defi}
A Drinfeld associator~$\Phi$ in~$\Anabc^+_\coeff(3)$ 
(resp.\ $\Anabc^m_\coeff(3)$ with $m=0,1$)
is a solution of 
equations (DA1)-(DA4) in $\Anabc^+_\coeff(n)$ (resp.\ $\Anabc^m_\coeff(n)$).
\end{defi}

Since the canonical map~$\tau_m$ 
from~$\Ac_\coeff(n)$ to~$\Anabc_\coeff^m(n)$ commutes with~$\Delta_i$ and~$\epsilon_i$
there exist even horizontal group-like 
Drinfeld associators in~$\Anabc_\coeff^m(3)$ by Fact~\ref{f:existPhi}. 
We use the notion of a twist~$T\in\Anabc_\coeff^m(2)$ ($m='+',0,1$) as defined
in Section~\ref{s:assoA}. For example, 

\begin{equation}
\exp(d_1^2 d_2^2)\quad \mbox{and}\quad \exp\left(d_1d_2(d_1+d_2)\setbox1=\hbox{\begin{picture}(0,0)%
\includegraphics{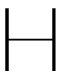}%
\end{picture}%
\setlength{\unitlength}{3947sp}%
\begingroup\makeatletter\ifx\SetFigFont\undefined%
\gdef\SetFigFont#1#2#3#4#5{%
  \reset@font\fontsize{#1}{#2pt}%
  \fontfamily{#3}\fontseries{#4}\fontshape{#5}%
  \selectfont}%
\fi\endgroup%
\begin{picture}(269,344)(579,-158)
\end{picture}
}\vcenter{\box1}\!\right) 
\end{equation}

are 
non-trivial twists in~$\Anabc_\coeff^+(2)$.
With the definition of even, horizontal, and group-like 
Drinfeld associators in
$\Anabc_\coeff^m(3)$ ($m=0,1$) as in Section~\ref{s:assoA}, we have the following
theorem.

\begin{theorem}\label{t:existF2}
(1a) There exists exactly one even horizontal group-like Drinfeld associator in
$\Anabc_\coeff^0(3)$.

(1b) There exists exactly one even group-like Drinfeld associator in $\Anabc_\coeff^1(3)$.

(2) The unique Drinfeld associator in (1a) and (1b) is equal to 

$$\exp\left(F\cdot [t^{12}, t^{23}]\right)\in\Anabc^m(3)\subset\Anabc^m_\coeff(3)
\quad (m=0,1)$$

for a unique $F\in\Q[[C, D, E]]\subset\widehat{\Lambda}(3)$ where the inclusion
is given by

$$ C=(d_1+d_2)^2+(d_2+d_3)^2, D=(d_1+d_2)^2(d_2+d_3)^2,
E=d_1d_3-d_2(d_1+d_2+d_3).$$
\end{theorem}
\begin{proof} {\em Existence of F :} Let~$\Phi$ be an even horizontal group-like Drinfeld
associator in~$\Ac(3)$. By~$(DA5)$ we have~$\Phi=\exp(p(t^{12}, t^{23}))$ for a Lie
series $p$ in~$A$ and $B$ that involves only terms of even degrees. The free Lie
algebra with two generators~$A$ and~$B$ is spanned linearly by a set~$K$ that is defined
recursively by~$A, B\in K$ and by~$\ad_A(c), \ad_B(c)\in K$ for $c\in K$. Therefore
Lemma~\ref{l:adjoinABC} below implies that 

$$\Phi_m=\tau_m(\Phi)=\exp(F\cdot [t^{12}, t^{23}])$$ 

for some $F\in\Q[[(d_1+d_2)^2, (d_2+d_3)^2, E]]$. 
Equation~(DA3) implies that $p$ satisfies
$p(t^{23}, t^{12})=p(t^{12}, t^{23})\in\Anabc(3)$. 
As a consequence, $F$ is invariant
under the permutation of $(d_1+d_2)^2$ and $(d_2+d_3)^2$. This implies $F\in\Q[[C,D,E]]$.


{\em Uniqueness of $\Phi_1$ :}
The proofs of Theorems~8 and~9 of~\cite{LeM}
can be adapted
to see that any even group-like
associator in~$\Anabc_\coeff^1(3)$ is related to~$\Phi_1$ by an even group-like twist
$T\in\Anabc_\coeff^1(2)$ (see also Lemma~4.17 of~\cite{BN2}). 
Corollary~\ref{c:Anabfree} implies $T=1$. Therefore~$\Phi_1$
is unique.


{\em Uniqueness of $F$ :} Corollary~\ref{c:Anabfree} implies that the $\Lambda^+(n)$-module
$$M=\Anab^+(n)/\Lambda^+(n)I_1\supset\Anab^1(n)$$ is torsion-free and~$0\not=[t^{12}, t^{23}]\in M$.
This implies that~$F\in\Lambda^+(3)$ is uniquely determined by~$\Phi_1$.
It is easy to see that $C,D,E$ are algebraically independent by using that 
$d_1+d_2, d_2+d_3, E$ are
algebraically independent. 
Therefore~$F$ is uniquely determined as a formal series in $C,D,E$.


{\em Uniqueness of~$\Phi_0$ :} Except for the uniqueness of the Lie series $p$ equation~$(DA5)$
holds also for horizontal group-like Drinfeld associators in~$\Anabc_\coeff^0(3)$. As in
the first part of the proof we see that an even horizontal group-like Drinfeld associator
$\Psi$ in~$\Anabc_\coeff^0(3)$ can be expressed as~$\Psi=\exp(F'\cdot [t^{12}, t^{23}])$
for some $F'\in\coeff[[C,D,E]]$. Let $\pi:\Anabc_\coeff^0(n)\longrightarrow\Anabc_\coeff^1(n)$
be the canonical projection. Then 

$$\pi(\Psi)=\exp(F'\cdot [t^{12}, t^{23}])=\Phi_1=
\exp(F\cdot [t^{12}, t^{23}])$$ 

by the uniqueness of~$\Phi_1$ which implies $F=F'$ by
the uniqueness of~$F$. This completes the proof.
\end{proof}

It follows from Theorem~\ref{t:existF2} that any even group-like associator 
in~$\Ac_\coeff(3)$ is mapped to an associator with rational coefficients 
in~$\Anabc_\coeff^1(3)$ by the canonical projection.
The next lemma was used in the proof of Theorem~\ref{t:existF2}.

\begin{lemma}\label{l:adjoinABC}
The following identities hold in $\Anab^0(3)$~:

\begin{eqnarray}
\mbox{$[t^{12},[t^{12},[t^{12},t^{23}]]]$} = (d_1+d_2)^2\,
[t^{12}, t^{23}],\ 
\mbox{$[t^{23},[t^{23},[t^{12},t^{23}]]]$} = (d_2+d_3)^2\,
[t^{12}, t^{23}], & & \label{e:ad12}
\\ 
\mbox{$[t^{12},[t^{23},[t^{12},t^{23}]]]$}
 = (d_1d_3-d_2(d_1+d_2+d_3))\, [t^{12}, t^{23}] =
\mbox{$[t^{23},[t^{12},[t^{12},t^{23}]]]$}.\ & & \mbox{}\label{e:ad3}
\end{eqnarray}
\end{lemma}
\begin{proof}
In the computation in~$\Anab^0(\Gamma)$ 
below, the first equality follows from
relations~$(LS)$ and~$(STU)$.
For the second equality we apply
relations~$(LS)$, $(AS)$, and $(CL2A)$.
The third equality follows from relation~$(CL1A)$.

\begin{eqnarray}
\setbox1=\hbox{\input{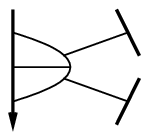}}\vcenter{\box1}\! & = &
\frac{1}{2}\left(\setbox1=\hbox{\input{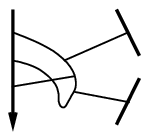}}\vcenter{\box1}\!+\setbox1=\hbox{\input{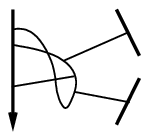}}\vcenter{\box1}\!+\setbox1=\hbox{\input{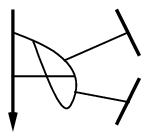}}\vcenter{\box1}\!\right)\label{e:CL3A}\\
& = &
\frac{1}{2}\left(\setbox1=\hbox{\input{y5}}\vcenter{\box1}\!-\setbox1=\hbox{\input{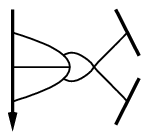}}\vcenter{\box1}\!-d_i^2\setbox1=\hbox{\input{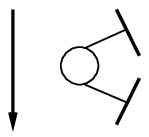}}\vcenter{\box1}\!\right)
= d_i^2\setbox1=\hbox{\input{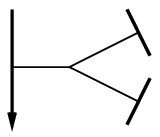}}\vcenter{\box1}\!-\frac{d_i^2}{2}\setbox1=\hbox{\input{emptybub}}\vcenter{\box1}\!\nonumber
\end{eqnarray}

Now we prove
equation~(\ref{e:hsymm}) by the computation in~$\Anab^0(3)$ below. The three equalities 
in this computation follow by applying
relations~$(CL2A)$ and~$(AS)$, equation~(\ref{e:CL3A}),
relations~$(STU)$, $(AS)$, and $(CL2A)$, respectively.

\begin{eqnarray}
\setbox1=\hbox{\begin{picture}(0,0)%
\includegraphics{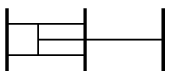}%
\end{picture}%
\setlength{\unitlength}{3947sp}%
\begingroup\makeatletter\ifx\SetFigFont\undefined%
\gdef\SetFigFont#1#2#3#4#5{%
  \reset@font\fontsize{#1}{#2pt}%
  \fontfamily{#3}\fontseries{#4}\fontshape{#5}%
  \selectfont}%
\fi\endgroup%
\begin{picture}(794,344)(579,-158)
\end{picture}
}\vcenter{\box1}\! & = &
-d_1^{-1}d_3\setbox1=\hbox{\begin{picture}(0,0)%
\includegraphics{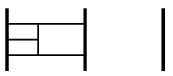}%
\end{picture}%
\setlength{\unitlength}{3947sp}%
\begingroup\makeatletter\ifx\SetFigFont\undefined%
\gdef\SetFigFont#1#2#3#4#5{%
  \reset@font\fontsize{#1}{#2pt}%
  \fontfamily{#3}\fontseries{#4}\fontshape{#5}%
  \selectfont}%
\fi\endgroup%
\begin{picture}(794,344)(579,-158)
\end{picture}
}\vcenter{\box1}\!=
-d_1^{-1}d_3\left(d_1^2\setbox1=\hbox{\begin{picture}(0,0)%
\includegraphics{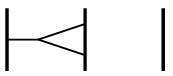}%
\end{picture}%
\setlength{\unitlength}{3947sp}%
\begingroup\makeatletter\ifx\SetFigFont\undefined%
\gdef\SetFigFont#1#2#3#4#5{%
  \reset@font\fontsize{#1}{#2pt}%
  \fontfamily{#3}\fontseries{#4}\fontshape{#5}%
  \selectfont}%
\fi\endgroup%
\begin{picture}(794,344)(579,-158)
\end{picture}
}\vcenter{\box1}\!-(1/2)d_1^2\setbox1=\hbox{\begin{picture}(0,0)%
\includegraphics{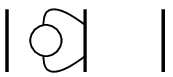}%
\end{picture}%
\setlength{\unitlength}{3947sp}%
\begingroup\makeatletter\ifx\SetFigFont\undefined%
\gdef\SetFigFont#1#2#3#4#5{%
  \reset@font\fontsize{#1}{#2pt}%
  \fontfamily{#3}\fontseries{#4}\fontshape{#5}%
  \selectfont}%
\fi\endgroup%
\begin{picture}(794,344)(579,-158)
\end{picture}
}\vcenter{\box1}\!\right)\nonumber\\
& = &
(1/2) (d_1+d_2)d_3\setbox1=\hbox{\begin{picture}(0,0)%
\includegraphics{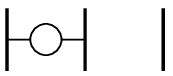}%
\end{picture}%
\setlength{\unitlength}{3947sp}%
\begingroup\makeatletter\ifx\SetFigFont\undefined%
\gdef\SetFigFont#1#2#3#4#5{%
  \reset@font\fontsize{#1}{#2pt}%
  \fontfamily{#3}\fontseries{#4}\fontshape{#5}%
  \selectfont}%
\fi\endgroup%
\begin{picture}(794,344)(579,-158)
\end{picture}
}\vcenter{\box1}\!\label{e:hsymm}
\end{eqnarray}

In the computation below, we apply relations $(STU)$ and $(AS)$ to obtain 
the first three equalities. The fourth equality follows from
equations~(\ref{e:CL3A}) and (\ref{e:hsymm}) and from
relation~$(CL2A)$ and~$(AS)$. The fifth equality is implied by equation~(\ref{e:CL3A}).

\begin{eqnarray*}
& & [t^{12},[t^{12},[t^{12},t^{23}]]]\\ & = &
\left[t^{12},\left[t^{12},\setbox1=\hbox{\begin{picture}(0,0)%
\includegraphics{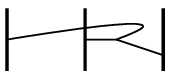}%
\end{picture}%
\setlength{\unitlength}{3947sp}%
\begingroup\makeatletter\ifx\SetFigFont\undefined%
\gdef\SetFigFont#1#2#3#4#5{%
  \reset@font\fontsize{#1}{#2pt}%
  \fontfamily{#3}\fontseries{#4}\fontshape{#5}%
  \selectfont}%
\fi\endgroup%
\begin{picture}(794,344)(579,-158)
\end{picture}
}\vcenter{\box1}\!\right]\right]\\ & = &
\left[t^{12},-\setbox1=\hbox{\begin{picture}(0,0)%
\includegraphics{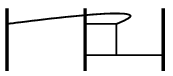}%
\end{picture}%
\setlength{\unitlength}{3947sp}%
\begingroup\makeatletter\ifx\SetFigFont\undefined%
\gdef\SetFigFont#1#2#3#4#5{%
  \reset@font\fontsize{#1}{#2pt}%
  \fontfamily{#3}\fontseries{#4}\fontshape{#5}%
  \selectfont}%
\fi\endgroup%
\begin{picture}(794,344)(579,-158)
\end{picture}
}\vcenter{\box1}\!+\setbox1=\hbox{\begin{picture}(0,0)%
\includegraphics{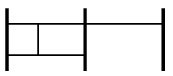}%
\end{picture}%
\setlength{\unitlength}{3947sp}%
\begingroup\makeatletter\ifx\SetFigFont\undefined%
\gdef\SetFigFont#1#2#3#4#5{%
  \reset@font\fontsize{#1}{#2pt}%
  \fontfamily{#3}\fontseries{#4}\fontshape{#5}%
  \selectfont}%
\fi\endgroup%
\begin{picture}(794,344)(579,-158)
\end{picture}
}\vcenter{\box1}\!\right]\\ & = &
\setbox1=\hbox{\begin{picture}(0,0)%
\includegraphics{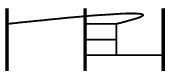}%
\end{picture}%
\setlength{\unitlength}{3947sp}%
\begingroup\makeatletter\ifx\SetFigFont\undefined%
\gdef\SetFigFont#1#2#3#4#5{%
  \reset@font\fontsize{#1}{#2pt}%
  \fontfamily{#3}\fontseries{#4}\fontshape{#5}%
  \selectfont}%
\fi\endgroup%
\begin{picture}(794,344)(579,-158)
\end{picture}
}\vcenter{\box1}\!+\setbox1=\hbox{\begin{picture}(0,0)%
\includegraphics{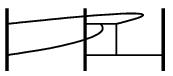}%
\end{picture}%
\setlength{\unitlength}{3947sp}%
\begingroup\makeatletter\ifx\SetFigFont\undefined%
\gdef\SetFigFont#1#2#3#4#5{%
  \reset@font\fontsize{#1}{#2pt}%
  \fontfamily{#3}\fontseries{#4}\fontshape{#5}%
  \selectfont}%
\fi\endgroup%
\begin{picture}(794,344)(579,-158)
\end{picture}
}\vcenter{\box1}\!+2\setbox1=\hbox{\begin{picture}(0,0)%
\includegraphics{hsymm.eps}%
\end{picture}%
\setlength{\unitlength}{3947sp}%
\begingroup\makeatletter\ifx\SetFigFont\undefined%
\gdef\SetFigFont#1#2#3#4#5{%
  \reset@font\fontsize{#1}{#2pt}%
  \fontfamily{#3}\fontseries{#4}\fontshape{#5}%
  \selectfont}%
\fi\endgroup%
\begin{picture}(794,344)(579,-158)
\end{picture}
}\vcenter{\box1}\!+\setbox1=\hbox{\begin{picture}(0,0)%
\includegraphics{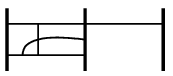}%
\end{picture}%
\setlength{\unitlength}{3947sp}%
\begingroup\makeatletter\ifx\SetFigFont\undefined%
\gdef\SetFigFont#1#2#3#4#5{%
  \reset@font\fontsize{#1}{#2pt}%
  \fontfamily{#3}\fontseries{#4}\fontshape{#5}%
  \selectfont}%
\fi\endgroup%
\begin{picture}(794,344)(579,-158)
\end{picture}
}\vcenter{\box1}\!+\setbox1=\hbox{\begin{picture}(0,0)%
\includegraphics{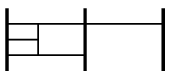}%
\end{picture}%
\setlength{\unitlength}{3947sp}%
\begingroup\makeatletter\ifx\SetFigFont\undefined%
\gdef\SetFigFont#1#2#3#4#5{%
  \reset@font\fontsize{#1}{#2pt}%
  \fontfamily{#3}\fontseries{#4}\fontshape{#5}%
  \selectfont}%
\fi\endgroup%
\begin{picture}(794,344)(579,-158)
\end{picture}
}\vcenter{\box1}\!\\
& = &
d_2^2\setbox1=\hbox{\begin{picture}(0,0)%
\includegraphics{comm.eps}%
\end{picture}%
\setlength{\unitlength}{3947sp}%
\begingroup\makeatletter\ifx\SetFigFont\undefined%
\gdef\SetFigFont#1#2#3#4#5{%
  \reset@font\fontsize{#1}{#2pt}%
  \fontfamily{#3}\fontseries{#4}\fontshape{#5}%
  \selectfont}%
\fi\endgroup%
\begin{picture}(794,344)(579,-158)
\end{picture}
}\vcenter{\box1}\!-(1/2)d_2d_3\setbox1=\hbox{\begin{picture}(0,0)%
\includegraphics{bub12b.eps}%
\end{picture}%
\setlength{\unitlength}{3947sp}%
\begingroup\makeatletter\ifx\SetFigFont\undefined%
\gdef\SetFigFont#1#2#3#4#5{%
  \reset@font\fontsize{#1}{#2pt}%
  \fontfamily{#3}\fontseries{#4}\fontshape{#5}%
  \selectfont}%
\fi\endgroup%
\begin{picture}(794,344)(579,-158)
\end{picture}
}\vcenter{\box1}\!+d_1d_2^{-1}\setbox1=\hbox{\begin{picture}(0,0)%
\includegraphics{cl1right.eps}%
\end{picture}%
\setlength{\unitlength}{3947sp}%
\begingroup\makeatletter\ifx\SetFigFont\undefined%
\gdef\SetFigFont#1#2#3#4#5{%
  \reset@font\fontsize{#1}{#2pt}%
  \fontfamily{#3}\fontseries{#4}\fontshape{#5}%
  \selectfont}%
\fi\endgroup%
\begin{picture}(794,344)(579,-158)
\end{picture}
}\vcenter{\box1}\!
+(d_1+d_2)d_3\setbox1=\hbox{\begin{picture}(0,0)%
\includegraphics{bub12b.eps}%
\end{picture}%
\setlength{\unitlength}{3947sp}%
\begingroup\makeatletter\ifx\SetFigFont\undefined%
\gdef\SetFigFont#1#2#3#4#5{%
  \reset@font\fontsize{#1}{#2pt}%
  \fontfamily{#3}\fontseries{#4}\fontshape{#5}%
  \selectfont}%
\fi\endgroup%
\begin{picture}(794,344)(579,-158)
\end{picture}
}\vcenter{\box1}\!\\
& &
+d_1^{-1}d_2\setbox1=\hbox{\begin{picture}(0,0)%
\includegraphics{cl1left.eps}%
\end{picture}%
\setlength{\unitlength}{3947sp}%
\begingroup\makeatletter\ifx\SetFigFont\undefined%
\gdef\SetFigFont#1#2#3#4#5{%
  \reset@font\fontsize{#1}{#2pt}%
  \fontfamily{#3}\fontseries{#4}\fontshape{#5}%
  \selectfont}%
\fi\endgroup%
\begin{picture}(794,344)(579,-158)
\end{picture}
}\vcenter{\box1}\!+d_1^2\setbox1=\hbox{\begin{picture}(0,0)%
\includegraphics{comm.eps}%
\end{picture}%
\setlength{\unitlength}{3947sp}%
\begingroup\makeatletter\ifx\SetFigFont\undefined%
\gdef\SetFigFont#1#2#3#4#5{%
  \reset@font\fontsize{#1}{#2pt}%
  \fontfamily{#3}\fontseries{#4}\fontshape{#5}%
  \selectfont}%
\fi\endgroup%
\begin{picture}(794,344)(579,-158)
\end{picture}
}\vcenter{\box1}\!-(1/2)d_1d_3\setbox1=\hbox{\begin{picture}(0,0)%
\includegraphics{bub12b.eps}%
\end{picture}%
\setlength{\unitlength}{3947sp}%
\begingroup\makeatletter\ifx\SetFigFont\undefined%
\gdef\SetFigFont#1#2#3#4#5{%
  \reset@font\fontsize{#1}{#2pt}%
  \fontfamily{#3}\fontseries{#4}\fontshape{#5}%
  \selectfont}%
\fi\endgroup%
\begin{picture}(794,344)(579,-158)
\end{picture}
}\vcenter{\box1}\!\\ &
= &
(d_1^2+d_2^2)[t^{12},t^{23}]+(1/2)(d_1+d_2)d_3\setbox1=\hbox{\begin{picture}(0,0)%
\includegraphics{bub12b.eps}%
\end{picture}%
\setlength{\unitlength}{3947sp}%
\begingroup\makeatletter\ifx\SetFigFont\undefined%
\gdef\SetFigFont#1#2#3#4#5{%
  \reset@font\fontsize{#1}{#2pt}%
  \fontfamily{#3}\fontseries{#4}\fontshape{#5}%
  \selectfont}%
\fi\endgroup%
\begin{picture}(794,344)(579,-158)
\end{picture}
}\vcenter{\box1}\!
+d_1d_2^{-1}\left(d_2^2\setbox1=\hbox{\begin{picture}(0,0)%
\includegraphics{comm.eps}%
\end{picture}%
\setlength{\unitlength}{3947sp}%
\begingroup\makeatletter\ifx\SetFigFont\undefined%
\gdef\SetFigFont#1#2#3#4#5{%
  \reset@font\fontsize{#1}{#2pt}%
  \fontfamily{#3}\fontseries{#4}\fontshape{#5}%
  \selectfont}%
\fi\endgroup%
\begin{picture}(794,344)(579,-158)
\end{picture}
}\vcenter{\box1}\!\right.\\
& & \left. -(1/2)d_2d_3\setbox1=\hbox{\begin{picture}(0,0)%
\includegraphics{bub12b.eps}%
\end{picture}%
\setlength{\unitlength}{3947sp}%
\begingroup\makeatletter\ifx\SetFigFont\undefined%
\gdef\SetFigFont#1#2#3#4#5{%
  \reset@font\fontsize{#1}{#2pt}%
  \fontfamily{#3}\fontseries{#4}\fontshape{#5}%
  \selectfont}%
\fi\endgroup%
\begin{picture}(794,344)(579,-158)
\end{picture}
}\vcenter{\box1}\!\right)+d_1^{-1}d_2
\left(d_1^2\setbox1=\hbox{\begin{picture}(0,0)%
\includegraphics{comm.eps}%
\end{picture}%
\setlength{\unitlength}{3947sp}%
\begingroup\makeatletter\ifx\SetFigFont\undefined%
\gdef\SetFigFont#1#2#3#4#5{%
  \reset@font\fontsize{#1}{#2pt}%
  \fontfamily{#3}\fontseries{#4}\fontshape{#5}%
  \selectfont}%
\fi\endgroup%
\begin{picture}(794,344)(579,-158)
\end{picture}
}\vcenter{\box1}\!-(1/2)d_1d_3\setbox1=\hbox{\begin{picture}(0,0)%
\includegraphics{bub12b.eps}%
\end{picture}%
\setlength{\unitlength}{3947sp}%
\begingroup\makeatletter\ifx\SetFigFont\undefined%
\gdef\SetFigFont#1#2#3#4#5{%
  \reset@font\fontsize{#1}{#2pt}%
  \fontfamily{#3}\fontseries{#4}\fontshape{#5}%
  \selectfont}%
\fi\endgroup%
\begin{picture}(794,344)(579,-158)
\end{picture}
}\vcenter{\box1}\!\right)\\ & =
& (d_1+d_2)^2[t^{12}, t^{23}]
\end{eqnarray*}

This proves the first part of equation~(\ref{e:ad12}). 
The second part of equation~(\ref{e:ad12}) follows
by applying the $\Q$-algebra automorphism
of~$\Anab^0(3)$ induced by interchanging the intervals~$1$ and~$3$ of the skeleton. The
proof of equation~(\ref{e:ad3}) is similar to the proof of
equation~(\ref{e:ad12}).
\end{proof}

Lemma~\ref{l:adjoinABC} and Theorem~\ref{t:existF2} imply 
that the denominator of the homogeneous part of degree~$2n$ 
in~$\exp (F\cdot x)\in\Q[[C,D,E,x]]$ ($\deg(x)=2$) is a
lower bound for this denominator in any even group-like associator in~$\Ac(3)$.
Now we come to P.\ Vogel's proof of the formula for~$F$ using Theorem~\ref{t:F}.

\smallskip

\begin{proof}[of equations~(\ref{e:fF}) to~(\ref{e:defX})]
There exists a unique solution~$\Psi\in\Q[[u,v]]$ of

\begin{equation}\label{e:eqG}
1+u=\cosh\left(\Psi\sqrt{uv}\right)+\sqrt{u/v}\sinh\left(\Psi\sqrt{uv}\right)=1+u\Psi+(1/2)uv\Psi^2+\ldots.
\end{equation}

By Theorem~\ref{t:F} the series~$\Psi$ is related to~$F$ by equations~(\ref{e:fF}) and~(\ref{e:defX}).
Equation~(\ref{e:eqG}) implies the following equation in~$R=\Q[[u,v]][(uv)^{-1/2}]$~:

$$\left(1+\sqrt{u/v}\right)e^{2\Psi\sqrt{uv}}-2(1+u)e^{\Psi\sqrt{uv}}+1-\sqrt{u/v}=0.$$

The solutions in~$R[\sqrt{v},(\sqrt{u}+\sqrt{v})^{-1}]$ of this quadratic equation are
$$e^{\Psi\sqrt{uv}}
={1+u\pm\sqrt{(1+u)^2-1+u/v}\over 1+\sqrt{u/v}}
={\sqrt{v}(1+u)\pm\sqrt{u}\sqrt{1+v(u+2)}\over\sqrt{u}+\sqrt{v}}$$

where we can determine the sign~$\pm$ because~$\Psi\in\Q[[u,v]]$.
Let $H\in R[[w]]$ be the unique solution of

$$e^{H\sqrt{uv}}={\sqrt{v}\sqrt{1+uw}+\sqrt{u}\sqrt{1+vw}\over\sqrt{u}+\sqrt{v}}
\in R[[w]].$$

Then $\Psi=H(u+2)$. 
%
Let $H'$ be the partial derivative of $H$ with respect to $w$. Then
$$H'\sqrt{uv}e^{H\sqrt{uv}}={1\over 2\left(\sqrt{u}+\sqrt{v}\right)}
\left({\sqrt{v}u\over\sqrt{1+uw}}+{\sqrt{u}v\over\sqrt{1+vw}}\right)$$
$$={\sqrt{uv}\over 2\left(\sqrt{u}+\sqrt{v}\right)}{\sqrt{v}\sqrt{1+uw}+\sqrt{u}\sqrt{1+vw}\over
\sqrt{(1+uw)(1+vw)}}.$$
Therefore
$$H'={1\over 2\sqrt{(1+uw)(1+vw)}}=
\sum_{p,q\geq0}\pmatrix{-1/2\cr p\cr}\pmatrix{-1/2\cr q\cr}{u^pv^qw^{p+q}\over2}$$
which implies
$$H=\sum_{p,q\geq0}\pmatrix{-1/2\cr p\cr}\pmatrix{-1/2\cr q\cr}{u^pv^qw^{p+q+1}\over2(p+q+1)}$$
because $H$ is~$0$ for $w=0$.
This implies equation~(\ref{e:defPsi}) by substituting~$w=u+2$.
\end{proof}

By solving equation~(\ref{e:eqG}) iteratively one sees that the series~$\Psi$ in 
equation~(\ref{e:defPsi}) satisfies~$\Psi\in\Q[u][[v]]\subset\Q[[u,v]]$ and
that the coefficient of~$v^q$ of~$\Psi$ is a polynomial in~$u$ of degree~$q$.

It would be interesting if the associator of Theorem~\ref{t:existF2} could be used
to investigate the coefficients of the image
of~$\Phi_{KZ}\in\Ac_\C(3)$ in~$\Anabc_\C^m(3)$ ($m=0,1$)
by the canonical projection. 
The simplest relation holds for~$m=1$~:

\begin{remark}\label{r:existG}
For every group-like Drinfeld associator $\Phi\in\Anabc_\coeff^1(3)$
there exists a unique series $G(\Phi)\in\coeff[[d_1+d_2, d_1d_2]]\subset\Lambda^+_\coeff(2)$ 
such that $\Phi$ is related to the associator 
$\exp(F\cdot [t^{12}, t^{23}])\in\Anabc_\coeff^1(3)$ of Theorem~\ref{t:existF2}
by the twist~$\exp\left(G(\Phi)\cdot \setbox1=\hbox{\begin{picture}(0,0)%
\includegraphics{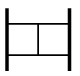}%
\end{picture}%
\setlength{\unitlength}{3947sp}%
\begingroup\makeatletter\ifx\SetFigFont\undefined%
\gdef\SetFigFont#1#2#3#4#5{%
  \reset@font\fontsize{#1}{#2pt}%
  \fontfamily{#3}\fontseries{#4}\fontshape{#5}%
  \selectfont}%
\fi\endgroup%
\begin{picture}(344,344)(579,-158)
\end{picture}
}\vcenter{\box1}\!\right)$.
\end{remark}

The existence of the twist in Remark~\ref{r:existG}
follows as in the proof of Theorem~\ref{t:existF2}
by using the structure of~$\Anab^1_\coeff(2)$. The uniqueness of the twist
uses the structure of~$\Anab^1_\coeff(1)$ and the triviality of~$H^2(K^\bullet)$ where the
cochain complex
$K^\bullet$ consists of the 
spaces~$K^n=\bigcap_{i=1}^n\Ker(\epsilon_i)\subset \Anab^1_\coeff(n)$ and the
coboundary maps $\delta^n=\sum_{i=0}^{n+1}(-1)^i\Delta_i$ (compare Proposition~4.4 and
Remark~4.9 of~\cite{BN2}).

\section{Tensor products and duality of $\gloo$-modules}\label{s:rep1}

Let $V=V_0\oplus V_1$ be a $\Z/(2)$-graded vector space. The
algebra~$\End(V)$ has a natural~$\Z/(2)$-grading with homogeneous
components~$\End(V)_0=\End(V_0)\oplus \End(V_1)$
and~$\End(V)_1=\Hom(V_0,V_1)\oplus \Hom(V_1,V_0)$. The
map~$\str:\End(V)\longrightarrow\Q$ defined
by
$$\str(\varphi)=\tr(\varphi_{\vert
\End(V_0)})-\tr(\varphi_{\vert \End(V_1)})$$ 
is called supertrace.
The number~$\sdim(V)=\dim(V_0)-\dim(V_1)=\str(\id_V)$ is called
the superdimension of~$V$. We consider~$\End(V)$ as a Lie
superalgebra~$\gl(V)$ with bracket~$[\cdot,\cdot]$ induced
by

$$[A,B]=A B-(-1)^{\deg(A)\deg (B)}B A$$

for homogeneous elements $A,B\in\End(V)=\gl(V)$. For~$\dim
V_0=\dim V_1=1$ the Lie superalgebra~$\gl(V)$ is isomorphic to
$\gloo$ which has a homogeneous basis given by the following four
matrices~:

\begin{equation}\label{e:HDEF}
 H=\left(\begin{array}{ll} 1 & 0\\ 0 & -1\end{array}\right),
\quad D=\left(\begin{array}{ll} 1 & 0\\ 0 & 1\end{array}\right),
\quad E=\left(\begin{array}{ll} 0 & 1\\ 0 & 0\end{array}\right),
\quad F=\left(\begin{array}{ll} 0 & 0\\ 1 & 0\end{array}\right).
\end{equation}

The non-vanishing brackets of basis elements are given by

$$ [H,E]=-[E,H]=2E,\quad [F,H]=-[H,F]=2F,\quad [E,F]=[F,E]=D. $$

For $t=(\lambda,\mu,\sigma)\in\Q^*\times\Q\times \Z/(2)$ there
exists a unique $2$-dimensional $\gloo$-module $\V$ with
$\sdim(\V)=0$ such that for all $v\in \V$ with $\deg(v)=\sigma$
we have $H\cdot v=(\mu+1) v$, $D\cdot v=\lambda v$ and $E\cdot
v=0$.

For each triple $t=(\lambda,\mu,\sigma)\in\Q^*\times\Q\times
\Z/(2)$ we fix a choice of a vector $0\not=\vv\in \V$ with
$\deg(\vv)=\sigma$. Denote $F\cdot \vv$ by $\ww$. Then $H\cdot
\ww=(\mu-1)\ww$, $D\cdot \ww=\lambda \ww$, $F\cdot \ww=0$
and~$E\cdot \ww=\lambda \vv$. In particular, the vectors~$\vv$
and~$\ww$ form a basis of~$\V$. It is easy to see that the
modules~$\V$ are simple.

For triples $t_i=(\lambda_i, \mu_i,
\sigma_i)\in\Q^*\times\Q\times \Z/(2)$ ($i=1,2$) with
$\lambda_1+\lambda_2\not=0$ and for $e_2=(0,1,0), e_3=(0,0,1)
\in\Q^*\times\Q\times \Z/(2)$ we define $\gloo$-linear maps

\begin{eqnarray}
\Ystd_{t_1, t_2} : V_{t_1+t_2+e_2}\longrightarrow V_{t_1}\otimes
V_{t_2} & \mbox{by} & \Ystd_{t_1,
t_2}(v_{t_1+t_2+e_2})=\vvo\otimes \vvp\mbox{, and}\\ \Ydot_{t_1,
t_2} : V_{t_1+t_2-e_2+e_3}\longrightarrow \Vo\otimes \Vp &
\mbox{by} & \Ydot_{t_1, t_2}(w_{t_1+t_2-e_2+e_3})=\wwo\otimes
\wwp.
\end{eqnarray}

Since~$\Ystd_{t_1, t_2}$ and~$\Ydot_{t_1, t_2}$ are non-trivial,
well-defined, and have non-isomorphic simple images,
equation~(\ref{e:tensordecomp}) holds for reasons of dimension.

\begin{equation}\label{e:tensordecomp} \Vo\otimes
\Vp\cong V_{t_1+t_2+e_2}\oplus V_{t_1+t_2-e_2+e_3} \quad\mbox{if
$\lambda_1+\lambda_2\not=0$}.
\end{equation}

The $\Z/(2)$-graded space $(\V)^*=\Hom_\Q(\V,\Q)$ becomes a~$\gloo$-module
by

$$ (a\cdot \beta)(v)=(-1)^{\deg(a)\deg(\beta)}\beta(-a\cdot v) $$

for all $v\in\V$ and homogeneous elements $a\in\gloo,
\beta\in(\V)^*$. Let~$\alpha\in (\V)^*$ be given
by~$\alpha(\ww)=1$, $\alpha(\vv)=0$. We have
$H\cdot\alpha=(-\mu+1)\alpha$, $D\cdot\alpha=-\lambda\alpha$,
$E\cdot\alpha=0$, and~$\deg(\alpha)=\sigma+1$. This implies

\begin{equation}\label{e:dualmod}
\left(\V\right)^*\cong V_{t^*}\quad\mbox{where}\quad
t^*=(-\lambda,-\mu,\sigma+1).
\end{equation}

For $(\mu, \sigma)\in\Q\times \Z/(2)$ there exists a unique
$1$-dimensional representation~$I_\mu^\sigma=\Q$ of~$\gloo$ with
$\sdim(I_\mu^\sigma)=(-1)^\sigma$ and~$H\cdot v=\mu v$ for
all~$v\in I_\mu^\sigma$. The formulas
$$D\cdot v=E\cdot v=F\cdot
v=0$$ 
hold for all $v\in I_\mu^\sigma$. 
The modules $I_\mu^\sigma$ and~$V_t$ form a complete set of isomorphism
types of simple $\gl(1\vert 1)$-modules up to isomorphisms of degree~$0$.
Define a $\gloo$-linear
map $\gcap_t:V_{t^*}\otimes V_t\longrightarrow I_0^0=\Q$ by

\begin{equation}
\gcap_t(v_{t^*}\otimes \ww)=1\quad\mbox{and}\quad
\gcap_t(w_{t^*}\otimes \vv)=(-1)^\sigma.
\end{equation}

The map~$\gcap_t$ can be obtained as the composition
of~$\iota\otimes\id$ with the evaluation map, where~$\iota$ is the
isomorphism that we used to prove equation~(\ref{e:dualmod}). We
define~$\gcup_t:I_0^0\longrightarrow V_{t^*}\otimes V_t$
by~$\gcup_t(1)=w_{t^*}\otimes\vv-(-1)^\sigma v_{t^*}\otimes\ww$.
Then we have

\begin{equation}\label{e:templieb}
(\gcap_{t^*}\otimes \I_t
)\circ(\I_t\otimes\gcup_t)=\I_t=(\I_t\otimes\gcap_t)\circ(\gcup_{t^*}\otimes
\I_t),
\end{equation}

where~$\I_t=\id_{V_t}$. There exist unique $\gloo$-linear
maps~$\Astd_{t_1,t_2}$, $\Adot_{t_1,t_2}$ satisfying

\begin{equation}\label{e:defA}
\Astd_{t_1,t_2}\circ
\Ystd_{t_1,t_2}=\I_{t_1+t_2+e_2}\quad\mbox{and}\quad
\Adot_{t_1,t_2}\circ\Ydot_{t_1,t_2}=\I_{t_1+t_2-e_2+e_3}.
\end{equation}

We have~$\Astd_{t_1,t_2}\circ\Ydot_{t_1,t_2}=0$
and~$\Adot_{t_1,t_2}\circ\Ystd_{t_1,t_2}=0$ because these maps are
homomorphisms between non-isomorphic simple modules. 
By equation~(\ref{e:tensordecomp}) and Schurs lemma (see \cite{Kac}) the elements

\begin{equation}\label{e:end0bas}
\Ystd_{t_1,t_2}\circ \Astd_{t_1,t_2}\quad ,\quad \Ydot_{t_1, t_2}\circ\Adot_{t_1, t_2}
\end{equation}

are a basis of the vector space~$\End_{\gloo}(V_{t_1}\otimes V_{t_2})$ of $\gloo$-linear
endomorphisms of \mbox{$V_{t_1}\otimes V_{t_2}$} of degree~$0$.
In the
following lemma we present some relations between the morphisms
introduced in this section.

\begin{lemma}\label{l:AY}
For $t_i=(\lambda_i,\mu_i,\sigma_i)\in\Q^*\times\Q\times\Z/(2)$
($i=1,\ldots,4$) with~$t_3=t_1+t_2+e_2$ and~$t_4=t_1+t_2-e_2+e_3$ we have

\begin{eqnarray*}
\Adot_{t_1, t_2} & = & (\I_{t_4}\otimes
\gcap_{t_2})\circ(\Ystd_{t_4, t_2^*}\otimes
\I_{t_2})=(\gcap_{t_1^*}\otimes
\I_{t_4})\circ(\I_{t_1}\otimes\Ystd_{t_1^*,t_4}),\\
\Astd_{t_1,t_2} & = &
(-1)^{\sigma_2}\frac{\lambda_1}{\lambda_3}(\I_{t_3}\otimes\gcap_{t_2})\circ(\Ydot_{t_3,
t_2^*} \otimes \I_{t_2})\\ & = &
(-1)^{\sigma_1}\frac{\lambda_2}{\lambda_3}(\gcap_{t_1^*}\otimes
\I_{t_3})\circ(\I_{t_1}\otimes\Ydot_{t_1^*,t_3}).
\end{eqnarray*}
\end{lemma}
\begin{proof}
We compare the three elements of~$\Hom_\gloo(V_{t_1}\otimes
V_{t_2}, V_{t_4})$ from the first equation by comparing the images
of~$v_{t_1}\otimes v_{t_2}$ and~$w_{t_1}\otimes w_{t_2}$ which
determine the morphisms. We obtain

\begin{eqnarray*}
& & (\I_{t_4}\otimes \gcap_{t_2})\circ(\Ystd_{t_4, t_2^*}\otimes
\I_{t_2})(v_{t_1}\otimes v_{t_2})=(\I_{t_4}\otimes
\gcap_{t_2})(v_{t_4}\otimes v_{t_2^*}\otimes v_{t_2})=0, \\ & &
\Adot_{t_1, t_2}(v_{t_1}\otimes v_{t_2})=0,\\
 & & (\gcap_{t_1^*}\otimes
\I_{t_4})\circ(\I_{t_1}\otimes\Ystd_{t_1^*,t_4})(v_{t_1}\otimes
v_{t_2})=(\gcap_{t_1^*}\otimes \I_{t_4})(v_{t_1}\otimes
v_{t_1^*}\otimes v_{t_4})=0
\end{eqnarray*}

as the value of~$v_{t_1}\otimes v_{t_2}$, and

\begin{eqnarray*}
 & & (\I_{t_4}\otimes
\gcap_{t_2})\circ(\Ystd_{t_4, t_2^*}\otimes
\I_{t_2})(w_{t_1}\otimes w_{t_2})=(\I_{t_4}\otimes
\gcap_{t_2})(F\cdot (v_{t_4}\otimes v_{t_2^*})\otimes w_{t_2})\\ &
& =(\I_{t_4}\otimes \gcap_{t_2})(w_{t_4}\otimes v_{t_2^*}\otimes
w_{t_2}+(-1)^{\sigma_4}v_{t_4}\otimes w_{t_2^*}\otimes
w_{t_2})=w_{t_4},\\
& & \Adot_{t_1, t_2}(w_{t_1}\otimes w_{t_2})=w_{t_4},
\\
 & & (\gcap_{t_1^*}\otimes
\I_{t_4})\circ(\I_{t_1}\otimes\Ystd_{t_1^*,t_4})(w_{t_1}\otimes
w_{t_2})=(\gcap_{t_1^*}\otimes \I_{t_4})(w_{t_1}\otimes
F\cdot(v_{t_1^*}\otimes v_{t_4}))\\ & & =(\gcap_{t_1^*}\otimes
\I_{t_4})(w_{t_1}\otimes w_{t_1^*}\otimes
v_{t_4}-(-1)^{\sigma_1}w_{t_1}\otimes v_{t_1^*}\otimes
w_{t_4})=w_{t_4}.
\end{eqnarray*}

 We illustrate a small difference in the proof
of the second equation by a sample computation.

\begin{eqnarray*}
(\I_{t_3}\otimes\gcap_{t_2})\circ(\Ydot_{t_3, t_2^*} \otimes
\I_{t_2})(v_{t_1}\otimes
v_{t_2})=\frac{1}{\lambda_1}(\I_{t_3}\otimes\gcap_{t_2})(E\cdot(w_{t_3}\otimes
w_{t_2^*})\otimes v_{t_2})\\
=(\I_{t_3}\otimes\gcap_{t_2})\left(\frac{\lambda_3}{\lambda_1}v_{t_3}\otimes
w_{t_2^*}\otimes
v_{t_2}+(-1)^{\sigma_3}\frac{\lambda_2}{\lambda_1}w_{t_3}\otimes
v_{t_2^*}\otimes
v_{t_2}\right)=(-1)^{\sigma_2}\frac{\lambda_3}{\lambda_1}
v_{t_3}.
\end{eqnarray*}

The rest of the proof is straightforward.
\end{proof}

\section{The tensor functor~$\Wnab_0$}\label{s:W0}

Let~$\Gamma$ be a unitrivalent graph with oriented edges and
vertices whose edges $i\in E(\Gamma)$ are colored by pairs~$(\lambda_i,
\mu_i)\in\Q^*\times \Q$.  
We say that the coloring of~$\Gamma$
is {\em admissible}, if the following two conditions hold at every
trivalent vertex~$v$ of~$\Gamma$ that is incident to $i,j,k\in E(\Gamma)$:

\begin{eqnarray}
s_{v,i}\lambda_i+s_{v,j}\lambda_j+s_{v,k}\lambda_k & = & 0,\label{e:admissibilitya}\\
s_{v,i}\mu_i+s_{v,j}\mu_j+s_{v,k}\mu_k & = & s_{v,i}s_{v,j}s_{v,k}.\label{e:admissibilityb}
\end{eqnarray}

When $\Gamma$ is admissibly colored then~$\Lambda(\Gamma)\not=0$ by
equation~(\ref{e:admissibilitya}) because~$\lambda_i\not=0$ for
all~$i\in E(\Gamma)$.
Let us define a category~$\Anabo$. Objects of~$\Anabo$ are finite
sequences $(c_1,\ldots, c_k)$ of triples
$$c_i=(\lambda_i,\mu_i,\sigma_i)\in \Q^*\times\Q\times\Z/(2).$$
Morphisms from~$c^0=(c_1^0,\ldots,c_{k_0}^0)$ to
$c^1=(c_1^1,\ldots,c_{k_1}^1)$ with $c_j^i=(\lambda_j^i, \mu_j^i,
\sigma_j^i)$ consist of unitrivalent graphs~$\Gamma$ with
admissible coloring, where the univalent vertices of~$\Gamma$ are
related to~$c^0$ and~$c^1$ as follows~: the boundary of~$\Gamma$
is decomposed into two disjoint
subsets~$\partial\Gamma=\partial_0\Gamma\cup\partial_1\Gamma$,
$\partial_i\Gamma$ is in bijection with~$\{1,\ldots,k_i\}$, the
$j$-th vertex of~$\partial_i\Gamma$ is incident to an edge
of~$\Gamma$ that is colored by~$(\lambda_j^i,\mu_j^i)$, and that
edge is directed towards the boundary iff~$\sigma_j^i+i\equiv
0\mod 2$.
We represent morphisms of~$\Anabo$
graphically by generic pictures of~$\Gamma$ in~$\R\times
[0,1]$ where the $i$-th boundary point of~$\partial_j\Gamma$ has
coordinates~$(i,j)$. The composition $b\circ
a\in\Hom_{\Anabo}(c^0,c^2)$ of morphisms $a\in\Hom_{\Anabo}(c^0,
c^1), b\in\Hom_{\Anabo}(c^1, c^2)$ is defined graphically by
placing~$b$ onto the top of~$a$ and by shrinking the result 
to~$\R\times [0,1]$. We define a tensor product of objects
of~$\Anabo$ by concatenation of sequences. For
$a\in\Hom_{\Anabo}(c^0, c^1), b\in\Hom_{\Anabo}(c^2, c^3)$ we
define $a\otimes b\in\Hom_{\Anabo}(c^0\otimes c^2, c^1\otimes
c^3)$ graphically by placing~$a$ to the left of~$b$.

To an object $c=(c_1, \ldots, c_k)$ of~$\Anabo$ with
$c_i=(\lambda_i,\mu_i, \sigma_i)$ we assign the~$\gloo$-module

\begin{equation}\label{e:Wot}
W_0(c)=V_{t_1}\otimes\ldots\otimes V_{t_k}\quad\mbox{with}\quad
t_i=((-1)^{\sigma_i}\lambda_i,(-1)^{\sigma_i}\mu_i,\sigma_i).
\end{equation}

Let $t=(\lambda, \mu, 0)$. We do not specify colors or
orientations of edges in pictures, when they are uniquely
determined by the context or when the formulas hold for arbitrary orientations. 
Assume that the
graphs in equations~(\ref{e:id}) to~(\ref{e:cup}) are colored by
$(\lambda, \mu)$. Then we define

\begin{eqnarray}
W_0\left(\setbox1=\hbox{\begin{picture}(0,0)%
\includegraphics{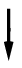}%
\end{picture}%
\setlength{\unitlength}{3947sp}%
\begingroup\makeatletter\ifx\SetFigFont\undefined%
\gdef\SetFigFont#1#2#3#4#5{%
  \reset@font\fontsize{#1}{#2pt}%
  \fontfamily{#3}\fontseries{#4}\fontshape{#5}%
  \selectfont}%
\fi\endgroup%
\begin{picture}(74,294)(339,217)
\end{picture}
}\vcenter{\box1}\!\right)=\I_t, & &
W_0\left(\setbox1=\hbox{\begin{picture}(0,0)%
\includegraphics{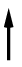}%
\end{picture}%
\setlength{\unitlength}{3947sp}%
\begingroup\makeatletter\ifx\SetFigFont\undefined%
\gdef\SetFigFont#1#2#3#4#5{%
  \reset@font\fontsize{#1}{#2pt}%
  \fontfamily{#3}\fontseries{#4}\fontshape{#5}%
  \selectfont}%
\fi\endgroup%
\begin{picture}(74,294)(339,267)
\end{picture}
}\vcenter{\box1}\!\right)=\I_{t^*},\label{e:id}\\
W_0\left(\setbox1=\hbox{\begin{picture}(0,0)%
\includegraphics{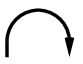}%
\end{picture}%
\setlength{\unitlength}{3947sp}%
\begingroup\makeatletter\ifx\SetFigFont\undefined%
\gdef\SetFigFont#1#2#3#4#5{%
  \reset@font\fontsize{#1}{#2pt}%
  \fontfamily{#3}\fontseries{#4}\fontshape{#5}%
  \selectfont}%
\fi\endgroup%
\begin{picture}(359,263)(579,217)
\end{picture}
}\vcenter{\box1}\!\right)=\gcap_t, & &
W_0\left(\setbox1=\hbox{\begin{picture}(0,0)%
\includegraphics{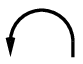}%
\end{picture}%
\setlength{\unitlength}{3947sp}%
\begingroup\makeatletter\ifx\SetFigFont\undefined%
\gdef\SetFigFont#1#2#3#4#5{%
  \reset@font\fontsize{#1}{#2pt}%
  \fontfamily{#3}\fontseries{#4}\fontshape{#5}%
  \selectfont}%
\fi\endgroup%
\begin{picture}(359,263)(564,217)
\end{picture}
}\vcenter{\box1}\!\right)=\gcap_{t^*},\\
W_0\left(\setbox1=\hbox{\begin{picture}(0,0)%
\includegraphics{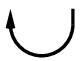}%
\end{picture}%
\setlength{\unitlength}{3947sp}%
\begingroup\makeatletter\ifx\SetFigFont\undefined%
\gdef\SetFigFont#1#2#3#4#5{%
  \reset@font\fontsize{#1}{#2pt}%
  \fontfamily{#3}\fontseries{#4}\fontshape{#5}%
  \selectfont}%
\fi\endgroup%
\begin{picture}(359,263)(564,298)
\end{picture}
}\vcenter{\box1}\!\right)=\gcup_t, & &
W_0\left(\setbox1=\hbox{\begin{picture}(0,0)%
\includegraphics{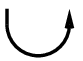}%
\end{picture}%
\setlength{\unitlength}{3947sp}%
\begingroup\makeatletter\ifx\SetFigFont\undefined%
\gdef\SetFigFont#1#2#3#4#5{%
  \reset@font\fontsize{#1}{#2pt}%
  \fontfamily{#3}\fontseries{#4}\fontshape{#5}%
  \selectfont}%
\fi\endgroup%
\begin{picture}(359,263)(579,298)
\end{picture}
}\vcenter{\box1}\!\right)=\gcup_{t^*}.\label{e:cup}
\end{eqnarray}

Let $c_i=(\lambda_i,\mu_i,\sigma_i)$. Consider orientations and
colors of the graphs~$\setbox1=\hbox{\begin{picture}(0,0)%
\includegraphics{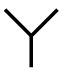}%
\end{picture}%
\setlength{\unitlength}{3947sp}%
\begingroup\makeatletter\ifx\SetFigFont\undefined%
\gdef\SetFigFont#1#2#3#4#5{%
  \reset@font\fontsize{#1}{#2pt}%
  \fontfamily{#3}\fontseries{#4}\fontshape{#5}%
  \selectfont}%
\fi\endgroup%
\begin{picture}(294,319)(304,217)
\end{picture}
}\vcenter{\box1}\!$ and~$\setbox1=\hbox{\begin{picture}(0,0)%
\includegraphics{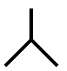}%
\end{picture}%
\setlength{\unitlength}{3947sp}%
\begingroup\makeatletter\ifx\SetFigFont\undefined%
\gdef\SetFigFont#1#2#3#4#5{%
  \reset@font\fontsize{#1}{#2pt}%
  \fontfamily{#3}\fontseries{#4}\fontshape{#5}%
  \selectfont}%
\fi\endgroup%
\begin{picture}(294,319)(304,-58)
\end{picture}
}\vcenter{\box1}\!$ such that

\begin{equation}\setbox1=\hbox{\begin{picture}(0,0)%
\includegraphics{y.eps}%
\end{picture}%
\setlength{\unitlength}{3947sp}%
\begingroup\makeatletter\ifx\SetFigFont\undefined%
\gdef\SetFigFont#1#2#3#4#5{%
  \reset@font\fontsize{#1}{#2pt}%
  \fontfamily{#3}\fontseries{#4}\fontshape{#5}%
  \selectfont}%
\fi\endgroup%
\begin{picture}(294,319)(304,217)
\end{picture}
}\vcenter{\box1}\!\in\Hom_{\Anabo}(c_1, (c_2, c_3))\quad\mbox{and}\quad
\setbox1=\hbox{\begin{picture}(0,0)%
\includegraphics{a.eps}%
\end{picture}%
\setlength{\unitlength}{3947sp}%
\begingroup\makeatletter\ifx\SetFigFont\undefined%
\gdef\SetFigFont#1#2#3#4#5{%
  \reset@font\fontsize{#1}{#2pt}%
  \fontfamily{#3}\fontseries{#4}\fontshape{#5}%
  \selectfont}%
\fi\endgroup%
\begin{picture}(294,319)(304,-58)
\end{picture}
}\vcenter{\box1}\!\in\Hom_{\Anabo}((c_1, c_2), c_3).
\end{equation}

Then we define
$t_i=((-1)^{\sigma_i}\lambda_i,(-1)^{\sigma_i}\mu_i,\sigma_i)$
and

\begin{eqnarray}
W_0\left(\setbox1=\hbox{\begin{picture}(0,0)%
\includegraphics{y.eps}%
\end{picture}%
\setlength{\unitlength}{3947sp}%
\begingroup\makeatletter\ifx\SetFigFont\undefined%
\gdef\SetFigFont#1#2#3#4#5{%
  \reset@font\fontsize{#1}{#2pt}%
  \fontfamily{#3}\fontseries{#4}\fontshape{#5}%
  \selectfont}%
\fi\endgroup%
\begin{picture}(294,319)(304,217)
\end{picture}
}\vcenter{\box1}\!\right) & = &
\left\{\begin{array}{ll}\Ystd_{t_2,t_3} & \mbox{if
$\sigma_1+\sigma_2+\sigma_3\equiv 0\mod 2$},\label{e:WoYo}\\
\lambda_1\Ydot_{t_2,t_3} & \mbox{if
$\sigma_1+\sigma_2+\sigma_3\equiv 1\mod
2$},\end{array}\right. \\ W_0\left(\setbox1=\hbox{\begin{picture}(0,0)%
\includegraphics{a.eps}%
\end{picture}%
\setlength{\unitlength}{3947sp}%
\begingroup\makeatletter\ifx\SetFigFont\undefined%
\gdef\SetFigFont#1#2#3#4#5{%
  \reset@font\fontsize{#1}{#2pt}%
  \fontfamily{#3}\fontseries{#4}\fontshape{#5}%
  \selectfont}%
\fi\endgroup%
\begin{picture}(294,319)(304,-58)
\end{picture}
}\vcenter{\box1}\!\right) & = &
\left\{\begin{array}{ll}\lambda_3\Astd_{t_1,t_2} & \mbox{if
$\sigma_1+\sigma_2+\sigma_3\equiv 0\mod 2$},\label{e:WoAo}\\
\Adot_{t_1,t_2} & \mbox{if $\sigma_1+\sigma_2+\sigma_3\equiv
1\mod 2$}.
\end{array}\right.
\end{eqnarray}

Let $\Xt_{t_1,t_2}\in\Hom_{\gloo}(W_0((c_1, c_2)),W_0((c_2, c_1)))$ be the
superpermutation of tensor factors induced
by~$\Xt_{t_1,t_2}(v\otimes w)=(-1)^{\deg(v)\deg(w)} w\otimes v$.
We define

\begin{equation}\label{e:WoX}
W_0\left(\setbox1=\hbox{\begin{picture}(0,0)%
\includegraphics{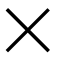}%
\end{picture}%
\setlength{\unitlength}{3947sp}%
\begingroup\makeatletter\ifx\SetFigFont\undefined%
\gdef\SetFigFont#1#2#3#4#5{%
  \reset@font\fontsize{#1}{#2pt}%
  \fontfamily{#3}\fontseries{#4}\fontshape{#5}%
  \selectfont}%
\fi\endgroup%
\begin{picture}(244,244)(329,267)
\end{picture}
}\vcenter{\box1}\!\right)=\Xt_{t_1,t_2}\quad\mbox{for}\quad
\setbox1=\hbox{\begin{picture}(0,0)%
\includegraphics{x.eps}%
\end{picture}%
\setlength{\unitlength}{3947sp}%
\begingroup\makeatletter\ifx\SetFigFont\undefined%
\gdef\SetFigFont#1#2#3#4#5{%
  \reset@font\fontsize{#1}{#2pt}%
  \fontfamily{#3}\fontseries{#4}\fontshape{#5}%
  \selectfont}%
\fi\endgroup%
\begin{picture}(244,244)(329,267)
\end{picture}
}\vcenter{\box1}\!\in\Hom_{\Anabo}((c_1, c_2), (c_2, c_1)).
\end{equation}

We have the following lemma.

\begin{lemma}\label{l:W0} The map~$W_0$ induces a tensor functor from~$\Anabo$ to $\gloo$-modules.
\end{lemma}
\begin{proof}
The tensor category~$\Anabo$ has a presentation by admissibly colored
generators\footnote{Some of these generators are superfluous, but
facilitate the statement of the defining relations of~$\Anabo$.}

$$ \setbox1=\hbox{\begin{picture}(0,0)%
\includegraphics{down.eps}%
\end{picture}%
\setlength{\unitlength}{3947sp}%
\begingroup\makeatletter\ifx\SetFigFont\undefined%
\gdef\SetFigFont#1#2#3#4#5{%
  \reset@font\fontsize{#1}{#2pt}%
  \fontfamily{#3}\fontseries{#4}\fontshape{#5}%
  \selectfont}%
\fi\endgroup%
\begin{picture}(74,294)(339,217)
\end{picture}
}\vcenter{\box1}\!, \setbox1=\hbox{\begin{picture}(0,0)%
\includegraphics{up.eps}%
\end{picture}%
\setlength{\unitlength}{3947sp}%
\begingroup\makeatletter\ifx\SetFigFont\undefined%
\gdef\SetFigFont#1#2#3#4#5{%
  \reset@font\fontsize{#1}{#2pt}%
  \fontfamily{#3}\fontseries{#4}\fontshape{#5}%
  \selectfont}%
\fi\endgroup%
\begin{picture}(74,294)(339,267)
\end{picture}
}\vcenter{\box1}\!, \setbox1=\hbox{\begin{picture}(0,0)%
\includegraphics{rightcap.eps}%
\end{picture}%
\setlength{\unitlength}{3947sp}%
\begingroup\makeatletter\ifx\SetFigFont\undefined%
\gdef\SetFigFont#1#2#3#4#5{%
  \reset@font\fontsize{#1}{#2pt}%
  \fontfamily{#3}\fontseries{#4}\fontshape{#5}%
  \selectfont}%
\fi\endgroup%
\begin{picture}(359,263)(579,217)
\end{picture}
}\vcenter{\box1}\!,
\setbox1=\hbox{\begin{picture}(0,0)%
\includegraphics{rightcup.eps}%
\end{picture}%
\setlength{\unitlength}{3947sp}%
\begingroup\makeatletter\ifx\SetFigFont\undefined%
\gdef\SetFigFont#1#2#3#4#5{%
  \reset@font\fontsize{#1}{#2pt}%
  \fontfamily{#3}\fontseries{#4}\fontshape{#5}%
  \selectfont}%
\fi\endgroup%
\begin{picture}(359,263)(579,298)
\end{picture}
}\vcenter{\box1}\!, \setbox1=\hbox{\begin{picture}(0,0)%
\includegraphics{leftcap.eps}%
\end{picture}%
\setlength{\unitlength}{3947sp}%
\begingroup\makeatletter\ifx\SetFigFont\undefined%
\gdef\SetFigFont#1#2#3#4#5{%
  \reset@font\fontsize{#1}{#2pt}%
  \fontfamily{#3}\fontseries{#4}\fontshape{#5}%
  \selectfont}%
\fi\endgroup%
\begin{picture}(359,263)(564,217)
\end{picture}
}\vcenter{\box1}\!, \setbox1=\hbox{\begin{picture}(0,0)%
\includegraphics{leftcup.eps}%
\end{picture}%
\setlength{\unitlength}{3947sp}%
\begingroup\makeatletter\ifx\SetFigFont\undefined%
\gdef\SetFigFont#1#2#3#4#5{%
  \reset@font\fontsize{#1}{#2pt}%
  \fontfamily{#3}\fontseries{#4}\fontshape{#5}%
  \selectfont}%
\fi\endgroup%
\begin{picture}(359,263)(564,298)
\end{picture}
}\vcenter{\box1}\!,
\setbox1=\hbox{\begin{picture}(0,0)%
\includegraphics{y.eps}%
\end{picture}%
\setlength{\unitlength}{3947sp}%
\begingroup\makeatletter\ifx\SetFigFont\undefined%
\gdef\SetFigFont#1#2#3#4#5{%
  \reset@font\fontsize{#1}{#2pt}%
  \fontfamily{#3}\fontseries{#4}\fontshape{#5}%
  \selectfont}%
\fi\endgroup%
\begin{picture}(294,319)(304,217)
\end{picture}
}\vcenter{\box1}\!, \setbox1=\hbox{\begin{picture}(0,0)%
\includegraphics{a.eps}%
\end{picture}%
\setlength{\unitlength}{3947sp}%
\begingroup\makeatletter\ifx\SetFigFont\undefined%
\gdef\SetFigFont#1#2#3#4#5{%
  \reset@font\fontsize{#1}{#2pt}%
  \fontfamily{#3}\fontseries{#4}\fontshape{#5}%
  \selectfont}%
\fi\endgroup%
\begin{picture}(294,319)(304,-58)
\end{picture}
}\vcenter{\box1}\!, \setbox1=\hbox{\begin{picture}(0,0)%
\includegraphics{x.eps}%
\end{picture}%
\setlength{\unitlength}{3947sp}%
\begingroup\makeatletter\ifx\SetFigFont\undefined%
\gdef\SetFigFont#1#2#3#4#5{%
  \reset@font\fontsize{#1}{#2pt}%
  \fontfamily{#3}\fontseries{#4}\fontshape{#5}%
  \selectfont}%
\fi\endgroup%
\begin{picture}(244,244)(329,267)
\end{picture}
}\vcenter{\box1}\! $$

(where for the last three graphs all orientations are possible) modulo the relations

\begin{eqnarray}
& \setbox1=\hbox{\input{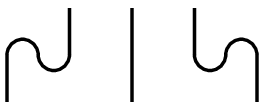}}\vcenter{\box1}\!,\quad
\setbox1=\hbox{\begin{picture}(0,0)%
\includegraphics{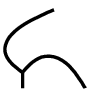}%
\end{picture}%
\setlength{\unitlength}{3947sp}%
\begingroup\makeatletter\ifx\SetFigFont\undefined%
\gdef\SetFigFont#1#2#3#4#5{%
  \reset@font\fontsize{#1}{#2pt}%
  \fontfamily{#3}\fontseries{#4}\fontshape{#5}%
  \selectfont}%
\fi\endgroup%
\begin{picture}(429,419)(119,367)
\end{picture}
}\vcenter{\box1}\!\!\!\!=\!\!\!\setbox1=\hbox{\begin{picture}(0,0)%
\includegraphics{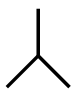}%
\end{picture}%
\setlength{\unitlength}{3947sp}%
\begingroup\makeatletter\ifx\SetFigFont\undefined%
\gdef\SetFigFont#1#2#3#4#5{%
  \reset@font\fontsize{#1}{#2pt}%
  \fontfamily{#3}\fontseries{#4}\fontshape{#5}%
  \selectfont}%
\fi\endgroup%
\begin{picture}(344,419)(729,-683)
\end{picture}
}\vcenter{\box1}\!\!\!\!=\!\!\!\setbox1=\hbox{\begin{picture}(0,0)%
\includegraphics{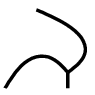}%
\end{picture}%
\setlength{\unitlength}{3947sp}%
\begingroup\makeatletter\ifx\SetFigFont\undefined%
\gdef\SetFigFont#1#2#3#4#5{%
  \reset@font\fontsize{#1}{#2pt}%
  \fontfamily{#3}\fontseries{#4}\fontshape{#5}%
  \selectfont}%
\fi\endgroup%
\begin{picture}(429,419)(1254,-533)
\end{picture}
}\vcenter{\box1}\!,
& \label{e:relverify}\\ & \setbox1=\hbox{\begin{picture}(0,0)%
\includegraphics{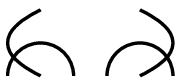}%
\end{picture}%
\setlength{\unitlength}{3947sp}%
\begingroup\makeatletter\ifx\SetFigFont\undefined%
\gdef\SetFigFont#1#2#3#4#5{%
  \reset@font\fontsize{#1}{#2pt}%
  \fontfamily{#3}\fontseries{#4}\fontshape{#5}%
  \selectfont}%
\fi\endgroup%
\begin{picture}(838,361)(332,-175)
\put(691,-47){\makebox(0,0)[lb]{\smash{\SetFigFont{12}{14.4}{\familydefault}{\mddefault}{\updefault}{\color[rgb]{0,0,0}$=$}%
}}}
\end{picture}
}\vcenter{\box1}\!,\ \setbox1=\hbox{\begin{picture}(0,0)%
\includegraphics{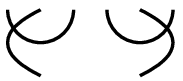}%
\end{picture}%
\setlength{\unitlength}{3947sp}%
\begingroup\makeatletter\ifx\SetFigFont\undefined%
\gdef\SetFigFont#1#2#3#4#5{%
  \reset@font\fontsize{#1}{#2pt}%
  \fontfamily{#3}\fontseries{#4}\fontshape{#5}%
  \selectfont}%
\fi\endgroup%
\begin{picture}(838,361)(332,-608)
\put(676,-511){\makebox(0,0)[lb]{\smash{\SetFigFont{12}{14.4}{\familydefault}{\mddefault}{\updefault}{\color[rgb]{0,0,0}$=$}%
}}}
\end{picture}
}\vcenter{\box1}\!,\quad
\setbox1=\hbox{\begin{picture}(0,0)%
\includegraphics{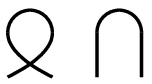}%
\end{picture}%
\setlength{\unitlength}{3947sp}%
\begingroup\makeatletter\ifx\SetFigFont\undefined%
\gdef\SetFigFont#1#2#3#4#5{%
  \reset@font\fontsize{#1}{#2pt}%
  \fontfamily{#3}\fontseries{#4}\fontshape{#5}%
  \selectfont}%
\fi\endgroup%
\begin{picture}(686,359)(183,-83)
\put(458, 46){\makebox(0,0)[lb]{\smash{\SetFigFont{12}{14.4}{\familydefault}{\mddefault}{\updefault}{\color[rgb]{0,0,0}$=$}%
}}}
\end{picture}
}\vcenter{\box1}\!,\ \setbox1=\hbox{\begin{picture}(0,0)%
\includegraphics{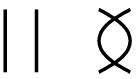}%
\end{picture}%
\setlength{\unitlength}{3947sp}%
\begingroup\makeatletter\ifx\SetFigFont\undefined%
\gdef\SetFigFont#1#2#3#4#5{%
  \reset@font\fontsize{#1}{#2pt}%
  \fontfamily{#3}\fontseries{#4}\fontshape{#5}%
  \selectfont}%
\fi\endgroup%
\begin{picture}(644,344)(354,217)
\put(611,314){\makebox(0,0)[lb]{\smash{\SetFigFont{12}{14.4}{\familydefault}{\mddefault}{\updefault}{\color[rgb]{0,0,0}$=$}%
}}}
\end{picture}
}\vcenter{\box1}\!,\ \setbox1=\hbox{\begin{picture}(0,0)%
\includegraphics{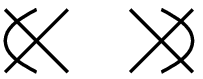}%
\end{picture}%
\setlength{\unitlength}{3947sp}%
\begingroup\makeatletter\ifx\SetFigFont\undefined%
\gdef\SetFigFont#1#2#3#4#5{%
  \reset@font\fontsize{#1}{#2pt}%
  \fontfamily{#3}\fontseries{#4}\fontshape{#5}%
  \selectfont}%
\fi\endgroup%
\begin{picture}(944,344)(204,217)
\put(601,314){\makebox(0,0)[lb]{\smash{\SetFigFont{12}{14.4}{\familydefault}{\mddefault}{\updefault}{\color[rgb]{0,0,0}$=$}%
}}}
\end{picture}
}\vcenter{\box1}\!,\ \setbox1=\hbox{\begin{picture}(0,0)%
\includegraphics{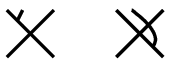}%
\end{picture}%
\setlength{\unitlength}{3947sp}%
\begingroup\makeatletter\ifx\SetFigFont\undefined%
\gdef\SetFigFont#1#2#3#4#5{%
  \reset@font\fontsize{#1}{#2pt}%
  \fontfamily{#3}\fontseries{#4}\fontshape{#5}%
  \selectfont}%
\fi\endgroup%
\begin{picture}(794,269)(279,67)
\put(601,129){\makebox(0,0)[lb]{\smash{\SetFigFont{12}{14.4}{\familydefault}{\mddefault}{\updefault}{\color[rgb]{0,0,0}$=$}%
}}}
\end{picture}
}\vcenter{\box1}\!,
& \label{e:relclear}
\end{eqnarray}

where the pictures in equations~(\ref{e:relverify}), (\ref{e:relclear}) 
represent words in the generators and~$\otimes,\circ$
 such that all compositions are
defined. 

The value of~$W_0$ on objects and generators of~$\Anabo$ has been
defined in equations~(\ref{e:Wot}) to~(\ref{e:WoX}). It remains to
verify the compatibility of~$\Wnab_0$ with the defining relations
of~$\Anabo$.

The relation on the left side of equation~(\ref{e:relverify})
follows from equation~(\ref{e:templieb}).

Assume that the graphs on the right side of
equation~(\ref{e:relverify}) are colored and oriented such that
they are mapped by~$\Wnab_0$ to elements
of~$\Hom_\gloo(V_{t_1}\otimes V_{t_2}, V_{t_3})$ for certain
parameters~$t_i=(\lambda_i,\mu_i,\sigma_i)\in\Q^*\times\Q\times\Z/(2)$.
We will distinguish two cases. For
$\sigma_1+\sigma_2+\sigma_3\equiv 1\mod 2$ we compute

$$ \Wnab_0\left(\setbox1=\hbox{\begin{picture}(0,0)%
\includegraphics{a.eps}%
\end{picture}%
\setlength{\unitlength}{3947sp}%
\begingroup\makeatletter\ifx\SetFigFont\undefined%
\gdef\SetFigFont#1#2#3#4#5{%
  \reset@font\fontsize{#1}{#2pt}%
  \fontfamily{#3}\fontseries{#4}\fontshape{#5}%
  \selectfont}%
\fi\endgroup%
\begin{picture}(294,319)(304,-58)
\end{picture}
}\vcenter{\box1}\!\right) =
\Adot_{t_1,t_2}=(\I_{t_4}\otimes \gcap_{t_2})\circ(\Ystd_{t_4,
t_2^*}\otimes \I_{t_2})=\Wnab_0\left(\setbox1=\hbox{\begin{picture}(0,0)%
\includegraphics{rel7a.eps}%
\end{picture}%
\setlength{\unitlength}{3947sp}%
\begingroup\makeatletter\ifx\SetFigFont\undefined%
\gdef\SetFigFont#1#2#3#4#5{%
  \reset@font\fontsize{#1}{#2pt}%
  \fontfamily{#3}\fontseries{#4}\fontshape{#5}%
  \selectfont}%
\fi\endgroup%
\begin{picture}(429,419)(119,367)
\end{picture}
}\vcenter{\box1}\!\right), $$

where we used the first equation of Lemma~\ref{l:AY} and
equations~(\ref{e:Wot}) to~(\ref{e:WoAo}). For
$\sigma_1+\sigma_2+\sigma_3\equiv 0\mod 2$ we compute

\begin{eqnarray*}
\Wnab_0\left(\setbox1=\hbox{\begin{picture}(0,0)%
\includegraphics{a.eps}%
\end{picture}%
\setlength{\unitlength}{3947sp}%
\begingroup\makeatletter\ifx\SetFigFont\undefined%
\gdef\SetFigFont#1#2#3#4#5{%
  \reset@font\fontsize{#1}{#2pt}%
  \fontfamily{#3}\fontseries{#4}\fontshape{#5}%
  \selectfont}%
\fi\endgroup%
\begin{picture}(294,319)(304,-58)
\end{picture}
}\vcenter{\box1}\!\right) & = &
(-1)^{\sigma_3}\lambda_3\Astd_{t_1,t_2}=
(-1)^{\sigma_2+\sigma_3}\lambda_1(\I_{t_3}\otimes\gcap_{t_2})\circ(\Ydot_{t_3,
t_2^*} \otimes \I_{t_2})\\ & = &
(\I_{t_3}\otimes\gcap_{t_2})\circ\left( (-1)^{\sigma_1}\lambda_1
\Ydot_{t_3, t_2^*} \otimes
\I_{t_2}\right)=\Wnab_0\left(\setbox1=\hbox{\begin{picture}(0,0)%
\includegraphics{rel7a.eps}%
\end{picture}%
\setlength{\unitlength}{3947sp}%
\begingroup\makeatletter\ifx\SetFigFont\undefined%
\gdef\SetFigFont#1#2#3#4#5{%
  \reset@font\fontsize{#1}{#2pt}%
  \fontfamily{#3}\fontseries{#4}\fontshape{#5}%
  \selectfont}%
\fi\endgroup%
\begin{picture}(429,419)(119,367)
\end{picture}
}\vcenter{\box1}\!\right),
\end{eqnarray*}

where we used the second equation of Lemma~\ref{l:AY} this time.
The equation

$$\Wnab_0\left(\setbox1=\hbox{\begin{picture}(0,0)%
\includegraphics{a.eps}%
\end{picture}%
\setlength{\unitlength}{3947sp}%
\begingroup\makeatletter\ifx\SetFigFont\undefined%
\gdef\SetFigFont#1#2#3#4#5{%
  \reset@font\fontsize{#1}{#2pt}%
  \fontfamily{#3}\fontseries{#4}\fontshape{#5}%
  \selectfont}%
\fi\endgroup%
\begin{picture}(294,319)(304,-58)
\end{picture}
}\vcenter{\box1}\!\right)=\Wnab_0\left(\setbox1=\hbox{\begin{picture}(0,0)%
\includegraphics{rel7c.eps}%
\end{picture}%
\setlength{\unitlength}{3947sp}%
\begingroup\makeatletter\ifx\SetFigFont\undefined%
\gdef\SetFigFont#1#2#3#4#5{%
  \reset@font\fontsize{#1}{#2pt}%
  \fontfamily{#3}\fontseries{#4}\fontshape{#5}%
  \selectfont}%
\fi\endgroup%
\begin{picture}(429,419)(1254,-533)
\end{picture}
}\vcenter{\box1}\!\right)$$

is verified similarly. For the verification of the compatibility
of~$\Wnab_0$ with the relations in equation~(\ref{e:relclear}) we
use that~$\Wnab_0$ maps graphs to morphisms of degree~$0$ together
with simple properties of the superpermutation~$\Xt_{t_1,t_2}$ and
the map~$\gcap_t$. This completes the proof.
\end{proof}

The morphisms~$\Ystd_{t_1,t_2}$ 
verify~$\Ystd_{t_1,t_2}=(-1)^{\sigma_1\sigma_2}
\Xt_{t_2, t_1}\Ystd_{t_2,t_1}$
whereas
for~$\Ydot_{t_1,t_2}$
we have $\Ydot_{t_1,t_2}=(-1)^{(\sigma_1+1)(\sigma_2+1)}
\Xt_{t_2, t_1}\circ\Ydot_{t_2,t_1}$. 
Using equation~(\ref{e:WoYo}) 
this translates as follows to trivalent vertices~$v$
of colored trivalent graphs~$\Gamma$ that are incident 
to~$i,j,k\in E(\Gamma)$. We define~$s_v\in\{\pm 1\}$ by
$s_v=1$ iff~$\vert s_{v,i}+s_{v,j}+s_{v,k}\vert=1$. 
Then we have

\begin{equation}\label{e:symantisym}
W_0\left(\setbox1=\hbox{\begin{picture}(0,0)%
\includegraphics{y.eps}%
\end{picture}%
\setlength{\unitlength}{3947sp}%
\begingroup\makeatletter\ifx\SetFigFont\undefined%
\gdef\SetFigFont#1#2#3#4#5{%
  \reset@font\fontsize{#1}{#2pt}%
  \fontfamily{#3}\fontseries{#4}\fontshape{#5}%
  \selectfont}%
\fi\endgroup%
\begin{picture}(294,319)(304,217)
\end{picture}
}\vcenter{\box1}\!\right)=s_vW_0\left(\setbox1=\hbox{\begin{picture}(0,0)%
\includegraphics{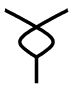}%
\end{picture}%
\setlength{\unitlength}{3947sp}%
\begingroup\makeatletter\ifx\SetFigFont\undefined%
\gdef\SetFigFont#1#2#3#4#5{%
  \reset@font\fontsize{#1}{#2pt}%
  \fontfamily{#3}\fontseries{#4}\fontshape{#5}%
  \selectfont}%
\fi\endgroup%
\begin{picture}(344,394)(204,-133)
\end{picture}
}\vcenter{\box1}\!\right).
\end{equation}

When~$s_v=-1$ we distinguish the cases 
$s_{v,i}+s_{v,j}+s_{v,k}$=-3 where $v$ is called a {\em source},
and $s_{v,i}+s_{v,j}+s_{v,k}=3$ where $v$ is called a {\em sink}.

\section{$\Lambda(\Gamma)$-linear weight systems}\label{s:Wtilde}

Maps from trivalent diagrams to modules are 
called weight systems. In this section we will combine a well-known construction 
of weight systems related to~$\gl(1\vert 1)$ with the results of the previous section.
Define a category~$\Anab$ that has the same objects
as~$\Anabo$. The set $\Hom_\Anab(c,c')$ is a direct sum of
modules~$\Anab(\Gamma)$
where~$\Gamma\in\Hom_{\Anabo}(c,c')$.
The graphical definition of the
composition of morphisms in~$\Anabo$ extends to trivalent diagrams and
induces $\Q$-linear maps

\begin{eqnarray*}
\beta:\Lambda(\Gamma_1)\otimes_\Q\Lambda(\Gamma_2)\longrightarrow\Lambda(\Gamma_1\circ\Gamma_2),
& \mbox{and} &
\alpha:\Anab(\Gamma_1)\otimes_\Q\Anab(\Gamma_2)\longrightarrow\Anab(\Gamma_1\circ\Gamma_2)
\end{eqnarray*}

with the property that for all $p_i\in\Lambda(\Gamma_i)$,
$D_i\in\Anab(\Gamma_i)$ ($i=1,2$) we have

\begin{eqnarray*} \alpha((p_1D_1)\otimes (p_2D_2)) & = & \beta(p_1\otimes
p_2)\alpha(D_1\otimes D_2).
\end{eqnarray*}

We
consider~$\Anabo$ as a subcategory of~$\Anab$ in the obvious way.
For objects~$c$ of~$\Anab$ we define a functor~$W$ by~$W(c)=W_0(c)$.
Now we extend~$\Wnab_0$ to morphisms of~$\Anab$. 
Define~$\omega\in\gloo^{\otimes 2}$ by

\begin{equation}\label{e:defomega}
 \omega=(1/2)(H\otimes D+D\otimes H)+F\otimes E-E\otimes F.
\end{equation}

Sometimes we use the notation $\omega=\sum_\nu
a_\nu\otimes b_\nu$ for~$\omega$. The elements~$\omega$ and~$D\otimes D$ are a
basis of the space of invariants in~$\gloo\otimes\gloo$ by the
adjoint representation.

Consider an object $c=(c_1,\ldots, c_k)$ of~$\Anab$ with
$c_i=(\lambda_i,\mu_i,\sigma_i)$. Let~$\Gamma=\id_c\in\Anab_0$.
Define $T_c\in\End_\Anab(c)$ as 
the trivalent diagram of degree~$1$ on~$\Gamma$ that connects
the first interval of~$\Gamma$ with the second interval. We
define $W(T_c)\in\End_\gloo(W(c))$ by

\begin{equation}\label{e:defWomega}
W(T_c)(v_1\otimes\ldots \otimes v_k)=\sum_\nu (-1)^{\deg
(v_1)\deg(b_\nu)+\sigma_1+\sigma_2}a_\nu\cdot v_1\otimes
b_\nu\cdot v_2\otimes v_3\otimes\ldots\otimes v_k,
\end{equation}

where $v_i\in W(c_i)$.\footnote{The explicit appearence of the factor 
$(-1)^{\sigma_1+\sigma_2}$ in 
equation~(\ref{e:defWomega}) is sometimes avoided by introducing an antisymmetry relation 
for univalent
vertices that are glued to~$\Gamma$.} For morphisms~$D$ of~$\Anabo$ we
define~$W(D)=W_0(D)$.
For any unitrivalent graph~$\Gamma$ with admissible coloring, we
consider~$\gloo$-modules as a~$\Lambda(\Gamma)$-module, where
$d_i$ acts by multiplication with~$\lambda_i$ when
$(\lambda_i,\mu_i)$ is the color of the edge~$i$. 
Equation~(\ref{e:admissibilitya}) implies that this definition is compatible with
the defining relations of~$\Lambda(\Gamma)$.
We have the
following lemma.

\begin{lemma}\label{l:Wnab}
The definition of~$\Wnab$ extends uniquely to a functor from~$\Anab$
to~$\gloo$-modules that is $\Lambda(\Gamma)$-linear
on~$\Anab(\Gamma)$ and~$\Q$-linear\ on general morphisms.
\end{lemma}
\begin{proof}
Using relation~$(STU)$ we see that the modules~$\Anab(\Gamma)$ are
generated by diagrams~$D$ such that~$D\setminus\Gamma$ contains no trivalent
vertex. For the diagrams~$D$ one can choose pictures such that horizontal stripes
around intervals in~$D\setminus\Gamma$ are equal to~$T_c^{12}$ for various objects~$c$ of
$\Anab$. Linearity of~$\Wnab$ now implies that the definition of~$\Wnab$ determines the
value of all morphisms of~$\Anab$.
We have to verify
that~$\Wnab$ is independent of the choices from above and
compatible with the defining relations of~$\Anab(\Gamma)$.

It follows from Lemma~\ref{l:W0} and 
general properties of~$\omega$ ($\omega$ is invariant, supersymmetric, and of degree~$0$) 
that the definition of~$\Wnab$ induces 
well-defined $\Q$-linear maps~$\Wnab'_\Gamma$ from $\A(\Gamma)$ to $\Q$ for all
admissibly colored unitrivalent graphs~$\Gamma$. 
It remains to show that
the~$\Lambda(\Gamma)$-linear 
maps~$\Wnab_\Gamma:\Anab(\Gamma)\longrightarrow\Q$ defined 
by~$\Wnab_\Gamma\circ \tau=\Wnab'_\Gamma$ are well-defined, where~$\tau$ denotes the canonical
map from~$\A(\Gamma)$ to~$\Anabc(\Gamma)$.
As in the proof of Lemma~\ref{l:Gnlinind2}
this follows from particular properties of~$\gl(1\vert 1)$ and~$\omega$~:
equations~(\ref{e:Wprop1}) and (\ref{e:Wprop2}) (see~\cite{FKV})
imply the compatibility of~$\Wnab_\Gamma$ with relations~$(CL1A)$ and~$(CL2A)$. 

\begin{eqnarray}
& & \Wnab_\Gamma\!\left(\setbox1=\hbox{\input{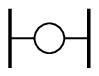}}\vcenter{\box1}\!\right)\, = \,  
-2d_id_j\,\Wnab_\Gamma\!\left(\setbox1=\hbox{\input{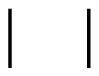}}\vcenter{\box1}\!\right) ,\ \
\Wnab_\Gamma\!\left(\setbox1=\hbox{\begin{picture}(0,0)%
\includegraphics{inty.eps}%
\end{picture}%
\setlength{\unitlength}{4144sp}%
\begingroup\makeatletter\ifx\SetFigFont\undefined%
\gdef\SetFigFont#1#2#3#4#5{%
  \reset@font\fontsize{#1}{#2pt}%
  \fontfamily{#3}\fontseries{#4}\fontshape{#5}%
  \selectfont}%
\fi\endgroup%
\begin{picture}(329,303)(439,98)
\end{picture}
}\vcenter{\box1}\!\right)\, =\, 0,\ \ 
\Wnab_\Gamma\!\left(\setbox1=\hbox{\begin{picture}(0,0)%
\includegraphics{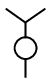}%
\end{picture}%
\setlength{\unitlength}{4144sp}%
\begingroup\makeatletter\ifx\SetFigFont\undefined%
\gdef\SetFigFont#1#2#3#4#5{%
  \reset@font\fontsize{#1}{#2pt}%
  \fontfamily{#3}\fontseries{#4}\fontshape{#5}%
  \selectfont}%
\fi\endgroup%
\begin{picture}(204,339)(889,-253)
\end{picture}
}\vcenter{\box1}\!\right)\, =\, 0,\ \mbox{}
\label{e:Wprop1}\\
& & \Wnab_\Gamma\!\left(\setbox1=\hbox{\begin{picture}(0,0)%
\includegraphics{vh.eps}%
\end{picture}%
\setlength{\unitlength}{4144sp}%
\begingroup\makeatletter\ifx\SetFigFont\undefined%
\gdef\SetFigFont#1#2#3#4#5{%
  \reset@font\fontsize{#1}{#2pt}%
  \fontfamily{#3}\fontseries{#4}\fontshape{#5}%
  \selectfont}%
\fi\endgroup%
\begin{picture}(294,294)(394,197)
\end{picture}
}\vcenter{\box1}\!\right)\ = \ \frac{1}{2}\,
\Wnab_\Gamma\!\left(\setbox1=\hbox{\begin{picture}(0,0)%
\includegraphics{vt3.eps}%
\end{picture}%
\setlength{\unitlength}{4144sp}%
\begingroup\makeatletter\ifx\SetFigFont\undefined%
\gdef\SetFigFont#1#2#3#4#5{%
  \reset@font\fontsize{#1}{#2pt}%
  \fontfamily{#3}\fontseries{#4}\fontshape{#5}%
  \selectfont}%
\fi\endgroup%
\begin{picture}(294,294)(664,-208)
\end{picture}
}\vcenter{\box1}\!+\setbox1=\hbox{\begin{picture}(0,0)%
\includegraphics{vt4.eps}%
\end{picture}%
\setlength{\unitlength}{4144sp}%
\begingroup\makeatletter\ifx\SetFigFont\undefined%
\gdef\SetFigFont#1#2#3#4#5{%
  \reset@font\fontsize{#1}{#2pt}%
  \fontfamily{#3}\fontseries{#4}\fontshape{#5}%
  \selectfont}%
\fi\endgroup%
\begin{picture}(294,294)(664,-388)
\end{picture}
}\vcenter{\box1}\!-\setbox1=\hbox{\begin{picture}(0,0)%
\includegraphics{vt1.eps}%
\end{picture}%
\setlength{\unitlength}{4144sp}%
\begingroup\makeatletter\ifx\SetFigFont\undefined%
\gdef\SetFigFont#1#2#3#4#5{%
  \reset@font\fontsize{#1}{#2pt}%
  \fontfamily{#3}\fontseries{#4}\fontshape{#5}%
  \selectfont}%
\fi\endgroup%
\begin{picture}(294,294)(664,197)
\end{picture}
}\vcenter{\box1}\!-\setbox1=\hbox{\begin{picture}(0,0)%
\includegraphics{vt2.eps}%
\end{picture}%
\setlength{\unitlength}{4144sp}%
\begingroup\makeatletter\ifx\SetFigFont\undefined%
\gdef\SetFigFont#1#2#3#4#5{%
  \reset@font\fontsize{#1}{#2pt}%
  \fontfamily{#3}\fontseries{#4}\fontshape{#5}%
  \selectfont}%
\fi\endgroup%
\begin{picture}(294,294)(664,17)
\end{picture}
}\vcenter{\box1}\!\right). \label{e:Wprop2}
\end{eqnarray}\nopagebreak{}\end{proof}
%
%
%

Now consider a trivalent diagram~$D$ on~$\Gamma$, where~$\Gamma$
is an admissibly colored trivalent graph. Let $p\in \Gamma\setminus V(D)$ 
be a point on an 
edge with color~$(\lambda_p,\mu_p)$. By cutting~$D$ at~$p$,
we obtain a diagram~$D_p\in\End_\Anab((\lambda_p,\mu_p,0))$.

\begin{lemma}\label{l:cutnotimp1}
%
%
There exists a unique linear map~$\Wnabt:\End_{\Anab}(\emptyset)\longrightarrow\Q$
such that for a trivalent diagram~$D$ on an admissibly colored trivalent graph~$\Gamma$
and any~$p\in\Gamma\setminus V(D)$ on an edge with color~$(\lambda_p,\mu_p)$ we have

$$
W\left(D_{p}\right)=\lambda_p\Wnabt(D)\I_{t_p}\ ,\ 
\mbox{where $t_p=(\lambda_p,\mu_p,0)$.}$$

The map~$\Wnabt$ is~$\Lambda(\Gamma)$-linear on 
$\Anab(\Gamma)\subset\End_{\Anab}(\emptyset)$.
\end{lemma}
\begin{proof}
We only have to show that the definition of~$\Wnabt(D)$ does not
depend on the choice of~$p$, because
then Lemma~\ref{l:Wnab}
implies the compatibility
of~$\Wnabt$ with the defining relations 
of~$\Anab(\Gamma)\subset\End_{\Anab}(\emptyset)$ and the~$\Lambda(\Gamma)$-linearity.

Consider two points~$p_1, p_2\in\Gamma\setminus V(D)$ on
edges with colors $(\lambda_i, \mu_i)$ ($i=1,2$).
We treat the case~$\lambda_1\not=\lambda_2$
in detail.
Let $t_i=(\lambda_i, \mu_i, 0)$.
Pictorial
representations of the diagrams~$D$, $D_{p_1}$,
and~$D_{p_2}$ are shown in equation~(\ref{e:2diags}),
where the box labeled $x$ represents a
morphism~$x\in\End_{\Anab}((t_1,t_2^*))$.

\begin{equation}\label{e:2diags}
D\ =\ \setbox1=\hbox{\input{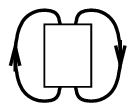}}\vcenter{\box1}\! \ \ ,\ \ D_{p_1}\ = \ \setbox1=\hbox{\input{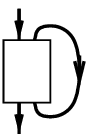}}\vcenter{\box1}\!
\ \ ,\ \ D_{p_2}\ =\ \setbox1=\hbox{\input{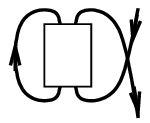}}\vcenter{\box1}\!.
\end{equation}

It follows from equations~(\ref{e:end0bas}), (\ref{e:WoYo}),
(\ref{e:WoAo}) and definitions that there exist $a,b\in\Q$ such
that

$$
W(x)=a\Ystd_{t_1,t_2^*}\circ\Astd_{t_1,t_2^*}+b\Ydot_{t_1,t_2^*}\circ\Adot_{t_1,t_2^*}=
\frac{a}{\lambda_1-\lambda_2}W\left(\setbox1=\hbox{\begin{picture}(0,0)%
\includegraphics{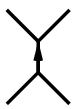}%
\end{picture}%
\setlength{\unitlength}{3947sp}%
\begingroup\makeatletter\ifx\SetFigFont\undefined%
\gdef\SetFigFont#1#2#3#4#5{%
  \reset@font\fontsize{#1}{#2pt}%
  \fontfamily{#3}\fontseries{#4}\fontshape{#5}%
  \selectfont}%
\fi\endgroup%
\begin{picture}(344,494)(579,-308)
\end{picture}
}\vcenter{\box1}\!\right)+\frac{b}{\lambda_1-\lambda_2}W\left(\setbox1=\hbox{\begin{picture}(0,0)%
\includegraphics{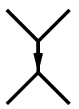}%
\end{picture}%
\setlength{\unitlength}{3947sp}%
\begingroup\makeatletter\ifx\SetFigFont\undefined%
\gdef\SetFigFont#1#2#3#4#5{%
  \reset@font\fontsize{#1}{#2pt}%
  \fontfamily{#3}\fontseries{#4}\fontshape{#5}%
  \selectfont}%
\fi\endgroup%
\begin{picture}(344,494)(579,-308)
\end{picture}
}\vcenter{\box1}\!\right)
$$


which implies

\begin{equation}\label{e:Wtilde} W(D_{p_i}) = \frac{a}{\lambda_1-\lambda_2}W\left(\setbox1=\hbox{\begin{picture}(0,0)%
\includegraphics{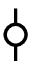}%
\end{picture}%
\setlength{\unitlength}{3947sp}%
\begingroup\makeatletter\ifx\SetFigFont\undefined%
\gdef\SetFigFont#1#2#3#4#5{%
  \reset@font\fontsize{#1}{#2pt}%
  \fontfamily{#3}\fontseries{#4}\fontshape{#5}%
  \selectfont}%
\fi\endgroup%
\begin{picture}(130,294)(761,-58)
\end{picture}
}\vcenter{\box1}\!\right)+
\frac{b}{\lambda_1-\lambda_2}W\left(\setbox1=\hbox{\begin{picture}(0,0)%
\includegraphics{o.eps}%
\end{picture}%
\setlength{\unitlength}{3947sp}%
\begingroup\makeatletter\ifx\SetFigFont\undefined%
\gdef\SetFigFont#1#2#3#4#5{%
  \reset@font\fontsize{#1}{#2pt}%
  \fontfamily{#3}\fontseries{#4}\fontshape{#5}%
  \selectfont}%
\fi\endgroup%
\begin{picture}(130,294)(761,-58)
\end{picture}
}\vcenter{\box1}\!\right)
=\lambda_i\frac{a+b}{\lambda_1-\lambda_2}\I_{t_i}.
\end{equation}

The
value~$\frac{a+b}{\lambda_1-\lambda_2}$ in
equation~(\ref{e:Wtilde}) does not depend on~$p_i$ ($i\in\{1,2\}$). 
It follows
that~$\Wnabt(D)$ is well-defined for diagrams~$D$ whose
colored skeleton has two edges colored by~$(\lambda'_i,\mu'_i)$ ($i=1,2$) 
with~$\lambda'_1\not=\lambda'_2$ because then we can treat the case~$\lambda_1=\lambda_2$
by a two-fold application of the argument above.

The remaining case is well-known from computations concerning the $1$-variable Alexander
polynomial of links.
Alternatively, there is a proof for $\lambda_1=\lambda_2$ 
(or, more generally, $\lambda_1\not=-\lambda_2$)
similar to the proof above, where a picture of~$D$ is chosen as in 
equation~(\ref{e:2diags}), but with
different orientations and with~$x\in\End_\Anab((t_1,t_2))$. 
\end{proof}

\section{The Kontsevich integral of unitrivalent graphs}

Recall from~\cite{LeM} 
the definition of the Kontsevich integral~$Z$ of 
framed $q$-tangles. 
We assume that~$Z$ is defined using an even group-like horizontal Drinfeld associator 
in~$\Ac(3)$. For a $q$-tangle~$T$ we have~$Z(T)\in\Ac(\Gamma(T))$ where~$\Gamma(T)$ is 
the underlying $1$-manifold of~$T$.
In this section we will study an extension of~$Z$ to graphs (compare~\cite{MuO}).
We consider a 
category~${\cal G^{\na}}$ whose morphisms are isotopy classes of
oriented \mbox{(half-)}framed unitrivalent graphs~$G$
with cyclically oriented vertices. By definition,
$G$ is properly embedded into~$\R\times [0,1]\times \R$ and 
we have

$$G\cap \R\times\{i\}\times \R=
\{1,\ldots,n_i\}\times\{i\}\times\{0\}$$ 

for $i=0,1$ and for certain~$n_0,n_1\geq 0$. 
We represent a strand of~$G$ 
with a right-handed half twist of the framing 
graphically by~$\setbox1=\hbox{\begin{picture}(0,0)%
\includegraphics{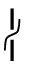}%
\end{picture}%
\setlength{\unitlength}{3947sp}%
\begingroup\makeatletter\ifx\SetFigFont\undefined%
\gdef\SetFigFont#1#2#3#4#5{%
  \reset@font\fontsize{#1}{#2pt}%
  \fontfamily{#3}\fontseries{#4}\fontshape{#5}%
  \selectfont}%
\fi\endgroup%
\begin{picture}(94,294)(779,-58)
\end{picture}
}\vcenter{\box1}\!$.
The objects of~${\cal G^{\na}}$ are non-associative words in 
the symbols~$+$ and~$-$. 
For example, $(-(+-))$ is an object of~${\cal G^\na}$.
Unitrivalent graphs~$G\in\Hom_{{\cal G}^\na}(w_0,w_1)$
are related to~$w_i$ as follows~:
the~$i$-th symbol~$a\in\{+,-\}$ of~$w_0$ (resp.\ $w_1$) 
corresponds to the~$i$-th lower (resp.\ upper) 
boundary point~$p=(i,0,0)$ (resp.\ $p=(i,1,0)$) of~$G$ 
where $a=+$ (resp.\ $a=-$) means that the
graph~$G$ must be oriented downwards (resp.\ upwards) at~$p$.

The category~$\Ac^\na$ 
has the same objects as~${\cal G^{\na}}$. 
The set $\Hom_{\Ac^\na}(w_0, w_1)$ is the direct sum of all modules~$\Ac(\Gamma)$
where~$\Gamma$ is a unitrivalent graph whose boundary is partitioned into two ordered
sets called lower and upper boundary as in Section~\ref{s:W0},
and $a\in\{+,-\}$ in~$w_0$ and~$w_1$ is related to the 
boundary points in a graphical representation of~$\Gamma$ in the same way
as above.
We consider the invariant~$Z$ of $q$-tangles from~\cite{LeM} as a tensor functor
from the subcategory~${\cal T}^\na\subset{\cal G}^\na$ of~$q$-tangles to~$\Ac^\na$.

Let us recall how~$Z$ depends on the orientations of edges of $q$-tangles.
Let~$\Gamma$ be an oriented~$1$-dimensional 
manifold with boundary and let~$\Gamma'\subset\Gamma$ be
a set of connected components of~$\Gamma$. Let~$D$ be a trivalent diagram on~$\Gamma$. 
Define~$S_{\Gamma'}(D)$ by inverting the orientation of all components 
of~$\Gamma'$ and by 
multiplying the result by~$(-1)^m$ where~$m$ is the number of univalent vertices
of~$D\setminus\Gamma$
that are glued to~$\Gamma'$.
Let~$T$ be a $q$-tangle and $T'\subset T$ be a set of connected components of~$T$.
Define~$S_{T'}(T)$ by inverting the orientation of all components of~$T'\subset T$.
With this notation we have

\begin{equation}\label{e:orinv}
Z\left(S_{T'}(T)\right)=S_{\Gamma(T')}\left(Z(T)\right).
\end{equation}

We omit the index of~$S$ when all components of a unitrivalent graph (resp.\ of the skeleton
of a trivalent diagram) are concerned. For~$a\in\End_{\Ac}(+)$ we 
define

\begin{equation}\label{e:invert}
S^*(a)=(\setbox1=\hbox{\begin{picture}(0,0)%
\includegraphics{rightcap.eps}%
\end{picture}%
\setlength{\unitlength}{3947sp}%
\begingroup\makeatletter\ifx\SetFigFont\undefined%
\gdef\SetFigFont#1#2#3#4#5{%
  \reset@font\fontsize{#1}{#2pt}%
  \fontfamily{#3}\fontseries{#4}\fontshape{#5}%
  \selectfont}%
\fi\endgroup%
\begin{picture}(359,263)(579,217)
\end{picture}
}\vcenter{\box1}\!\otimes\setbox1=\hbox{\begin{picture}(0,0)%
\includegraphics{up.eps}%
\end{picture}%
\setlength{\unitlength}{3947sp}%
\begingroup\makeatletter\ifx\SetFigFont\undefined%
\gdef\SetFigFont#1#2#3#4#5{%
  \reset@font\fontsize{#1}{#2pt}%
  \fontfamily{#3}\fontseries{#4}\fontshape{#5}%
  \selectfont}%
\fi\endgroup%
\begin{picture}(74,294)(339,267)
\end{picture}
}\vcenter{\box1}\!)\circ(\setbox1=\hbox{\begin{picture}(0,0)%
\includegraphics{up.eps}%
\end{picture}%
\setlength{\unitlength}{3947sp}%
\begingroup\makeatletter\ifx\SetFigFont\undefined%
\gdef\SetFigFont#1#2#3#4#5{%
  \reset@font\fontsize{#1}{#2pt}%
  \fontfamily{#3}\fontseries{#4}\fontshape{#5}%
  \selectfont}%
\fi\endgroup%
\begin{picture}(74,294)(339,267)
\end{picture}
}\vcenter{\box1}\!\otimes a\otimes 
\setbox1=\hbox{\begin{picture}(0,0)%
\includegraphics{up.eps}%
\end{picture}%
\setlength{\unitlength}{3947sp}%
\begingroup\makeatletter\ifx\SetFigFont\undefined%
\gdef\SetFigFont#1#2#3#4#5{%
  \reset@font\fontsize{#1}{#2pt}%
  \fontfamily{#3}\fontseries{#4}\fontshape{#5}%
  \selectfont}%
\fi\endgroup%
\begin{picture}(74,294)(339,267)
\end{picture}
}\vcenter{\box1}\!)\circ(\setbox1=\hbox{\begin{picture}(0,0)%
\includegraphics{up.eps}%
\end{picture}%
\setlength{\unitlength}{3947sp}%
\begingroup\makeatletter\ifx\SetFigFont\undefined%
\gdef\SetFigFont#1#2#3#4#5{%
  \reset@font\fontsize{#1}{#2pt}%
  \fontfamily{#3}\fontseries{#4}\fontshape{#5}%
  \selectfont}%
\fi\endgroup%
\begin{picture}(74,294)(339,267)
\end{picture}
}\vcenter{\box1}\!\otimes\setbox1=\hbox{\begin{picture}(0,0)%
\includegraphics{rightcup.eps}%
\end{picture}%
\setlength{\unitlength}{3947sp}%
\begingroup\makeatletter\ifx\SetFigFont\undefined%
\gdef\SetFigFont#1#2#3#4#5{%
  \reset@font\fontsize{#1}{#2pt}%
  \fontfamily{#3}\fontseries{#4}\fontshape{#5}%
  \selectfont}%
\fi\endgroup%
\begin{picture}(359,263)(579,298)
\end{picture}
}\vcenter{\box1}\!).
\end{equation}

It is unknown 
if
there exist elements~$a\in\Ac(1)\subset \End_{\Ac}(+)$ with~$S(a)\not=S^*(a)$.

In diagrams of morphisms
$G$ of ${\cal G}^\na$ we use 
projections of generic representatives of~$G$ to the first two coordinates such that 
the cyclic order at trivalent vertices is counterclockwise in the projection. 
For a unitrivalent graph~$\Gamma$ we define~$\Gamma^\op$ by inverting the cyclic
order of all trivalent vertices of~$\Gamma$.
This induces maps
$$
\Ac(\Gamma)\longrightarrow\Ac(\Gamma^\op),\ a\mapsto a^\op
\mbox{ and }
\Hom_{\cal G^\na}(w_0, w_1)\longrightarrow\Hom_{\cal G^\na}(w_0, w_1),\ G\mapsto G^\op.
$$

We denote the Kontsevich integral of the trivial knot by~$\nu$ and regard
it as an element of~$\Ac(1)$. 
Since~$\nu$ is equal to~$1$ in degree~$0$ the element~$\nu$ is invertible and
there exist unique roots of~$\nu^{n}$ ($n\in\Z$) that are equal to~$1$ in degree~$0$.
The elements~$\nu^k$ ($k\in\Q$) 
satisfy~$S(\nu^k)=S^*(\nu^k)$. In the formulas below the box labeled~$\nu^k$
represents the element~$\nu^k$ or~$S(\nu^k)$ according
to the orientation of the interval with the box.

\begin{theorem}\label{t:Ztriv}
(1) For any choice of

$$a,b\in\End_{\Ac^\na}(+)\quad ,\quad 
c\in\Hom_{\Ac^\na}(-,(++))\quad\mbox{,\ and}\quad
d\in\Hom_{\Ac^\na}(+,(--))$$ 

there exists a unique extension of~$Z$ to a tensor functor 
from~$\G^\na$ to~$\Ac^\na$ satisfying

\begin{eqnarray*}
Z\left(\setbox1=\hbox{\begin{picture}(0,0)%
\includegraphics{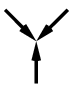}%
\end{picture}%
\setlength{\unitlength}{3947sp}%
\begingroup\makeatletter\ifx\SetFigFont\undefined%
\gdef\SetFigFont#1#2#3#4#5{%
  \reset@font\fontsize{#1}{#2pt}%
  \fontfamily{#3}\fontseries{#4}\fontshape{#5}%
  \selectfont}%
\fi\endgroup%
\begin{picture}(344,396)(404,116)
\end{picture}
}\vcenter{\box1}\!\right)=\setbox1=\hbox{\input{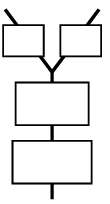}}\vcenter{\box1}\!\ , & &
Z\left(\setbox1=\hbox{\begin{picture}(0,0)%
\includegraphics{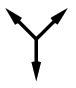}%
\end{picture}%
\setlength{\unitlength}{3947sp}%
\begingroup\makeatletter\ifx\SetFigFont\undefined%
\gdef\SetFigFont#1#2#3#4#5{%
  \reset@font\fontsize{#1}{#2pt}%
  \fontfamily{#3}\fontseries{#4}\fontshape{#5}%
  \selectfont}%
\fi\endgroup%
\begin{picture}(344,396)(404,116)
\end{picture}
}\vcenter{\box1}\!\right)=\setbox1=\hbox{\input{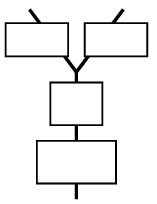}}\vcenter{\box1}\! \ ,\\
Z\left(\setbox1=\hbox{\begin{picture}(0,0)%
\includegraphics{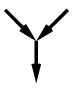}%
\end{picture}%
\setlength{\unitlength}{3947sp}%
\begingroup\makeatletter\ifx\SetFigFont\undefined%
\gdef\SetFigFont#1#2#3#4#5{%
  \reset@font\fontsize{#1}{#2pt}%
  \fontfamily{#3}\fontseries{#4}\fontshape{#5}%
  \selectfont}%
\fi\endgroup%
\begin{picture}(344,396)(404,116)
\end{picture}
}\vcenter{\box1}\!\right)=c\ , & &
Z\left(\setbox1=\hbox{\begin{picture}(0,0)%
\includegraphics{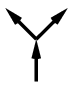}%
\end{picture}%
\setlength{\unitlength}{3947sp}%
\begingroup\makeatletter\ifx\SetFigFont\undefined%
\gdef\SetFigFont#1#2#3#4#5{%
  \reset@font\fontsize{#1}{#2pt}%
  \fontfamily{#3}\fontseries{#4}\fontshape{#5}%
  \selectfont}%
\fi\endgroup%
\begin{picture}(344,396)(404,116)
\end{picture}
}\vcenter{\box1}\!\right)=d\ , \\[3mm]
\mbox{and} \quad Z\left(\setbox1=\hbox{\begin{picture}(0,0)%
\includegraphics{htwist.eps}%
\end{picture}%
\setlength{\unitlength}{3947sp}%
\begingroup\makeatletter\ifx\SetFigFont\undefined%
\gdef\SetFigFont#1#2#3#4#5{%
  \reset@font\fontsize{#1}{#2pt}%
  \fontfamily{#3}\fontseries{#4}\fontshape{#5}%
  \selectfont}%
\fi\endgroup%
\begin{picture}(94,294)(779,-58)
\end{picture}
}\vcenter{\box1}\!\right) & = & \exp\left(\setbox1=\hbox{\begin{picture}(0,0)%
\includegraphics{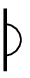}%
\end{picture}%
\setlength{\unitlength}{3947sp}%
\begingroup\makeatletter\ifx\SetFigFont\undefined%
\gdef\SetFigFont#1#2#3#4#5{%
  \reset@font\fontsize{#1}{#2pt}%
  \fontfamily{#3}\fontseries{#4}\fontshape{#5}%
  \selectfont}%
\fi\endgroup%
\begin{picture}(106,344)(-996,1117)
\end{picture}
}\vcenter{\box1}\!/4\right)\ .
\end{eqnarray*}

(2)  When~$\setbox1=\hbox{\begin{picture}(0,0)%
\includegraphics{x.eps}%
\end{picture}%
\setlength{\unitlength}{3947sp}%
\begingroup\makeatletter\ifx\SetFigFont\undefined%
\gdef\SetFigFont#1#2#3#4#5{%
  \reset@font\fontsize{#1}{#2pt}%
  \fontfamily{#3}\fontseries{#4}\fontshape{#5}%
  \selectfont}%
\fi\endgroup%
\begin{picture}(244,244)(329,267)
\end{picture}
}\vcenter{\box1}\!\circ c=c^{\rm op}$ and $\setbox1=\hbox{\begin{picture}(0,0)%
\includegraphics{x.eps}%
\end{picture}%
\setlength{\unitlength}{3947sp}%
\begingroup\makeatletter\ifx\SetFigFont\undefined%
\gdef\SetFigFont#1#2#3#4#5{%
  \reset@font\fontsize{#1}{#2pt}%
  \fontfamily{#3}\fontseries{#4}\fontshape{#5}%
  \selectfont}%
\fi\endgroup%
\begin{picture}(244,244)(329,267)
\end{picture}
}\vcenter{\box1}\!\circ d=d^{\rm op}$ holds 
then we have~$Z(\mbox{$G$}^\op)=\mbox{$Z(G)$}^\op$ for all morphisms~$G$ of~${\cal G}^\na$. 
%
%
%
%
\end{theorem}
\begin{skoproof}
First we consider oriented graphs
whose vertices are oriented boxes with a distinguished lower
boundary (called {\em coupons}). 
A category~${\cal G'}^\na$ is defined 
in the same way as~$\G^\na$ except that morphisms
are embedded framed 
graphs with coupons

$$
C_a=\mbox{$\setbox1=\hbox{\input{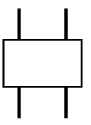}}\vcenter{\box1}\!$}\in\Hom_{{\cal G'}^\na}(s,t)
$$

that are colored by elements~$a\in\Hom_{\Ac^\na}(s, t)$.
It is easy to see that there exists a unique extension of~$Z$ to a 
tensor functor from~${\cal G'}^\na$ to~$\Ac^\na$ that verifies~$Z(C_a)=a$ for all 
coupons~$C_a$ as above.
This general construction implies that part~(1) of the theorem holds iff

$$
Z\left(\setbox1=\hbox{\begin{picture}(0,0)%
\includegraphics{tyddu.eps}%
\end{picture}%
\setlength{\unitlength}{3947sp}%
\begingroup\makeatletter\ifx\SetFigFont\undefined%
\gdef\SetFigFont#1#2#3#4#5{%
  \reset@font\fontsize{#1}{#2pt}%
  \fontfamily{#3}\fontseries{#4}\fontshape{#5}%
  \selectfont}%
\fi\endgroup%
\begin{picture}(344,396)(404,116)
\end{picture}
}\vcenter{\box1}\!\right)=Z\left(\setbox1=\hbox{\begin{picture}(0,0)%
\includegraphics{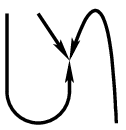}%
\end{picture}%
\setlength{\unitlength}{3947sp}%
\begingroup\makeatletter\ifx\SetFigFont\undefined%
\gdef\SetFigFont#1#2#3#4#5{%
  \reset@font\fontsize{#1}{#2pt}%
  \fontfamily{#3}\fontseries{#4}\fontshape{#5}%
  \selectfont}%
\fi\endgroup%
\begin{picture}(569,585)(579,-233)
\end{picture}
}\vcenter{\box1}\!\right)\quad\mbox{and}\quad 
Z\left(\setbox1=\hbox{\begin{picture}(0,0)%
\includegraphics{tyuud.eps}%
\end{picture}%
\setlength{\unitlength}{3947sp}%
\begingroup\makeatletter\ifx\SetFigFont\undefined%
\gdef\SetFigFont#1#2#3#4#5{%
  \reset@font\fontsize{#1}{#2pt}%
  \fontfamily{#3}\fontseries{#4}\fontshape{#5}%
  \selectfont}%
\fi\endgroup%
\begin{picture}(344,396)(404,116)
\end{picture}
}\vcenter{\box1}\!\right)=Z\left(\setbox1=\hbox{\begin{picture}(0,0)%
\includegraphics{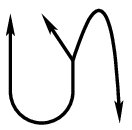}%
\end{picture}%
\setlength{\unitlength}{3947sp}%
\begingroup\makeatletter\ifx\SetFigFont\undefined%
\gdef\SetFigFont#1#2#3#4#5{%
  \reset@font\fontsize{#1}{#2pt}%
  \fontfamily{#3}\fontseries{#4}\fontshape{#5}%
  \selectfont}%
\fi\endgroup%
\begin{picture}(595,585)(564,-1733)
\end{picture}
}\vcenter{\box1}\!\right).$$

These two equations follow from 
a well-known identity in~$\Hom_{\Ac}(((++)-),+)$~:

\begin{equation}
Z\left(\setbox1=\hbox{\begin{picture}(0,0)%
\includegraphics{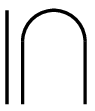}%
\end{picture}%
\setlength{\unitlength}{3947sp}%
\begingroup\makeatletter\ifx\SetFigFont\undefined%
\gdef\SetFigFont#1#2#3#4#5{%
  \reset@font\fontsize{#1}{#2pt}%
  \fontfamily{#3}\fontseries{#4}\fontshape{#5}%
  \selectfont}%
\fi\endgroup%
\begin{picture}(419,494)(429,-383)
\end{picture}
}\vcenter{\box1}\!\right)=\setbox1=\hbox{\input{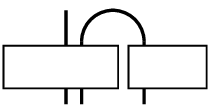}}\vcenter{\box1}\!.
\end{equation}

\smallskip

(2) The isotopy invariance of~$Z$ implies

$$
Z\left(\setbox1=\hbox{\begin{picture}(0,0)%
\includegraphics{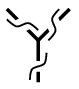}%
\end{picture}%
\setlength{\unitlength}{3947sp}%
\begingroup\makeatletter\ifx\SetFigFont\undefined%
\gdef\SetFigFont#1#2#3#4#5{%
  \reset@font\fontsize{#1}{#2pt}%
  \fontfamily{#3}\fontseries{#4}\fontshape{#5}%
  \selectfont}%
\fi\endgroup%
\begin{picture}(344,394)(729,-133)
\end{picture}
}\vcenter{\box1}\!\right)=Z\left(\mbox{$\setbox1=\hbox{\begin{picture}(0,0)%
\includegraphics{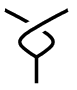}%
\end{picture}%
\setlength{\unitlength}{3947sp}%
\begingroup\makeatletter\ifx\SetFigFont\undefined%
\gdef\SetFigFont#1#2#3#4#5{%
  \reset@font\fontsize{#1}{#2pt}%
  \fontfamily{#3}\fontseries{#4}\fontshape{#5}%
  \selectfont}%
\fi\endgroup%
\begin{picture}(344,394)(204,-133)
\end{picture}
}\vcenter{\box1}\!$}^{\rm op}\right).
$$

When the upper two strands at a trivalent vertex
are both oriented downwards or both oriented
upwards then the equation

$$
\setbox1=\hbox{\begin{picture}(0,0)%
\includegraphics{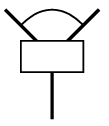}%
\end{picture}%
\setlength{\unitlength}{3947sp}%
\begingroup\makeatletter\ifx\SetFigFont\undefined%
\gdef\SetFigFont#1#2#3#4#5{%
  \reset@font\fontsize{#1}{#2pt}%
  \fontfamily{#3}\fontseries{#4}\fontshape{#5}%
  \selectfont}%
\fi\endgroup%
\begin{picture}(494,569)(504,-383)
\end{picture}
}\vcenter{\box1}\!\ = \ \frac{1}{2}\left(\setbox1=\hbox{\begin{picture}(0,0)%
\includegraphics{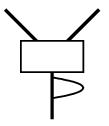}%
\end{picture}%
\setlength{\unitlength}{3947sp}%
\begingroup\makeatletter\ifx\SetFigFont\undefined%
\gdef\SetFigFont#1#2#3#4#5{%
  \reset@font\fontsize{#1}{#2pt}%
  \fontfamily{#3}\fontseries{#4}\fontshape{#5}%
  \selectfont}%
\fi\endgroup%
\begin{picture}(494,569)(504,-383)
\end{picture}
}\vcenter{\box1}\!-\setbox1=\hbox{\begin{picture}(0,0)%
\includegraphics{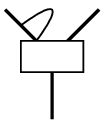}%
\end{picture}%
\setlength{\unitlength}{3947sp}%
\begingroup\makeatletter\ifx\SetFigFont\undefined%
\gdef\SetFigFont#1#2#3#4#5{%
  \reset@font\fontsize{#1}{#2pt}%
  \fontfamily{#3}\fontseries{#4}\fontshape{#5}%
  \selectfont}%
\fi\endgroup%
\begin{picture}(494,569)(504,-383)
\end{picture}
}\vcenter{\box1}\!
-\setbox1=\hbox{\begin{picture}(0,0)%
\includegraphics{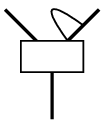}%
\end{picture}%
\setlength{\unitlength}{3947sp}%
\begingroup\makeatletter\ifx\SetFigFont\undefined%
\gdef\SetFigFont#1#2#3#4#5{%
  \reset@font\fontsize{#1}{#2pt}%
  \fontfamily{#3}\fontseries{#4}\fontshape{#5}%
  \selectfont}%
\fi\endgroup%
\begin{picture}(494,569)(1104,-383)
\end{picture}
}\vcenter{\box1}\!\right)
$$

and the symmetry properties of~$c$ and~$d$ imply

$$
Z\left(\setbox1=\hbox{\begin{picture}(0,0)%
\includegraphics{y2.eps}%
\end{picture}%
\setlength{\unitlength}{3947sp}%
\begingroup\makeatletter\ifx\SetFigFont\undefined%
\gdef\SetFigFont#1#2#3#4#5{%
  \reset@font\fontsize{#1}{#2pt}%
  \fontfamily{#3}\fontseries{#4}\fontshape{#5}%
  \selectfont}%
\fi\endgroup%
\begin{picture}(344,394)(204,-133)
\end{picture}
}\vcenter{\box1}\!\right)^{\rm op}=Z\left(\setbox1=\hbox{\begin{picture}(0,0)%
\includegraphics{y1.eps}%
\end{picture}%
\setlength{\unitlength}{3947sp}%
\begingroup\makeatletter\ifx\SetFigFont\undefined%
\gdef\SetFigFont#1#2#3#4#5{%
  \reset@font\fontsize{#1}{#2pt}%
  \fontfamily{#3}\fontseries{#4}\fontshape{#5}%
  \selectfont}%
\fi\endgroup%
\begin{picture}(344,394)(729,-133)
\end{picture}
}\vcenter{\box1}\!\right).
$$

Part~(2) of the theorem now follows from part~(1) and
the equations above.
\end{skoproof}

From now on we will fix the choices below in the definition of the 
extension of~$Z$ to unitrivalent graphs.

\begin{eqnarray}
Z\left(\setbox1=\hbox{\begin{picture}(0,0)%
\includegraphics{tyddd.eps}%
\end{picture}%
\setlength{\unitlength}{3947sp}%
\begingroup\makeatletter\ifx\SetFigFont\undefined%
\gdef\SetFigFont#1#2#3#4#5{%
  \reset@font\fontsize{#1}{#2pt}%
  \fontfamily{#3}\fontseries{#4}\fontshape{#5}%
  \selectfont}%
\fi\endgroup%
\begin{picture}(344,396)(404,116)
\end{picture}
}\vcenter{\box1}\!\right)=\Y[0,0,0]  & \qquad &
Z\left(\setbox1=\hbox{\begin{picture}(0,0)%
\includegraphics{tyuuu.eps}%
\end{picture}%
\setlength{\unitlength}{3947sp}%
\begingroup\makeatletter\ifx\SetFigFont\undefined%
\gdef\SetFigFont#1#2#3#4#5{%
  \reset@font\fontsize{#1}{#2pt}%
  \fontfamily{#3}\fontseries{#4}\fontshape{#5}%
  \selectfont}%
\fi\endgroup%
\begin{picture}(344,396)(404,116)
\end{picture}
}\vcenter{\box1}\!\right)=\Y[0,0,-1] \label{e:ZYf}\\
Z\left(\setbox1=\hbox{\begin{picture}(0,0)%
\includegraphics{tyddu.eps}%
\end{picture}%
\setlength{\unitlength}{3947sp}%
\begingroup\makeatletter\ifx\SetFigFont\undefined%
\gdef\SetFigFont#1#2#3#4#5{%
  \reset@font\fontsize{#1}{#2pt}%
  \fontfamily{#3}\fontseries{#4}\fontshape{#5}%
  \selectfont}%
\fi\endgroup%
\begin{picture}(344,396)(404,116)
\end{picture}
}\vcenter{\box1}\!\right)=\Y[-1/2,-1/2,-1]  & \qquad &
Z\left(\setbox1=\hbox{\begin{picture}(0,0)%
\includegraphics{tyuud.eps}%
\end{picture}%
\setlength{\unitlength}{3947sp}%
\begingroup\makeatletter\ifx\SetFigFont\undefined%
\gdef\SetFigFont#1#2#3#4#5{%
  \reset@font\fontsize{#1}{#2pt}%
  \fontfamily{#3}\fontseries{#4}\fontshape{#5}%
  \selectfont}%
\fi\endgroup%
\begin{picture}(344,396)(404,116)
\end{picture}
}\vcenter{\box1}\!\right)=\Y[1/2,1/2,0] \label{e:ZYb}
\end{eqnarray}

With this definition we have~$Z(G)\in\Ac(\Gamma(G))$ where~$\Gamma(G)$ is 
the underlying graph of~$G$.
For explicit computations in the following sections we list more values of~$Z$. 

\begin{equation}\label{e:Zremval}
\begin{array}{lll}
Z\left(\setbox1=\hbox{\begin{picture}(0,0)%
\includegraphics{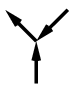}%
\end{picture}%
\setlength{\unitlength}{3947sp}%
\begingroup\makeatletter\ifx\SetFigFont\undefined%
\gdef\SetFigFont#1#2#3#4#5{%
  \reset@font\fontsize{#1}{#2pt}%
  \fontfamily{#3}\fontseries{#4}\fontshape{#5}%
  \selectfont}%
\fi\endgroup%
\begin{picture}(344,396)(404,116)
\end{picture}
}\vcenter{\box1}\!\right)=\Y[1/2,0,-1/2]  & &
Z\left(\setbox1=\hbox{\begin{picture}(0,0)%
\includegraphics{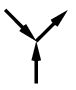}%
\end{picture}%
\setlength{\unitlength}{3947sp}%
\begingroup\makeatletter\ifx\SetFigFont\undefined%
\gdef\SetFigFont#1#2#3#4#5{%
  \reset@font\fontsize{#1}{#2pt}%
  \fontfamily{#3}\fontseries{#4}\fontshape{#5}%
  \selectfont}%
\fi\endgroup%
\begin{picture}(344,396)(404,116)
\end{picture}
}\vcenter{\box1}\!\right)=\Y[0,1/2,-1/2]\\ 
Z\left(\setbox1=\hbox{\begin{picture}(0,0)%
\includegraphics{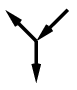}%
\end{picture}%
\setlength{\unitlength}{3947sp}%
\begingroup\makeatletter\ifx\SetFigFont\undefined%
\gdef\SetFigFont#1#2#3#4#5{%
  \reset@font\fontsize{#1}{#2pt}%
  \fontfamily{#3}\fontseries{#4}\fontshape{#5}%
  \selectfont}%
\fi\endgroup%
\begin{picture}(344,396)(404,116)
\end{picture}
}\vcenter{\box1}\!\right)=\Y[0,-1/2,-1/2]  & &
Z\left(\setbox1=\hbox{\begin{picture}(0,0)%
\includegraphics{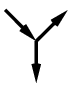}%
\end{picture}%
\setlength{\unitlength}{3947sp}%
\begingroup\makeatletter\ifx\SetFigFont\undefined%
\gdef\SetFigFont#1#2#3#4#5{%
  \reset@font\fontsize{#1}{#2pt}%
  \fontfamily{#3}\fontseries{#4}\fontshape{#5}%
  \selectfont}%
\fi\endgroup%
\begin{picture}(344,396)(404,116)
\end{picture}
}\vcenter{\box1}\!\right)=\Y[-1/2,0,-1/2]\\
Z\left(\setbox1=\hbox{\begin{picture}(0,0)%
\includegraphics{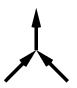}%
\end{picture}%
\setlength{\unitlength}{3947sp}%
\begingroup\makeatletter\ifx\SetFigFont\undefined%
\gdef\SetFigFont#1#2#3#4#5{%
  \reset@font\fontsize{#1}{#2pt}%
  \fontfamily{#3}\fontseries{#4}\fontshape{#5}%
  \selectfont}%
\fi\endgroup%
\begin{picture}(344,396)(404,-334)
\end{picture}
}\vcenter{\box1}\!\right)=\Aa[0,0,0]  & &
Z\left(\setbox1=\hbox{\begin{picture}(0,0)%
\includegraphics{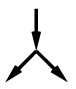}%
\end{picture}%
\setlength{\unitlength}{3947sp}%
\begingroup\makeatletter\ifx\SetFigFont\undefined%
\gdef\SetFigFont#1#2#3#4#5{%
  \reset@font\fontsize{#1}{#2pt}%
  \fontfamily{#3}\fontseries{#4}\fontshape{#5}%
  \selectfont}%
\fi\endgroup%
\begin{picture}(344,396)(404,-334)
\end{picture}
}\vcenter{\box1}\!\right)=\Aa[-1,0,0]\\
Z\left(\setbox1=\hbox{\begin{picture}(0,0)%
\includegraphics{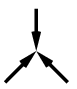}%
\end{picture}%
\setlength{\unitlength}{3947sp}%
\begingroup\makeatletter\ifx\SetFigFont\undefined%
\gdef\SetFigFont#1#2#3#4#5{%
  \reset@font\fontsize{#1}{#2pt}%
  \fontfamily{#3}\fontseries{#4}\fontshape{#5}%
  \selectfont}%
\fi\endgroup%
\begin{picture}(344,396)(404,-334)
\end{picture}
}\vcenter{\box1}\!\right)=\Aa[-1,-1/2,-1/2]  & &
Z\left(\setbox1=\hbox{\begin{picture}(0,0)%
\includegraphics{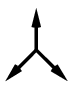}%
\end{picture}%
\setlength{\unitlength}{3947sp}%
\begingroup\makeatletter\ifx\SetFigFont\undefined%
\gdef\SetFigFont#1#2#3#4#5{%
  \reset@font\fontsize{#1}{#2pt}%
  \fontfamily{#3}\fontseries{#4}\fontshape{#5}%
  \selectfont}%
\fi\endgroup%
\begin{picture}(344,396)(404,-334)
\end{picture}
}\vcenter{\box1}\!\right)=\Aa[0,1/2,1/2]\\
Z\left(\setbox1=\hbox{\begin{picture}(0,0)%
\includegraphics{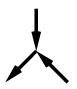}%
\end{picture}%
\setlength{\unitlength}{3947sp}%
\begingroup\makeatletter\ifx\SetFigFont\undefined%
\gdef\SetFigFont#1#2#3#4#5{%
  \reset@font\fontsize{#1}{#2pt}%
  \fontfamily{#3}\fontseries{#4}\fontshape{#5}%
  \selectfont}%
\fi\endgroup%
\begin{picture}(344,396)(404,-334)
\end{picture}
}\vcenter{\box1}\!\right)=\Aa[-1/2,1/2,0]  & &
Z\left(\setbox1=\hbox{\begin{picture}(0,0)%
\includegraphics{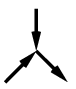}%
\end{picture}%
\setlength{\unitlength}{3947sp}%
\begingroup\makeatletter\ifx\SetFigFont\undefined%
\gdef\SetFigFont#1#2#3#4#5{%
  \reset@font\fontsize{#1}{#2pt}%
  \fontfamily{#3}\fontseries{#4}\fontshape{#5}%
  \selectfont}%
\fi\endgroup%
\begin{picture}(344,396)(404,-334)
\end{picture}
}\vcenter{\box1}\!\right)=\Aa[-1/2,0,1/2]\\
Z\left(\setbox1=\hbox{\begin{picture}(0,0)%
\includegraphics{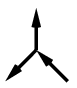}%
\end{picture}%
\setlength{\unitlength}{3947sp}%
\begingroup\makeatletter\ifx\SetFigFont\undefined%
\gdef\SetFigFont#1#2#3#4#5{%
  \reset@font\fontsize{#1}{#2pt}%
  \fontfamily{#3}\fontseries{#4}\fontshape{#5}%
  \selectfont}%
\fi\endgroup%
\begin{picture}(344,396)(404,-334)
\end{picture}
}\vcenter{\box1}\!\right)=\Aa[-1/2,0,-1/2]  & &
Z\left(\setbox1=\hbox{\begin{picture}(0,0)%
\includegraphics{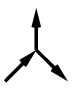}%
\end{picture}%
\setlength{\unitlength}{3947sp}%
\begingroup\makeatletter\ifx\SetFigFont\undefined%
\gdef\SetFigFont#1#2#3#4#5{%
  \reset@font\fontsize{#1}{#2pt}%
  \fontfamily{#3}\fontseries{#4}\fontshape{#5}%
  \selectfont}%
\fi\endgroup%
\begin{picture}(344,396)(404,-334)
\end{picture}
}\vcenter{\box1}\!\right)=\Aa[-1/2,-1/2,0]
\end{array}
\end{equation}

The computation of these values can be simplified by first 
generalizing symmetry properties
of~$Z$ from $q$-tangles (Proposition~3.1 of~\cite{LM2}) to trivalent graphs.

\section{The Alexander series of a tetrahedron}\label{s:tetra1}

For a trivalent diagram~$D$ 
on an admissibly colored trivalent graph~$\Gamma$ we
define

\begin{equation}\label{e:WstWt}
\Wst(D)=\Wnabt(D)\,h^{(1/2)\# V(D)-1}\in
h^{-1}\Q[[h]]\subset\Q[[h]][h^{-1}].
\end{equation}

This definition induces a continous linear map
$\Wst:\Anabc^0(\Gamma)\longrightarrow h^{-1}\Q[[h]]$.
Notice that~$(1/2)\#V(D)>\deg(D)$ when~$\Gamma\subset D$ has trivalent vertices.
For an admissibly colored framed unitrivalent graph~$G\subset \R^2\times I$ we define
$\Zb(G)=\tau(Z(G))$ where
the continous linear map~$\tau:\Ac(\Gamma(G))\longrightarrow \Anabc^0(\Gamma(G))$ is 
defined on trivalent diagrams~$D$ by~$\tau(D)=D$ and the skeleton~$\Gamma(G)$ of~$D$
is colored according to~$G$.

\begin{defi}
The invariant $\nabh(G)=\Wst\left(\Zb(G)\right)$ of an admissibly colored
framed trivalent graph~$G$ is called the {\em Alexander
series} of~$G$.
\end{defi}

In the following three sections we will compute the Alexander
series of a trivially embedded colored tetrahedron~$T$
in two different ways. With the first computation we will
determine the value of the Alexander series~$\nabh(T)$, and
with the second computation we will use this value to prove
Theorem~\ref{t:F}.

Let $\setbox1=\hbox{\begin{picture}(0,0)%
\includegraphics{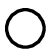}%
\end{picture}%
\setlength{\unitlength}{3947sp}%
\begingroup\makeatletter\ifx\SetFigFont\undefined%
\gdef\SetFigFont#1#2#3#4#5{%
  \reset@font\fontsize{#1}{#2pt}%
  \fontfamily{#3}\fontseries{#4}\fontshape{#5}%
  \selectfont}%
\fi\endgroup%
\begin{picture}(210,209)(651,-41)
\put(856,-41){\makebox(0,0)[lb]{\smash{\SetFigFont{6}{7.2}{\familydefault}{\mddefault}{\updefault}{\color[rgb]{0,0,0}$c$}%
}}}
\end{picture}
}\vcenter{\box1}\!$ be the trivial knot with color~$c=(\lambda,\mu)$ and with $0$-framing. 
It follows from
Section~2.2 of~\cite{BNG} that for $c=(1,0)$ we have
$\nabh(\setbox1=\hbox{\begin{picture}(0,0)%
\includegraphics{oc.eps}%
\end{picture}%
\setlength{\unitlength}{3947sp}%
\begingroup\makeatletter\ifx\SetFigFont\undefined%
\gdef\SetFigFont#1#2#3#4#5{%
  \reset@font\fontsize{#1}{#2pt}%
  \fontfamily{#3}\fontseries{#4}\fontshape{#5}%
  \selectfont}%
\fi\endgroup%
\begin{picture}(210,209)(651,-41)
\put(856,-41){\makebox(0,0)[lb]{\smash{\SetFigFont{6}{7.2}{\familydefault}{\mddefault}{\updefault}{\color[rgb]{0,0,0}$c$}%
}}}
\end{picture}
}\vcenter{\box1}\!)=1/(e^{h/2}-e^{-h/2})$. For~$c=(\lambda,\mu)$
we obtain from this value the more general formula

\begin{equation}\label{e:nabhOc}
\nabh\left(\setbox1=\hbox{\begin{picture}(0,0)%
\includegraphics{oc.eps}%
\end{picture}%
\setlength{\unitlength}{3947sp}%
\begingroup\makeatletter\ifx\SetFigFont\undefined%
\gdef\SetFigFont#1#2#3#4#5{%
  \reset@font\fontsize{#1}{#2pt}%
  \fontfamily{#3}\fontseries{#4}\fontshape{#5}%
  \selectfont}%
\fi\endgroup%
\begin{picture}(210,209)(651,-41)
\put(856,-41){\makebox(0,0)[lb]{\smash{\SetFigFont{6}{7.2}{\familydefault}{\mddefault}{\updefault}{\color[rgb]{0,0,0}$c$}%
}}}
\end{picture}
}\vcenter{\box1}\!\right)=\frac{1}{e^{\lambda
h/2}-e^{-\lambda
h/2}}.
\end{equation}

by using that the 
left side of equation~(\ref{e:Wprop1}) suffices to compute
$\Wst(\Zb(\setbox1=\hbox{\begin{picture}(0,0)%
\includegraphics{oc.eps}%
\end{picture}%
\setlength{\unitlength}{3947sp}%
\begingroup\makeatletter\ifx\SetFigFont\undefined%
\gdef\SetFigFont#1#2#3#4#5{%
  \reset@font\fontsize{#1}{#2pt}%
  \fontfamily{#3}\fontseries{#4}\fontshape{#5}%
  \selectfont}%
\fi\endgroup%
\begin{picture}(210,209)(651,-41)
\put(856,-41){\makebox(0,0)[lb]{\smash{\SetFigFont{6}{7.2}{\familydefault}{\mddefault}{\updefault}{\color[rgb]{0,0,0}$c$}%
}}}
\end{picture}
}\vcenter{\box1}\!))$ because~$\setbox1=\hbox{\begin{picture}(0,0)%
\includegraphics{oc.eps}%
\end{picture}%
\setlength{\unitlength}{3947sp}%
\begingroup\makeatletter\ifx\SetFigFont\undefined%
\gdef\SetFigFont#1#2#3#4#5{%
  \reset@font\fontsize{#1}{#2pt}%
  \fontfamily{#3}\fontseries{#4}\fontshape{#5}%
  \selectfont}%
\fi\endgroup%
\begin{picture}(210,209)(651,-41)
\put(856,-41){\makebox(0,0)[lb]{\smash{\SetFigFont{6}{7.2}{\familydefault}{\mddefault}{\updefault}{\color[rgb]{0,0,0}$c$}%
}}}
\end{picture}
}\vcenter{\box1}\!$ is $0$-framed and we know the structure
of~$\Anab(1)$
(see~\cite{Thu} for a more general result).
For a trivalent diagram~$D$ on~$\Gamma\in\Hom_{{\cal C}_0}(c,c')$ the 
definition

\begin{equation}
\Ws(D)=W(D)h^{\deg D}\in \Hom_{\gl(1\vert 1)}(W(c), W(c'))[[h]]
\end{equation}
 
induces a continous linear map~$\Ws:\Anabc^0(\Gamma)\longrightarrow
\Hom_{\gl(1\vert 1)}(W(c), W(c'))[[h]]$. Recall the 
definition~$\varphi(x)=(e^{x/2}-e^{-x/2})/x=2\sinh(x/2)/x$ from the introduction.
Equation~(\ref{e:nabhOc}) and Lemma~\ref{l:cutnotimp1} imply that
for the trivalent diagram~$\nu$ on a skeleton~$\I_t\in 
\End_{\Anabc_0}(t)$ with~$t=(\lambda,\mu,\sigma)$ we have

\begin{equation}\label{e:Wnu}
\Ws(\nu)=(1/\varphi(\lambda h))\,\I_t.
\end{equation}

Let $G\in\End_{{\cal G}^\na}(+)$. 
Then~$H=\setbox1=\hbox{\begin{picture}(0,0)%
\includegraphics{rightcap.eps}%
\end{picture}%
\setlength{\unitlength}{3947sp}%
\begingroup\makeatletter\ifx\SetFigFont\undefined%
\gdef\SetFigFont#1#2#3#4#5{%
  \reset@font\fontsize{#1}{#2pt}%
  \fontfamily{#3}\fontseries{#4}\fontshape{#5}%
  \selectfont}%
\fi\endgroup%
\begin{picture}(359,263)(579,217)
\end{picture}
}\vcenter{\box1}\!\circ(\id_-\otimes G)\circ\setbox1=\hbox{\begin{picture}(0,0)%
\includegraphics{leftcup.eps}%
\end{picture}%
\setlength{\unitlength}{3947sp}%
\begingroup\makeatletter\ifx\SetFigFont\undefined%
\gdef\SetFigFont#1#2#3#4#5{%
  \reset@font\fontsize{#1}{#2pt}%
  \fontfamily{#3}\fontseries{#4}\fontshape{#5}%
  \selectfont}%
\fi\endgroup%
\begin{picture}(359,263)(564,298)
\end{picture}
}\vcenter{\box1}\!
\in\End_{{\cal G}^\na}(\emptyset)$ is called
the {\em closure} (or {\em trace}) of~$G$. 
When~$H$ is admissibly colored and the upper (or lower) edge of~$G$ 
has color~$(\lambda,\mu)$ 
then the definition of~$Z$, Lemma~\ref{l:cutnotimp1}, 
and equations~(\ref{e:WstWt}) and~(\ref{e:Wnu}) imply that for~$t=(\lambda,\mu,0)$
we have

\begin{equation}\label{e:Wclosure}
\nabh(H)\I_t=
\frac{h^{(1/2)\#V(\Gamma(H))}}{\lambda h\varphi(\lambda h)} \Ws\left(\Zb(G)\right).
\end{equation}

For $i=1,2,3$ let~$t_i=(\lambda_i, \mu_i,\sigma_i)\in\Q^*\times\Q\times\Z/(2)$ with

\begin{equation}
\lambda_1+\lambda_2\ ,\ \lambda_2+\lambda_3\ ,\
\lambda_1+\lambda_2+\lambda_3\ \in\ \Q^*.\label{e:lambdai}
\end{equation}

We will use in the following sections the triples
$t_i=(\lambda_i,\mu_i,\sigma_i)\in\Q^*\times \Q\times\Z/(2)$ ($i=4,\ldots,11$) 
defined by

\begin{equation}
t_4=t_1+t_2+e_2\ ,\ t_5=t_4+t_3+e_2\ ,\ t_6=t_2+t_3+e_2\ ,
t_{i}=t_{i-3}-2e_2+e_3\ (i\geq 7)\label{e:ti}
\end{equation}

and the colors~$c_i=((-1)^{\sigma_i}\lambda_i, (-1)^{\sigma_i}\mu_i)$ ($i=1,\ldots, 11$).
Equation~(\ref{e:tensordecomp}) implies that the tensor product of
three modules decomposes as

\begin{equation}\label{e:threedecomp}
V_{t_1}\otimes V_{t_2}\otimes V_{t_3}\cong V_{t_5}\oplus
V_{t_8}^{\oplus 2}\oplus V_{t_{11}}.
\end{equation}

The following lemma concerns simple modules of multiplicity~$1$ in equation~(\ref{e:threedecomp}).

\begin{lemma}\label{l:assotriv}
With $t_i$ as in equation~(\ref{e:ti}) the following formulas hold.

\begin{eqnarray*}
(\I_{t_1}\otimes \Ystd_{t_2,t_3})\circ\Ystd_{t_1,t_6} & = & 
(\Ystd_{t_1,t_2}\otimes \I_{t_3})\circ\Ystd_{t_4,t_3},\\
(\I_{t_1}\otimes \Ydot_{t_2,t_3})\circ\Ydot_{t_1,t_9} & = & 
(\Ydot_{t_1,t_2}\otimes \I_{t_3})\circ\Ydot_{t_7,t_3},\\
\Astd_{t_1,t_6}\circ(\I_{t_1}\otimes \Astd_{t_2,t_3}) & = & 
\Astd_{t_4, t_3}\circ (\Astd_{t_1,t_2}\otimes \I_{t_3}),\\
\Adot_{t_1, t_9}\circ (\I_{t_1}\otimes \Adot_{t_2,t_3})& = & 
\Adot_{t_7, t_3}\circ (\Ydot_{t_1,t_2}\otimes \I_{t_3}).
\end{eqnarray*}
\end{lemma}
\begin{proof}
We compute

$$
(\I_{t_1}\otimes \Ystd_{t_2,t_3})\circ\Ystd_{t_1,t_6}(v_{t_5}) = v_{t_1}\otimes 
v_{t_2}\otimes v_{t_3} =
(\Ystd_{t_1,t_2}\otimes \I_{t_3})\circ\Ystd_{t_4,t_3}(t_{t_5}).
$$

This implies the first equation. The second equation is proved similarly. 
The remaining two equations follow from equation~(\ref{e:threedecomp}) and  
similar computations.
\end{proof}

\begin{figure}[!h]
\centering 
$$
T_1\ =\quad
\mbox{$\setbox1=\hbox{\input{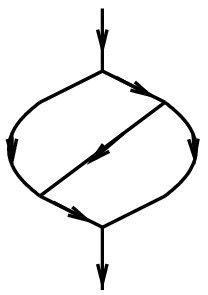}}\vcenter{\box1}\!$} \quad , \quad
T_2\ =\quad
\mbox{$\setbox1=\hbox{\input{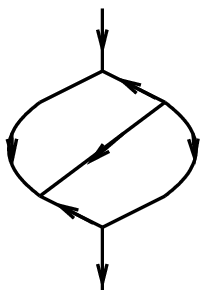}}\vcenter{\box1}\!$} \quad , \quad
T_3\ =\quad
\mbox{$\setbox1=\hbox{\input{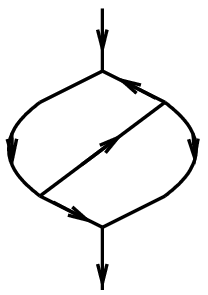}}\vcenter{\box1}\!$} \quad .
$$
\caption{Three colored diagrams whose
closures are planar tetrahedra} \label{f:threetetra}
\end{figure}

Consider the three
diagrams~$T_\indx$ ($\indx=1,2,3$) shown in Figure~\ref{f:threetetra}.

Fix a choice of~$\indx\in\{1,2,3\}$. 
Define~$\sigma_i=\sigma_i(\indx)$ by~$\sigma_i=1$ iff the edge 
of~$T_\indx$ labeled by~$c_i$ points downwards. Then 
equations~(\ref{e:lambdai}) and~(\ref{e:ti}) ensure that the 
coloring of~$T_\indx$ is admissible. 
Let~$b_i=((-1)^{\sigma_i}\lambda_i, (-1)^{\sigma_i}\mu_i, 
\sigma_i)$. Define~$b=b_5$ if $\indx=1$, $b=b_{11}$ if $k=2$, and
$b=b_8$ if $k=3$.
Then $T_\indx=U_\indx\circ A_\indx\circ L_\indx$ where 
$U_\indx\in\Hom_\G((b_1(b_2b_3)),b)$ consists of the upper half 
of~$T_\indx$, $L_\indx\in\Hom_\G(b,((b_1b_2)b_3))$ consists of the 
lower half of~$T_\indx$, 
and~$A_\indx\in\Hom_\G(((b_1b_2)b_3),(b_1(b_2b_3)))$ consists of 
three vertical colored strands. 
It follows from equations~(\ref{e:ZYf}), (\ref{e:WoYo}), and Lemma~\ref{l:assotriv}) that
the graph~$L_1$ is mapped by~$\Ws\circ \Zb$ to the
following morphism from~$V_{t_5}[[h]]$ to~$(V_{t_1}\otimes
V_{t_2}\otimes V_{t_3})[[h]]$:

\begin{equation}\label{e:WZL}
\Ws\left(\Zb(L_1)\right)=\Wnab(\Gamma(L_1))
=(\Ystd_{t_1, t_2}\otimes \I_{t_3})\circ
\Ystd_{t_4, t_3}
=(\I_{t_1}\otimes \Ystd_{t_2, t_3})\circ
\Ystd_{t_1, t_6}.
\end{equation}

The graph~$U_1$ is mapped by
$\Ws\circ \Zb$ to the following morphism from~$(V_{t_1}\otimes
V_{t_2}\otimes V_{t_3})[[h]]$ to~$V_{t_5}[[h]]$ (see 
equations~(\ref{e:Zremval}),~(\ref{e:Wnu}),~(\ref{e:WoAo}))~:

\begin{eqnarray}
\Ws\left(\Zb(U_1)\right)& = &
\varphi(\lambda_5 h)\varphi(\lambda_6 h)
\Wnab(\Gamma(U_1))\nonumber\\
& = &
\lambda_5\lambda_6\varphi(\lambda_5 h)\varphi(\lambda_6h)
\Astd_{t_1, t_6}\circ(\I_{t_1}\otimes\Astd_{t_2, t_3}).\label{e:WZU}
\end{eqnarray}


The associator~$\Phi$ is a series with constant
term~$1$, and all higher order terms of~$\Phi$ involve a
commutator (see (DA5)). Therefore, by equation~(\ref{e:threedecomp}) and
Schurs lemma, the action of~$\Phi$ on~$(V_{t_1}\otimes
V_{t_2}\otimes V_{t_3})[[h]]$ restricts to

\begin{equation}\label{e:Phitriv}
\I_{t_5}\in\End_\gloo(V_{t_5})\subset\End_\gloo(V_{t_1}\otimes
V_{t_2}\otimes V_{t_3})[[h]].
\end{equation}

This implies

\begin{equation}\label{e:Phitrivb}
\Astd_{t_1, t_6}\circ(\I_{t_1}\otimes \Astd_{t_2, t_3})
\circ \Ws(\Zb(A_1))\circ(\I_{t_1}\circ \Ystd_{t_2, t_3})
\circ
\Ystd_{t_1, t_6}=\I_{t_5}.
\end{equation}

Let $S_1$ be the closure of~$T_1$. Equations~(\ref{e:Wclosure}), 
(\ref{e:WZL}), (\ref{e:WZU}), 
and~(\ref{e:Phitrivb}) allow to
compute $\nabh(S_1)$ without knowing~$\Phi$~:

\begin{equation}\label{e:WtZT}
\nabh(S_1)\I_{t_5}=
\frac{h}{\lambda_5\varphi(\lambda_5 h)}\Ws\left(\Zb(T_1)\right)
 =  \lambda_6 h\varphi(\lambda_6 h)\I_{t_5}=(e^{\lambda_6h/2}-e^{-\lambda_6h/2})\I_{t_5}.
\end{equation}

%

In general, the Alexander series of
a planar tetrahedron is given by the following lemma.

\begin{lemma}\label{l:nabhT}
 Let~$T$ be a planar tetrahedron with admissible
coloring and blackboard framing. There exists a unique
edge~$e$ of~$T$ such that by reversing the orientation of~$e$ we
obtain a tetrahedron with one source and without a sink.
Let~$(\lambda,\mu)$ be the color of~$e$. Then we
have

$$
\nabh(T)=e^{\lambda h/2}-e^{-\lambda h/2}.$$
 \end{lemma}
\begin{proof} By Lemma~7.2.A of~\cite{Vir} there are
four isotopy classes of planar oriented tetrahedra with blackboard framing.
By Lemma~7.2.B of~\cite{Vir} two of these tetrahedra
do not have an admissible coloring.
By Lemma~7.2.C of~\cite{Vir} the remaining two oriented tetrahedra 
have a unique edge~$e$ as in the lemma.
For one admissibly colored tetrahedron~$S_1$ we have 
computed~$\nabh(S_1)$ in equation~(\ref{e:WtZT}).
The second tetrahedron is the closure~$S_2$ of~$T_2$ (see Figure~\ref{f:threetetra}). The
computation of~$\nabh(S_2)$ 
proceeds along the same lines as the computation of~$\nabh(S_1)$.
\end{proof}


\section{Associativity and $\gloo$-modules}\label{s:lemmas}

For the diagram~$T_3$ in Figure~\ref{f:threetetra}
we will see in Section~\ref{s:tetra2} that the contribution
of the associator~$\Phi$ in the computation of~$\Ws(\Zb(T_3))$ is non-trivial.
We will use
Lemma~\ref{l:assobas} and Corollary~\ref{c:assobas} 
below for a similar
purpose in this computation as we used Lemma~\ref{l:assotriv} to
deduce equation~(\ref{e:WZL}). 

\begin{lemma}\label{l:assobas}
For $t_i=(\lambda_i, \mu_i, \sigma_i)$ as in equations~(\ref{e:lambdai}), (\ref{e:ti}) and
$\kappa=\frac{\lambda_2(\lambda_1+\lambda_2+\lambda_3)}
{(\lambda_1+\lambda_2)(\lambda_2+\lambda_3)}$ we
have two bases 
of~$\Hom_{\gloo}(V_{t_8},\Vo\otimes\Vp\otimes\Vpp)$
that are related by

\begin{equation}\label{e:assobas} \left(\begin{array}{l} (\I_{t_1}\otimes \Ydot_{t_2,
t_3})\circ\Ystd_{t_1,t_9}\\ (\I_{t_1}\otimes \Ystd_{t_2,
t_3})\circ\Ydot_{t_1,t_6}\end{array}\right)=
\left(\begin{array}{cc}
\frac{(-1)^{\sigma_2}\lambda_3}{-\lambda_2-\lambda_3} &
\kappa\\
1 & \frac{(-1)^{\sigma_2}\lambda_1}{\lambda_1+\lambda_2}
\end{array}\right)
\left(\begin{array}{l} (\Ydot_{t_1, t_2}\otimes
\I_{t_3})\circ\Ystd_{t_7,t_3}\\ (\Ystd_{t_1,t_2}\otimes \I_{t_3}
)\circ\Ydot_{t_4, t_3}\end{array}\right). \end{equation}
\end{lemma}
\begin{proof}
We compute

\begin{eqnarray*}
& & (\I_{t_1}\otimes\Ydot_{t_2, t_3})(\Ystd_{t_1, t_9}(w_{t_8})) =
(\I_{t_1}\otimes\Ydot_{t_2,t_3})(F\cdot(v_{t_1}\otimes v_{t_9}))\\
& & = w_{t_1}\otimes\Ydot_{t_2, t_3}(v_{t_9})+(-1)^{\sigma_1}
v_{t_1}\otimes\Ydot_{t_2, t_3}(w_{t_9})\\ & &
=\frac{1}{\lambda_2+\lambda_3} w_{t_1}\otimes E\cdot(
w_{t_2}\otimes w_{t_3})+(-1)^{\sigma_1} v_{t_1}\otimes
w_{t_2}\otimes w_{t_3}\\ & &
=\frac{\lambda_2}{\lambda_2+\lambda_3} w_{t_1}\otimes
v_{t_2}\otimes
w_{t_3}-\frac{(-1)^{\sigma_2}\lambda_3}{\lambda_2+\lambda_3}
w_{t_1}\otimes w_{t_2}\otimes v_{t_3}+(-1)^{\sigma_1}
v_{t_1}\otimes w_{t_2}\otimes w_{t_3}
\end{eqnarray*}

and similarly (or simpler)

\begin{eqnarray*}
(\I_{t_1}\otimes\Ystd_{t_2, t_3})(\Ydot_{t_1, t_6}(w_{t_8})) & = &
(-1)^{\sigma_2}w_{t_1}\otimes v_{t_2}\otimes
w_{t_3}+w_{t_1}\otimes w_{t_2}\otimes v_{t_3},\\
(\Ydot_{t_1,t_2}\otimes \I_{t_3})(\Ystd_{t_7, t_3}(w_{t_8})) & = &
\frac{(-1)^{\sigma_2}\lambda_2}{\lambda_1+\lambda_2}
w_{t_1}\otimes v_{t_2}\otimes w_{t_3}+w_{t_1}\otimes
w_{t_2}\otimes v_{t_3}\\\nopagebreak{} & &
-\frac{(-1)^{\sigma_1+\sigma_2}\lambda_1}{\lambda_1+\lambda_2}v_{t_1}\otimes
w_{t_2}\otimes w_{t_3},\\ (\Ystd_{t_1,t_2}\otimes
\I_{t_3})(\Ydot_{t_4, t_3}(w_{t_8})) & = & w_{t_1}\otimes
v_{t_2}\otimes w_{t_3}+(-1)^{\sigma_1} v_{t_1}\otimes
w_{t_2}\otimes w_{t_3}.
\end{eqnarray*}

We see that the vectors
$(\I_{t_1}\otimes\Ydot_{t_2,t_3})(\Ystd_{t_1,t_9}(w_{t_8}))$ and
$(\I_{t_1}\otimes\Ystd_{t_2,t_3})(\Ydot_{t_1, t_6}(w_{t_8}))$ as well 
as~$(\Ydot_{t_1, t_2}\otimes
\I_{t_3})(\Ystd_{t_7,t_3}(w_{t_8}))$ and~$(\Ystd_{t_1, t_2}\otimes
\I_{t_3})(\Ydot_{t_4,t_3}(w_{t_8}))$ are linearly independent.
Equation~(\ref{e:threedecomp}) implies
that~$\Hom_{\gloo}(V_{t_8},\Vo\otimes\Vp\otimes\Vpp)$ is
two-dimensional, so we have
found two bases of that space. Verify that equation~(\ref{e:assobas}) 
holds when we evaluate the morphisms in this equation on~$w_{t_8}$. This
implies the lemma because~$w_{t_8}$ generates $V_{t_8}$.
\end{proof}

Dually to Lemma~\ref{l:assobas} we have the following corollary.

\begin{coro}\label{c:assobas}
With $t_i$ and~$\kappa$ as in Lemma~\ref{l:assobas} two bases 
of
$$\Hom_{\gloo}(\Vo\otimes\Vp\otimes\Vpp,
V_{t_8})$$
are related by

$$ \left(\begin{array}{l}  \Astd_{t_1, t_9}\circ(\I_{t_1}\otimes
\Adot_{t_2,t_3})\\ \Adot_{t_1, t_6}\circ(\I_{t_1}\otimes
\Astd_{t_2,t_3})\end{array}\right)=\left(\begin{array}{cc}
\frac{(-1)^{\sigma_2}\lambda_1}{-\lambda_1-\lambda_2} & 1\\
\kappa
& \frac{(-1)^{\sigma_2}\lambda_3}{\lambda_2+\lambda_3}
\end{array}\right)\left(\begin{array}{l} \Astd_{t_7,t_3}\circ(\Adot_{t_1,t_2}\otimes \I_{t_3})\\
\Adot_{t_4,t_3}\circ(\Astd_{t_1,t_2}\otimes
\I_{t_3})\end{array}\right). $$
\end{coro}
\begin{proof}
By Lemma~\ref{l:assobas} linear maps~$f_i$, $g_i$ 
are related by
$(f_1,
f_2)^T=A(g_1, g_2)^T$ for a certain matrix $A$. 
Therefore, maps~$f_i'$, $g_i'$ satisfying
$f_i'\circ f_j=\delta_{ij} \I_{t_8}$
and $g_i'\circ g_j=\delta_{ij} \I_{t_8}$ are related by
$(f_1', f_2')^T=({A^T})^{-1}(g_1', g_2')^T$. This implies the
corollary.
\end{proof}

In the rest of this section
we prepare the analogue of equation~(\ref{e:Phitrivb}) 
for the computation in Section~\ref{s:tetra2}. 
By the definition of the~$\gloo$-module structure
of~$V_{t_1}\otimes V_{t_2}$ we have

\begin{eqnarray}
F\cdot(v_{t_1}\otimes v_{t_2}) & = & w_{t_1}\otimes
v_{t_2}+(-1)^{\sigma_1} v_{t_1}\otimes
w_{t_2}\quad\mbox{and}\label{e:bc1}\\ E\cdot(w_{t_1}\otimes
w_{t_2}) & = & -(-1)^{\sigma_1}\lambda_2 w_{t_1}\otimes
v_{t_2}+\lambda_1 v_{t_1}\otimes w_{t_2}.\label{e:bc2}
\end{eqnarray}

For~$\lambda_1\not=-\lambda_2$
equations~(\ref{e:bc1}) and~(\ref{e:bc2}) are formulas for a
change of bases in the two dimensional eigenspace of~$H$
on~$\Vo\otimes \Vp$ and imply

\begin{eqnarray}
(\lambda_1+\lambda_2)w_{t_1}\otimes v_{t_2} & = & \lambda_1
F\cdot(v_{t_1}\otimes v_{t_2})-(-1)^{\sigma_1}
E\cdot(w_{t_1}\otimes w_{t_2}),\label{e:wvFE}\\
(\lambda_1+\lambda_2)v_{t_1}\otimes w_{t_2} & = &
(-1)^{\sigma_1}\lambda_2 F\cdot(v_{t_1}\otimes v_{t_2})+
E\cdot(w_{t_1}\otimes w_{t_2}).
\end{eqnarray}

\begin{lemma}\label{l:HI}
With $t_i$ and~$\kappa$ as in Lemma~\ref{l:assobas} the following
formulas hold~:

\begin{eqnarray*}
(\Astd_{t_1,t_2}\otimes \I_{t_3})\circ(\I_{t_1}\otimes
\Ystd_{t_2,t_3}) & = &
\Ystd_{t_4,t_3}\circ\Astd_{t_1,t_6}+\frac{(-1)^{\sigma_2}\lambda_1}{\lambda_1+\lambda_2}
\Ydot_{t_4,t_3}\circ\Adot_{t_1,t_6},\\ (\Adot_{t_1,t_2}\otimes
\I_{t_3})\circ(\I_{t_1}\otimes \Ydot_{t_2, t_3}) & = &
\Ydot_{t_7,t_3}\circ\Adot_{t_1,t_9}-\frac{(-1)^{\sigma_2}\lambda_3}{\lambda_2+\lambda_3}
\Ystd_{t_7,t_3}\circ\Astd_{t_1,t_9},\\ (\Adot_{t_1,t_2}\otimes
\I_{t_3} )\circ(\I_{t_1}\otimes\Ystd_{t_2, t_3}) & = &
\Ystd_{t_7,t_3}\circ\Adot_{t_1,t_6},
\\
(\Astd_{t_1,t_2}\otimes \I_{t_3})\circ(\I_{t_1}\otimes\Ydot_{t_2,
t_3}) & = & \kappa\Ydot_{t_4,t_3}\circ\Astd_{t_1,t_9}.
\end{eqnarray*}
\end{lemma}
\begin{proof}
We have

\begin{eqnarray*} & & (\Astd_{t_1,t_2}\otimes \I_{t_3})\circ(\I_{t_1}\otimes
\Ystd_{t_2,t_3})(v_{t_1}\otimes v_{t_6})\\ & &
=(\Astd_{t_1,t_2}\otimes \I_{t_3})(v_{t_1}\otimes v_{t_2}\otimes
v_{t_3})=v_{t_4}\otimes v_{t_3}\\ & & =
\left(\Ystd_{t_4,t_3}\circ\Astd_{t_1,t_6}+\frac{(-1)^{\sigma_2}\lambda_1}{\lambda_1+\lambda_2}
\Ydot_{t_4,t_3}\circ\Adot_{t_1,t_6}\right)(v_{t_1}\otimes
v_{t_6}).
\end{eqnarray*}

Using equation~(\ref{e:wvFE}) we see that

\begin{eqnarray*} & & (\Astd_{t_1,t_2}\otimes \I_{t_3})\circ(\I_{t_1}\otimes
\Ystd_{t_2,t_3})(w_{t_1}\otimes w_{t_6})\\ & &
=(\Astd_{t_1,t_2}\otimes \I_{t_3})(w_{t_1}\otimes
F\cdot(v_{t_2}\otimes v_{t_3}))\\ & & =(\Astd_{t_1,t_2}\otimes
\I_{t_3})(w_{t_1}\otimes w_{t_2}\otimes
v_{t_3}+(-1)^{\sigma_2}w_{t_1}\otimes v_{t_2}\otimes w_{t_3})\\
& & =(\Astd_{t_1,t_2}\otimes
\I_{t_3})\left(\frac{(-1)^{\sigma_2}\lambda_1}{\lambda_1+\lambda_2}
F\!\cdot\!(v_{t_1}\otimes v_{t_2})\otimes w_{t_3}
-\frac{(-1)^{\sigma_1+\sigma_2}}{\lambda_1+\lambda_2}
E\!\cdot\!(w_{t_1}\otimes w_{t_2})\otimes w_{t_3}\right)\\ & &
=\frac{(-1)^{\sigma_2}\lambda_1}{\lambda_1+\lambda_2}
w_{t_4}\otimes w_{t_3}\\ & & =
\left(\Ystd_{t_4,t_3}\circ\Astd_{t_1,t_6}+\frac{(-1)^{\sigma_2}\lambda_1}{\lambda_1+\lambda_2}
\Ydot_{t_4,t_3}\circ\Adot_{t_1,t_6}\right)(w_{t_1}\otimes
w_{t_6}).
\end{eqnarray*}

Since an element of~$\Hom_\gloo(V_{t_1}\otimes V_{t_6},
V_{t_4}\otimes V_{t_3})$ is determined by the images of
$v_{t_1}\otimes v_{t_6}$ and $w_{t_1}\otimes w_{t_6}$ this implies
the first equation of the lemma. The remaining three equations are
proved similarly.
\end{proof}

The following corollary holds for reasons of symmetry.

\begin{coro}\label{c:HI}
With $t_i$ and~$\kappa$ as in Lemma~\ref{l:assobas} the following
formulas hold~:

\begin{eqnarray*}
(\I_{t_1}\otimes \Astd_{t_2,t_3})\circ(\Ystd_{t_1,t_2}\otimes
\I_{t_3}) & = &
\Ystd_{t_1,t_6}\circ\Astd_{t_4,t_3}+\frac{(-1)^{\sigma_2}\lambda_3}{\lambda_2+\lambda_3}
\Ydot_{t_1,t_6}\circ\Adot_{t_4,t_3},\\ (\I_{t_1}\otimes
\Adot_{t_2,t_3})\circ(\Ydot_{t_1,t_2}\otimes \I_{t_3}) & = &
\Ydot_{t_1,t_9}\circ\Adot_{t_7,t_3}-\frac{(-1)^{\sigma_2}\lambda_1}{\lambda_1+\lambda_2}
\Ystd_{t_1,t_9}\circ\Astd_{t_7,t_3},\\
(\I_{t_1}\otimes\Adot_{t_2,t_3})\circ(\Ystd_{t_1,t_2}\otimes
\I_{t_3}) & = & \Ystd_{t_1,t_9}\circ\Adot_{t_4,t_3},\\
(\I_{t_1}\otimes\Astd_{t_2,t_3})\circ(\Ydot_{t_1,t_2}\otimes
\I_{t_3}) & = & \kappa\Ydot_{t_1,t_6}\circ\Astd_{t_7,t_3}.
\end{eqnarray*}
\end{coro}
\begin{proof}
Let~$\tau_{V,W}\in\Hom(V\otimes W, W\otimes V)$ be the linear map induced by
the permutation of tensor factors. When we interchange the labels~$t_{3n-2}$ and~$t_{3n}$ 
($n=1,2,3$) in the
equations of Lemma~\ref{l:HI}, 
replace the equations~$X=Y\in\Hom(V\otimes
W,V'\otimes W')$ by~$\tau_{V',W'} X\tau_{V,W}=\tau_{V',W'}
Y\tau_{V,W}$, and apply the properties preceding equation~(\ref{e:symantisym})
and similar equations for~$\Astd$ and~$\Adot$, then
we obtain the equations of the corollary.
\end{proof}

By equation~(\ref{e:threedecomp}) commutators of elements 
of~$\End_{\gl(1\vert 1)}(V_{t_1}\otimes V_{t_2}\otimes V_{t_3})$
lie in the subspace~$\End_{\gl(1\vert 1)}(V_{t_8})$ 
of~$\End_{\gl(1\vert 1)}(V_{t_1}\otimes V_{t_2}\otimes V_{t_3})$.
We make some explicit computations.

\begin{lemma}\label{l:brIII} 
With $t_i$ and~$\kappa$ as in Lemma~\ref{l:assobas}
we have

\begin{eqnarray*}
& & [(\Ystd_{t_1,t_2}\circ\Astd_{t_1,t_2})\otimes \I_{t_3},
\I_{t_1}\otimes (\Ystd_{t_2,t_3}\circ\Astd_{t_2,t_3})]\\ & = &
[(\Ydot_{t_1,t_2}\circ\Adot_{t_1,t_2})\otimes \I_{t_3},
\I_{t_1}\otimes (\Ydot_{t_2,t_3}\circ\Adot_{t_2,t_3})]\\ & = &
[\I_{t_1}\otimes(\Ydot_{t_2,t_3}\circ\Adot_{t_2,t_3}),
(\Ystd_{t_1,t_2}\circ\Astd_{t_1,t_2}\otimes \I_{t_3})]\\ & = &
[\I_{t_1}\otimes (\Ystd_{t_2,t_3}\circ\Astd_{t_2,t_3}),
(\Ydot_{t_1,t_2}\circ\Adot_{t_1,t_2})\otimes \I_{t_3}]\\ & = &
\frac{(-1)^{\sigma_2}\lambda_1\kappa}
{\lambda_1+\lambda_2}(\Ystd_{t_1,t_2}\otimes
\I_{t_3})\circ \Ydot_{t_4, t_3}\circ\Astd_{t_7,
t_3}\circ(\Adot_{t_1,t_2}\otimes \I_{t_3})\\ & &
-\frac{(-1)^{\sigma_2}\lambda_3}{\lambda_2+\lambda_3}(\Ydot_{t_1,t_2}\otimes
\I_{t_3})\circ \Ystd_{t_7, t_3}\circ\Adot_{t_4,
t_3}\circ(\Astd_{t_1,t_2}\otimes \I_{t_3}).
\end{eqnarray*}
\end{lemma}
\begin{proof} For the first commutator we compute
\begin{eqnarray*} 
& & 
[(\Ystd_{t_1,t_2}\circ\Astd_{t_1,t_2})\otimes \I_{t_3},
\I_{t_1}\otimes (\Ystd_{t_2,t_3}\circ\Astd_{t_2,t_3})]
\\ & = & (\Ystd_{t_1,t_2}\otimes \I_{t_3}) \circ\Ystd_{t_4,t_3}\circ \Astd_{t_1,t_6}\circ
(\I_{t_1}\otimes \Astd_{t_2,t_3})\\ & &
+\frac{(-1)^{\sigma_2}\lambda_1}{\lambda_1+\lambda_2}
(\Ystd_{t_1,t_2}\otimes \I_{t_3}) \circ\Ydot_{t_4,t_3}\circ \Adot_{t_1,t_6}\circ
(\I_{t_1}\otimes \Astd_{t_2,t_3})\\
& & -(\I_{t_1}\otimes \Ystd_{t_2, t_3}) \circ\Ystd_{t_1,t_6}\circ \Astd_{t_4,t_3}\circ
(\Astd_{t_1,t_2}\otimes \I_{t_3})\\ & &
-\frac{(-1)^{\sigma_2}\lambda_3}{\lambda_2+\lambda_3}(\I_{t_1}\otimes \Ystd_{t_2,t_3})
\circ\Ydot_{t_1,t_6}\circ\Adot_{t_4,t_3}\circ(\Astd_{t_1,t_2}\otimes \I_{t_3})
\\ & = & 
\frac{(-1)^{\sigma_2}\lambda_1}{\lambda_1+\lambda_2}
(\Ystd_{t_1,t_2}\otimes \I_{t_3}) \circ\Ydot_{t_4,t_3}\circ \Adot_{t_1,t_6}\circ
(\I_{t_1}\otimes \Astd_{t_2,t_3})\\
& & -\frac{(-1)^{\sigma_2}\lambda_3}{\lambda_2+\lambda_3}(\I_{t_1}\otimes \Ystd_{t_2,t_3})
\circ\Ydot_{t_1,t_6}\circ\Adot_{t_4,t_3}\circ(\Astd_{t_1,t_2}\otimes \I_{t_3})
\\ 
& = & \frac{(-1)^{\sigma_2}\lambda_1}{\lambda_1+\lambda_2}\left(\kappa
(\Ystd_{t_1, t_2}\otimes \I_{t_3})\circ \Ydot_{t_4,t_3}\circ\Astd_{t_7,t_3}\circ 
(\Adot_{t_1,t_2}\otimes \I_{t_3})\right.\\
& & +\frac{(-1)^{\sigma_2}\lambda_3}{\lambda_2+\lambda_3}\left.
(\Ystd_{t_1, t_2}\otimes \I_{t_3})\circ \Ydot_{t_4,t_3}\circ\Adot_{t_4,t_3}\circ 
(\Astd_{t_1,t_2}\otimes \I_{t_3})\right)\\
& & -\frac{(-1)^{\sigma_2}\lambda_3}{\lambda_2+\lambda_3}\left(
(\Ydot_{t_1, t_2}\otimes \I_{t_3})\circ \Ystd_{t_7,t_3}\circ\Adot_{t_4,t_3}\circ 
(\Astd_{t_1,t_2}\otimes \I_{t_3})\right.\\
& & + \frac{(-1)^{\sigma_2}\lambda_1}{\lambda_1+\lambda_2}\left.
(\Ystd_{t_1, t_2}\otimes \I_{t_3})\circ \Ydot_{t_4,t_3}\circ\Adot_{t_4,t_3}\circ 
(\Astd_{t_1,t_2}\otimes \I_{t_3})\right)\\
& = & \frac{(-1)^{\sigma_2}\lambda_1\kappa}
{\lambda_1+\lambda_2}(\Ystd_{t_1,t_2}\otimes
\I_{t_3})\circ \Ydot_{t_4, t_3}\circ\Astd_{t_7,
t_3}\circ(\Adot_{t_1,t_2}\otimes \I_{t_3})\\ & &
-\frac{(-1)^{\sigma_2}\lambda_3}{\lambda_2+\lambda_3}(\Ydot_{t_1,t_2}\otimes
\I_{t_3})\circ \Ystd_{t_7, t_3}\circ\Adot_{t_4,
t_3}\circ(\Astd_{t_1,t_2}\otimes \I_{t_3}),
\end{eqnarray*}

where the first equality follows from Lemma~\ref{l:HI} and Corollary~\ref{c:HI}, the second equality
is a consequence of Lemma~\ref{l:assotriv}, and the third equality is implied by 
Corollary~\ref{c:assobas}
and Lemma~\ref{l:assobas}.
The remaining equations can be proven similarly.
\end{proof}

The following lemma will be used to express the action of~$\Phi$
on~$(V_{t_1}\otimes V_{t_2}
\otimes V_{t_3})[[h]]$ in terms of the commutators of our basis elements from 
Lemma~\ref{l:brIII}.

\begin{lemma}\label{l:WY}
With $t_i=(\lambda_i,\mu_i,\sigma_i)$ as in equations~(\ref{e:lambdai}), (\ref{e:ti}), and

$$
Y\ = \ (-1)^{\sigma_1+\sigma_2+\sigma_3}\setbox1=\hbox{\begin{picture}(0,0)%
\includegraphics{comm.eps}%
\end{picture}%
\setlength{\unitlength}{3947sp}%
\begingroup\makeatletter\ifx\SetFigFont\undefined%
\gdef\SetFigFont#1#2#3#4#5{%
  \reset@font\fontsize{#1}{#2pt}%
  \fontfamily{#3}\fontseries{#4}\fontshape{#5}%
  \selectfont}%
\fi\endgroup%
\begin{picture}(794,344)(579,-158)
\end{picture}
}\vcenter{\box1}\!\in\End_{\Anab}((b_1,b_2,b_3))
\quad \mbox{($b_i=((-1)^{\sigma_i}\lambda_i,(-1)^{\sigma_i}\mu_i,\sigma_i))$}
$$ 

we have

$$W(Y)=(\lambda_1+\lambda_2)(\lambda_2+\lambda_3)
\left[(\Ystd_{t_1,t_2}\circ\Astd_{t_1,t_2})\otimes \I_{t_3},
\I_{t_1}\otimes (\Ystd_{t_2,t_3}\circ\Astd_{t_2,t_3})\right].$$
\end{lemma}
\begin{proof}
With $\omega=\sum_\nu a_\nu\otimes b_\nu$ as in equation~(\ref{e:defomega}) we compute

\begin{eqnarray*}
2\,\sum_\nu (-1)^{\deg(v_{t_1})\deg(b_\nu)} 
a_\nu\cdot v_{t_1}\otimes b_\nu\cdot v_{t_2} &  = & 
(\lambda_1(\mu_2+1)+\lambda_2(\mu_1+1))
v_{t_1}\otimes v_{t_2},\\
2\, \sum_\nu (-1)^{\deg(w_{t_1})\deg(b_\nu)} 
a_\nu\cdot w_{t_1}\otimes b_\nu\cdot w_{t_2}
&  = & (\lambda_1(\mu_2-1)+\lambda_2(\mu_1-1))
w_{t_1}\otimes w_{t_2}.
\end{eqnarray*}

Equation~(\ref{e:defWomega}) implies 

\begin{eqnarray*}
2\,W\!\!\left((-1)^{\sigma_1+\sigma_2}\setbox1=\hbox{\begin{picture}(0,0)%
\includegraphics{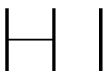}%
\end{picture}%
\setlength{\unitlength}{3947sp}%
\begingroup\makeatletter\ifx\SetFigFont\undefined%
\gdef\SetFigFont#1#2#3#4#5{%
  \reset@font\fontsize{#1}{#2pt}%
  \fontfamily{#3}\fontseries{#4}\fontshape{#5}%
  \selectfont}%
\fi\endgroup%
\begin{picture}(494,344)(579,-158)
\end{picture}
}\vcenter{\box1}\!\right) & \!\!=\!\! &
(a+b)(\Ystd_{t_1,t_2}\circ\Astd_{t_1,t_2})\otimes\I_{t_3}+
(a-b)(\Ydot_{t_1,t_2}\circ\Adot_{t_1,t_2})\otimes \I_{t_3},\\
2\,W\!\!\left((-1)^{\sigma_2+\sigma_3}\setbox1=\hbox{\begin{picture}(0,0)%
\includegraphics{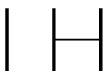}%
\end{picture}%
\setlength{\unitlength}{3947sp}%
\begingroup\makeatletter\ifx\SetFigFont\undefined%
\gdef\SetFigFont#1#2#3#4#5{%
  \reset@font\fontsize{#1}{#2pt}%
  \fontfamily{#3}\fontseries{#4}\fontshape{#5}%
  \selectfont}%
\fi\endgroup%
\begin{picture}(494,344)(1179,-158)
\end{picture}
}\vcenter{\box1}\!\right) & \!\!=\!\! &
(c+d)\I_{t_1}\otimes(\Ystd_{t_2,t_3}\circ\Astd_{t_2,t_3})+
(c-d)\I_{t_1}\otimes(\Ydot_{t_1,t_2}\circ\Adot_{t_1,t_2}),
\end{eqnarray*}

where $a=\lambda_1\mu_2+\lambda_2\mu_1$, $b=\lambda_1+\lambda_2$, 
$c=\lambda_2\mu_3+\lambda_3\mu_2$, and $d=\lambda_2+\lambda_3$.
By Lemma~\ref{l:brIII} we have

\begin{eqnarray*}
W(Y) & = & \left[W\!\!\left((-1)^{\sigma_1+\sigma_2}\setbox1=\hbox{\begin{picture}(0,0)%
\includegraphics{t12.eps}%
\end{picture}%
\setlength{\unitlength}{3947sp}%
\begingroup\makeatletter\ifx\SetFigFont\undefined%
\gdef\SetFigFont#1#2#3#4#5{%
  \reset@font\fontsize{#1}{#2pt}%
  \fontfamily{#3}\fontseries{#4}\fontshape{#5}%
  \selectfont}%
\fi\endgroup%
\begin{picture}(494,344)(579,-158)
\end{picture}
}\vcenter{\box1}\!\right),
W\!\!\left((-1)^{\sigma_2+\sigma_3}\setbox1=\hbox{\begin{picture}(0,0)%
\includegraphics{t23.eps}%
\end{picture}%
\setlength{\unitlength}{3947sp}%
\begingroup\makeatletter\ifx\SetFigFont\undefined%
\gdef\SetFigFont#1#2#3#4#5{%
  \reset@font\fontsize{#1}{#2pt}%
  \fontfamily{#3}\fontseries{#4}\fontshape{#5}%
  \selectfont}%
\fi\endgroup%
\begin{picture}(494,344)(1179,-158)
\end{picture}
}\vcenter{\box1}\!\right)
\right]\\
& = & (1/4)\left((a+b)(c+d)+(a-b)(c-d)
-(a+b)(c-d)\right.\\
& & \left.-(a-b)(c+d)\right) \left[(\Ystd_{t_1,t_2}\circ\Astd_{t_1,t_2})\otimes \I_{t_3},
\I_{t_1}\otimes (\Ystd_{t_2,t_3}\circ\Astd_{t_2,t_3})\right]\\
& = & bd
[(\Ystd_{t_1,t_2}\circ\Astd_{t_1,t_2})\otimes \I_{t_3},
\I_{t_1}\otimes (\Ystd_{t_2,t_3}\circ\Astd_{t_2,t_3})].
\end{eqnarray*}

This completes the proof.
\end{proof}

\section{Proof of Theorem~\ref{t:F}}\label{s:tetra2}

Let $T_3=U_3\circ A_3\circ L_3$ be the colored graph in
Figure~\ref{f:threetetra}.
The upper half~$U_3$ of~$T_3$ is mapped by
$\Ws\circ \Zb$ to the following morphism from~$(V_{t_1}\otimes
V_{t_2}\otimes V_{t_3})[[h]]$ to~$V_{t_8}[[h]]$ (see
equations~(\ref{e:ZYf}), (\ref{e:Zremval}), (\ref{e:Wnu}), (\ref{e:WoAo}))~:

\begin{eqnarray}
\Ws\left(\Zb(U_3)\right)& = &
\varphi(\lambda_8h)^{1/2}\varphi(\lambda_1h)^{-1/2}\varphi(\lambda_6h)^{1/2}\varphi(\lambda_2h)^{1/2}
\Wnab(\Gamma(U_3))\label{e:WZU3}\\
& = &
-\lambda_6\varphi(\lambda_8h)^{1/2}\varphi(\lambda_1h)^{-1/2}\varphi(\lambda_6h)^{1/2}
\varphi(\lambda_2h)^{1/2}
\Adot_{t_1, t_6}\circ(\I_{t_1}\circ\Astd_{t_2, t_3}).\nonumber
\end{eqnarray}

For similar reasons
the lower half~$L_3$ of~$T_3$ is mapped by~$\Ws\circ \Zb$ to the
following morphism from~$V_{t_8}[[h]]$ to~$(V_{t_1}\otimes
V_{t_2}\otimes V_{t_3})[[h]]$~:

\begin{eqnarray}
\Ws\left(\Zb(L_3)\right) & = &
\varphi(\lambda_1h)^{1/2}\varphi(\lambda_7h)^{1/2}\Wnab(\Gamma(L_3))\nonumber\\
& = &
\lambda_7\varphi(\lambda_1h)^{1/2}\varphi(\lambda_7h)^{1/2}
(\Ydot_{t_1, t_2})\otimes\I_{t_3})\circ 
\Ystd_{t_7, t_3}.\label{e:WZL3}
\end{eqnarray}

Since $\End_{\gloo}(V_{t_8})=\Q\,\I_{t_8}$ there exists a formal power series $x\in\Q[[h]]$
satisfying

\begin{equation}\label{e:defx}
x\, \I_{t_8}=
\Adot_{t_1, t_6}\circ(\I_{t_1}\otimes\Astd_{t_2, t_3})
\circ\Ws\left(\Zb(A_3)\right)\circ
(\Ydot_{t_1, t_2}\otimes\I_{t_3})\circ 
\Ystd_{t_7, t_3}.
\end{equation}

We will determine~$x$ in two different ways.
Using equations~(\ref{e:WZU3}),~(\ref{e:WZL3}), and~(\ref{e:defx}) we compute

\begin{eqnarray}
\nabh(S_3)\I_{t_4} & = & \frac{h}{\lambda_8\varphi(\lambda_8h)}
\Ws(\Zb(T_3))
\nonumber\\ & = &
-\frac{x\lambda_6\lambda_7\varphi(\lambda_6h)^{1/2}\varphi(\lambda_2h)^{1/2}\varphi(\lambda_7h)^{1/2}h}
{\lambda_8\varphi(\lambda_8h)^{1/2}} \I_{t_4}.\label{e:WtZT3}
\end{eqnarray}

By Lemma~\ref{l:nabhT} we have~$\nabh(S_3)=e^{-\lambda_2h/2}-e^{\lambda_2h/2}=
-\lambda_2 h\varphi(\lambda_2h)$ for the 
closure~$S_3$ of~$T_3$. Equation~(\ref{e:WtZT3}) implies

\begin{equation}\label{e:x1}
x=\frac{\lambda_2\lambda_8}{\lambda_6\lambda_7}
\sqrt{\frac{\varphi(\lambda_2h)\varphi(\lambda_8h)}{\varphi(\lambda_6h)\varphi(\lambda_7h)}}
=\kappa
\sqrt{\frac{\varphi(v)\varphi(u+v+w)}
{\varphi(u+v)\varphi(v+w)}},
\end{equation}

with~$\kappa=\lambda_2\lambda_8/(\lambda_6\lambda_7)$ as in Lemma~\ref{l:assobas}
and $u=\lambda_1 h, v=\lambda_2h, w=\lambda_3h$.
Now we use equation~(\ref{e:defx}) directly to derive an equation for~$x$ that
depends on a Drinfeld associator. We
start with a general remark. Let~$R$ be a commutative ring
with~$1$. 
Let $M_2(R[[h]])$ be the algebra of $2\times 2$-matrices over~$R[[h]]$.
Let~$a,b\in hR[[h]]$ be elements of the augmentation
ideal of~$R[[h]]$. Then we have

\begin{equation}\label{e:expA}
\exp\left(\left(\begin{array}{ll}0&a\\b&0\end{array}\right)\right)=
\left(\begin{array}{ll}\cosh(c)
& a\sinh(c)/c\\ b\sinh(c)/c & \cosh(c)\end{array}\right)\in M_2(R[[h]])
,\end{equation}

where $c^2=ab$ (the result does not depend on the choice of a
(formal) root~$c$ of~$c^2$ because~$\cosh(c)$ and~$\sinh(c)/c$
are even power series in~$c$). 
Notice that the bases of Lemma~\ref{l:assobas} and Corollary~\ref{c:assobas} 
establish isomorphisms

\begin{equation}\label{e:m2iso} 
M_2(\Q[[h]])\cong\End_{\gl(1\vert 1)}\left(V_{t_8}^{\oplus 2}\right)[[h]]\subset 
\End(V_{t_1}\otimes V_{t_2}\otimes V_{t_3})[[h]].\end{equation}

Let~$F\in\Q[[C,D,E]]\subset\Q[[d_1,d_2,d_3]]$ be the formal power series of 
Theorem~\ref{t:existF2}. Let $Y$ be as in Lemma~\ref{l:WY} with~$\sigma_2=1$.
Lemmas~\ref{l:brIII} and~\ref{l:WY} imply

\begin{eqnarray}
F(u,v,w)\Ws(Y) & = & a\, (\Ydot_{t_1,t_2}\otimes
\I_{t_3})\circ \Ystd_{t_7, t_3}\circ\Adot_{t_4,
t_3}\circ(\Astd_{t_1,t_2}\otimes \I_{t_3})\nonumber\\
& & + \,b\, (\Ystd_{t_1,t_2}\otimes
\I_{t_3})\circ \Ydot_{t_4, t_3}\circ\Astd_{t_7,
t_3}\circ(\Adot_{t_1,t_2}\otimes \I_{t_3}),\\
\mbox{where}\quad a & = & (u+v)\,w\, F(u,v,w),\\
\mbox{and}\quad b & = & -u(v+w)\,\kappa \,
F(u,v,w).\label{e:defb}
\end{eqnarray}

By equations~(\ref{e:defx}), (\ref{e:orinv}),
Corollary~\ref{c:assobas}, and equations~(\ref{e:expA}) to~(\ref{e:defb}) with

\begin{equation}
c^2 = ab = 
-uvw(u+v+w)F(u,v,w)^2
\end{equation}

we have

\begin{eqnarray}
x\, \I_{t_8} & = & \Adot_{t_1, t_6}\circ(\I_{t_1}\otimes\Astd_{t_2, t_3})
\circ\Ws\left(\exp(F(d_1, -d_2, d_3)\cdot Y)\right)\circ
(\Ydot_{t_1, t_2}\otimes\I_{t_3})\circ 
\Ystd_{t_7, t_3}\nonumber\\
& = &
\left(\kappa\Astd_{t_7,t_3}\circ(\Adot_{t_1, t_2}\otimes\I_{t_3})-
\lambda_3/
(\lambda_2+\lambda_3)\Adot_{t_4,t_3}\circ(\Astd_{t_1, t_2}\otimes\I_{t_3})\right)
\nonumber\\
& & \
\circ\exp(F(u,v,w)\Ws(Y))\circ
(\Ydot_{t_1, t_2}\otimes\I_{t_3})\circ 
\Ystd_{t_7, t_3}\nonumber\\ & = &
(\kappa\cosh(c)
-(\lambda_3b/(\lambda_2+\lambda_3))\sinh(c)/c)\I_{t_8}\nonumber\\
& = & \kappa(\cosh(c)+
uw F(u,v,w)\sinh(c)/c)\I_{t_8}.\label{e:x2}
\end{eqnarray}

Equations~(\ref{e:x1}) and~(\ref{e:x2}) imply 

\begin{equation}
\cosh(c)+uw\sinh(c)/(c/F(u,v,w))=\sqrt{\frac{\varphi(v)\varphi(u+v+w)}{\varphi(u+v)\varphi(v+w)}}.
\end{equation}

Since this formula holds for arbitrary values of~$u, v, w\in\Q^*h$
with $u+v\not=0$, $v+w\not=0$, and
$u+v+w\not=0$ we have proven the equation between formal
power series stated in Theorem~\ref{t:F}.

\section{$\nabh$ and Viro's Alexander
invariant}\label{s:MultAlexViro}

In this section we relate the Alexander series~$\nabh$ to Viro's Alexander 
invariant~$\underline{\Delta}^1$. Proofs that are direct 
translations of proofs from~\cite{Vir} will only be sketched in what follows.
We start with deriving relations between values of~$\nabh$ on different
trivalent framed graphs with admissible coloring.
When the colors and orientations 
of the lower three edges of~$\setbox1=\hbox{\begin{picture}(0,0)%
\includegraphics{o.eps}%
\end{picture}%
\setlength{\unitlength}{3947sp}%
\begingroup\makeatletter\ifx\SetFigFont\undefined%
\gdef\SetFigFont#1#2#3#4#5{%
  \reset@font\fontsize{#1}{#2pt}%
  \fontfamily{#3}\fontseries{#4}\fontshape{#5}%
  \selectfont}%
\fi\endgroup%
\begin{picture}(130,294)(761,-58)
\end{picture}
}\vcenter{\box1}\!$ are fixed, then
by equation~(\ref{e:admissibilitya})
there are two possible colors of the upper edge in an admissible coloring
of a graph. 
In the first case the upper and lower edge point into the same
direction, and we 
compute

\begin{equation}\label{e:nabhbub}
\nabh\left(\setbox1=\hbox{\input{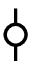}}\vcenter{\box1}\!\right)=\nabh(\setbox1=\hbox{\begin{picture}(0,0)%
\includegraphics{oc.eps}%
\end{picture}%
\setlength{\unitlength}{3947sp}%
\begingroup\makeatletter\ifx\SetFigFont\undefined%
\gdef\SetFigFont#1#2#3#4#5{%
  \reset@font\fontsize{#1}{#2pt}%
  \fontfamily{#3}\fontseries{#4}\fontshape{#5}%
  \selectfont}%
\fi\endgroup%
\begin{picture}(210,209)(651,-41)
\put(856,-41){\makebox(0,0)[lb]{\smash{\SetFigFont{6}{7.2}{\familydefault}{\mddefault}{\updefault}{\color[rgb]{0,0,0}$c$}%
}}}
\end{picture}
}\vcenter{\box1}\!)^{-1}\nabh\left(\setbox1=\hbox{\begin{picture}(0,0)%
\includegraphics{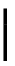}%
\end{picture}%
\setlength{\unitlength}{3947sp}%
\begingroup\makeatletter\ifx\SetFigFont\undefined%
\gdef\SetFigFont#1#2#3#4#5{%
  \reset@font\fontsize{#1}{#2pt}%
  \fontfamily{#3}\fontseries{#4}\fontshape{#5}%
  \selectfont}%
\fi\endgroup%
\begin{picture}(52,294)(804,-58)
\put(856,-41){\makebox(0,0)[lb]{\smash{\SetFigFont{6}{7.2}{\familydefault}{\mddefault}{\updefault}{\color[rgb]{0,0,0}$c$}%
}}}
\end{picture}
}\vcenter{\box1}\!\right)
\end{equation}

by using equations~(\ref{e:ZYf}), (\ref{e:Zremval}), 
(\ref{e:WstWt}), (\ref{e:nabhOc}), (\ref{e:Wnu}), and 
Lemma~\ref{l:cutnotimp1}. In the second case $\widehat{W}\circ Z_{\widehat{\C}}$ 
maps~$\setbox1=\hbox{\begin{picture}(0,0)%
\includegraphics{o.eps}%
\end{picture}%
\setlength{\unitlength}{3947sp}%
\begingroup\makeatletter\ifx\SetFigFont\undefined%
\gdef\SetFigFont#1#2#3#4#5{%
  \reset@font\fontsize{#1}{#2pt}%
  \fontfamily{#3}\fontseries{#4}\fontshape{#5}%
  \selectfont}%
\fi\endgroup%
\begin{picture}(130,294)(761,-58)
\end{picture}
}\vcenter{\box1}\!$ to~$x\,\psi$ where~$x\in\Q[[h]]$ and~$\psi$ is a 
morphism between non-isomorphic simple $\gl(1\vert 1)$-modules. 
Therefore, we obtain 

\begin{equation}\label{e:nabhbub0}
\nabh\left(\setbox1=\hbox{\input{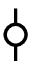}}\vcenter{\box1}\!\right)=0\quad\mbox{if $c\not=d$.}
\end{equation}

Consider two parallel strands $\setbox1=\hbox{\input{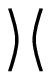}}\vcenter{\box1}\!$ colored by
$a=(\lambda,\mu)$, $b=(\lambda',\mu')$ and define $s\in\{\pm 1\}$ (resp.\
$s'\in\{\pm 1\}$) iff the left (resp.\ right) strand
points downwards. Then we have

\begin{equation}\label{e:nabhii}
\nabh\left(\setbox1=\hbox{\input{iiab}}\vcenter{\box1}\!\right)=\sum_c\nabh\left(\setbox1=\hbox{\begin{picture}(0,0)%
\includegraphics{oc.eps}%
\end{picture}%
\setlength{\unitlength}{3947sp}%
\begingroup\makeatletter\ifx\SetFigFont\undefined%
\gdef\SetFigFont#1#2#3#4#5{%
  \reset@font\fontsize{#1}{#2pt}%
  \fontfamily{#3}\fontseries{#4}\fontshape{#5}%
  \selectfont}%
\fi\endgroup%
\begin{picture}(210,209)(651,-41)
\put(856,-41){\makebox(0,0)[lb]{\smash{\SetFigFont{6}{7.2}{\familydefault}{\mddefault}{\updefault}{\color[rgb]{0,0,0}$c$}%
}}}
\end{picture}
}\vcenter{\box1}\!\right)\nabh\left(\setbox1=\hbox{\input{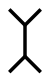}}\vcenter{\box1}\!\right)\quad
\mbox{if $s\lambda+s'\lambda'\not=0$,}
\end{equation}

where the sum runs over the two colors~$c$ such that the coloring of
$\setbox1=\hbox{\input{iab}}\vcenter{\box1}\!$ is admissible. In proofs of 
equation~(\ref{e:nabhii}) and equation~(\ref{e:nabhH}) below, 
we use equation~(\ref{e:tensordecomp}) to show that the left sides of these
equations are equal to linear combinations

$$
\sum_p x_p\nabh\left(\setbox1=\hbox{\begin{picture}(0,0)%
\includegraphics{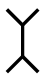}%
\end{picture}%
\setlength{\unitlength}{3947sp}%
\begingroup\makeatletter\ifx\SetFigFont\undefined%
\gdef\SetFigFont#1#2#3#4#5{%
  \reset@font\fontsize{#1}{#2pt}%
  \fontfamily{#3}\fontseries{#4}\fontshape{#5}%
  \selectfont}%
\fi\endgroup%
\begin{picture}(194,344)(279,-83)
\put(406, 59){\makebox(0,0)[lb]{\smash{\SetFigFont{6}{7.2}{\familydefault}{\mddefault}{\updefault}{\color[rgb]{0,0,0}$p$}%
}}}
\end{picture}
}\vcenter{\box1}\!\right),$$ 

for certain~$x_p\in\Q[[h]]$
and we use equation~(\ref{e:nabhbub}) 
to determine the coefficients~$x_p$ (see
the proof of~9.2.A in~\cite{Vir} for more details).

Now consider $\setbox1=\hbox{\input{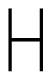}}\vcenter{\box1}\!$ with colors~$a,b,c,d,e\in\Q^*\times
\Q$. Let $a=(\lambda, \mu)$, $b=(\lambda',\mu')$, and define $s,
s'\in\{\pm 1\}$ as above. We assume that
$s\lambda+s'\lambda'\not=0$. Then

\begin{equation}\label{e:nabhH}
\nabh\left(\setbox1=\hbox{\input{hae}}\vcenter{\box1}\!\right)=\sum_f\nabh\left(\setbox1=\hbox{\input{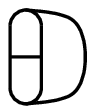}}\vcenter{\box1}\!\right)\nabh\left(\setbox1=\hbox{\begin{picture}(0,0)%
\includegraphics{oc.eps}%
\end{picture}%
\setlength{\unitlength}{3947sp}%
\begingroup\makeatletter\ifx\SetFigFont\undefined%
\gdef\SetFigFont#1#2#3#4#5{%
  \reset@font\fontsize{#1}{#2pt}%
  \fontfamily{#3}\fontseries{#4}\fontshape{#5}%
  \selectfont}%
\fi\endgroup%
\begin{picture}(210,209)(651,-41)
\put(856,-41){\makebox(0,0)[lb]{\smash{\SetFigFont{6}{7.2}{\familydefault}{\mddefault}{\updefault}{\color[rgb]{0,0,0}$c$}%
}}}
\end{picture}
}\vcenter{\box1}\!\right)\nabh\left(\setbox1=\hbox{\input{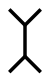}}\vcenter{\box1}\!\right),
\end{equation}

where the sum runs over the one or two colors~$f$ such that the
coloring of~$\setbox1=\hbox{\input{iaf}}\vcenter{\box1}\!$ is admissible and~$\setbox1=\hbox{\input{hae}}\vcenter{\box1}\!\subset\setbox1=\hbox{\input{tetaf}}\vcenter{\box1}\!$
coincides with~$\setbox1=\hbox{\input{hae}}\vcenter{\box1}\!$ as oriented colored graph.
Let~$i$ be the edge colored~$c$ in~$\setbox1=\hbox{\input{hae}}\vcenter{\box1}\!$.
When the restriction $s\lambda+s'\lambda'\not=0$ is satisfied (resp.\ violated) we
say that equation~(\ref{e:nabhH}) can (resp.\ cannot) be applied to the edge~$i$.

For a strand colored by~$c=(\lambda, \mu)$ with a
right-handed half-twist, it follows by a direct computation that

\begin{equation}\label{e:nabhhtwist}
\nabh\left(\setbox1=\hbox{\begin{picture}(0,0)%
\includegraphics{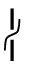}%
\end{picture}%
\setlength{\unitlength}{3947sp}%
\begingroup\makeatletter\ifx\SetFigFont\undefined%
\gdef\SetFigFont#1#2#3#4#5{%
  \reset@font\fontsize{#1}{#2pt}%
  \fontfamily{#3}\fontseries{#4}\fontshape{#5}%
  \selectfont}%
\fi\endgroup%
\begin{picture}(94,294)(779,-58)
\put(856,-41){\makebox(0,0)[lb]{\smash{\SetFigFont{6}{7.2}{\familydefault}{\mddefault}{\updefault}{\color[rgb]{0,0,0}$c$}%
}}}
\end{picture}
}\vcenter{\box1}\!\right)=e^{\lambda\mu h/4}\nabh\left(\setbox1=\hbox{\begin{picture}(0,0)%
\includegraphics{ic.eps}%
\end{picture}%
\setlength{\unitlength}{3947sp}%
\begingroup\makeatletter\ifx\SetFigFont\undefined%
\gdef\SetFigFont#1#2#3#4#5{%
  \reset@font\fontsize{#1}{#2pt}%
  \fontfamily{#3}\fontseries{#4}\fontshape{#5}%
  \selectfont}%
\fi\endgroup%
\begin{picture}(52,294)(804,-58)
\put(856,-41){\makebox(0,0)[lb]{\smash{\SetFigFont{6}{7.2}{\familydefault}{\mddefault}{\updefault}{\color[rgb]{0,0,0}$c$}%
}}}
\end{picture}
}\vcenter{\box1}\!\right).
\end{equation}

At a trivalent vertex~$v$ of~$G$ we
obtain from part~(2) of Theorem~\ref{t:Ztriv}
and equation~(\ref{e:symantisym}) that

\begin{equation}
\nabh\left(\setbox1=\hbox{\begin{picture}(0,0)%
\includegraphics{y1.eps}%
\end{picture}%
\setlength{\unitlength}{3947sp}%
\begingroup\makeatletter\ifx\SetFigFont\undefined%
\gdef\SetFigFont#1#2#3#4#5{%
  \reset@font\fontsize{#1}{#2pt}%
  \fontfamily{#3}\fontseries{#4}\fontshape{#5}%
  \selectfont}%
\fi\endgroup%
\begin{picture}(344,394)(729,-133)
\end{picture}
}\vcenter{\box1}\!\right)=s_v\nabh\left(\setbox1=\hbox{\begin{picture}(0,0)%
\includegraphics{y2.eps}%
\end{picture}%
\setlength{\unitlength}{3947sp}%
\begingroup\makeatletter\ifx\SetFigFont\undefined%
\gdef\SetFigFont#1#2#3#4#5{%
  \reset@font\fontsize{#1}{#2pt}%
  \fontfamily{#3}\fontseries{#4}\fontshape{#5}%
  \selectfont}%
\fi\endgroup%
\begin{picture}(344,394)(204,-133)
\end{picture}
}\vcenter{\box1}\!\right).
\label{e:nabhsymasym}
\end{equation}

We have collected all
properties of~$\nabh$ that are needed to state the following proposition.

\begin{prop}\label{p:nabunique}
The invariant~$\nabh$ of embedded colored framed trivalent
graphs~$G$ is uniquely determined by its value on the trivial knot 
in equation~(\ref{e:nabhOc}), its values on planar tetrahedra in 
Lemma~\ref{l:nabhT}, and by the skein relations in equations~(\ref{e:nabhbub})
to~(\ref{e:nabhsymasym}).
\end{prop}
\begin{proof}
Consider a diagram of~$G$. We use equations~(\ref{e:nabhii}) and
(\ref{e:nabhsymasym}) to
replace each crossing in that diagram by a planar graph with
two trivalent vertices (notice that despite the restriction in equation~(\ref{e:nabhii})
this is always possible).
By equation~(\ref{e:nabhhtwist}) we may assume that the resulting planar 
graph has blackboard framing.

Let $n$ be the number of connected 
components plus the number of trivalent vertices of a planar trivalent graph~$G$.
The connected components~$F$ of~$\R^2\setminus
G$ are called the faces of the diagram, and the trivalent vertices
(resp.\ the edges) in the closure of a face~$F$ are called the
vertices (resp.\ edges) of that face. Let~$\ell$ be the minimal
number of vertices of a face of~$G$. 
\footnote{For reasons of Euler
characteristic we have~$\ell\leq 5$ (see~\cite{BeS}),
but we do not need this in the proof.} 

We will prove the proposition by induction on the pairs~$(n,\ell)\in\N\times\N_0$ 
with lexicographical order.
For~$\ell=0$ and~$n=1$ the graph~$G$ is a trivial knot.
For~$\ell=0$ and~$n>1$, we use
equation~(\ref{e:nabhii}), the
property~$\sum_c\nabh(\setbox1=\hbox{\begin{picture}(0,0)%
\includegraphics{oc.eps}%
\end{picture}%
\setlength{\unitlength}{3947sp}%
\begingroup\makeatletter\ifx\SetFigFont\undefined%
\gdef\SetFigFont#1#2#3#4#5{%
  \reset@font\fontsize{#1}{#2pt}%
  \fontfamily{#3}\fontseries{#4}\fontshape{#5}%
  \selectfont}%
\fi\endgroup%
\begin{picture}(210,209)(651,-41)
\put(856,-41){\makebox(0,0)[lb]{\smash{\SetFigFont{6}{7.2}{\familydefault}{\mddefault}{\updefault}{\color[rgb]{0,0,0}$c$}%
}}}
\end{picture}
}\vcenter{\box1}\!)=0$ where the sum runs over the
same values of~$c$ as the sum in equation~(\ref{e:nabhii}), and
equation~(\ref{e:nabhbub}) to show that~$\nabh(G)=0$ in this case.
For~$\ell=1$ the graph~$G$
cannot have an admissible coloring, so we do not need to consider
that case. In the case~$\ell=2$ we can apply
equation~(\ref{e:nabhbub}) which will reduce~$n$ by~$2$, or we can apply
equation~(\ref{e:nabhbub0}) to show that~$\nabh(G)=0$.
Now let $\ell\geq 3$.

{\em Case 1:} Assume that equation~(\ref{e:nabhH}) can be applied to an edge of~$F$.
Then we can reduce~$\ell$ by one while preserving~$n$.

{\em Case 2:}
When equation~(\ref{e:nabhH}) 
cannot be applied to an edge of~$F$
\footnote{In that case~$\ell$ must be~$\geq 4$ and even.}
we choose a trivalent
vertex~$v$ of~$G$ that is connected to a vertex of the face~$F$ by
an edge~$e$, and that is not itself a vertex of~$F$. Such a vertex~$v$ 
exists because~$\ell>2$ was minimal.

{\em Case 2a:\ } We assume that equation~(\ref{e:nabhH}) can be applied to~$e$. 
Then we proceed as shown schematically in Figure~\ref{f:p1m2}.
Equation~(\ref{e:nabhH}) can be applied to the two edges~$i$, $j$ 
in the second step in Figure~\ref{f:p1m2} because 
equation~(\ref{e:nabhH}) could not be applied to an edge of~$F$ 
in the first picture.
This way we again decrease~$\ell$ by~$1$ while preserving~$n$. 

\begin{figure}[!h]
\centering 
\mbox{$\setbox1=\hbox{\input{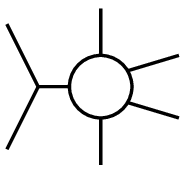}}\vcenter{\box1}\!\ \leadsto\ 
\setbox1=\hbox{\input{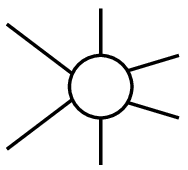}}\vcenter{\box1}\!\ \leadsto\
\setbox1=\hbox{\input{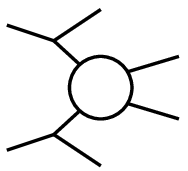}}\vcenter{\box1}\!$}
\caption{Reducing the number of vertices of a face} \label{f:p1m2}
\end{figure}

{\em Case 2b:\ } When equation~(\ref{e:nabhH}) cannot be applied to~$e$, we proceed
as follows.
We first use equation~(\ref{e:nabhbub}) to add a bubble to
one edge of~$F$. This will increase the number of edges and vertices of~$F$ by two.
Let~$e$ be the new edge of~$F$ that belongs to the bubble, let~$(\lambda,0)$ be
the color of~$e$, and let~$v$ be one vertex of~$e$.
As shown in Figure~\ref{f:gencol}, our plan is to
apply equation~(\ref{e:nabhH})~$\ell$ times to
push~$v$ around~$F$, and then to
use equation~(\ref{e:nabhbub}) again to remove the bubble. 

\begin{figure}[!h]
\centering 
\mbox{$\setbox1=\hbox{\input{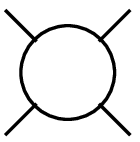}}\vcenter{\box1}\!\ \leadsto\ 
\setbox1=\hbox{\input{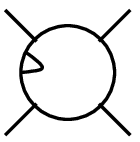}}\vcenter{\box1}\!\ \leadsto\
\setbox1=\hbox{\input{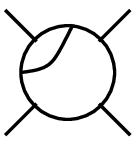}}\vcenter{\box1}\!\leadsto\ \ldots\ \leadsto
\setbox1=\hbox{\input{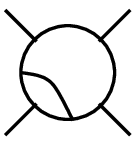}}\vcenter{\box1}\!\ \leadsto\
\setbox1=\hbox{\input{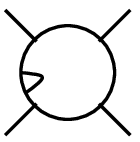}}\vcenter{\box1}\!\ \leadsto\
\setbox1=\hbox{\input{chcol0}}\vcenter{\box1}\!$}
\caption{Changing the colors of the edges of a face} \label{f:gencol}
\end{figure}

This will
express~$\nabh(G)$ as a linear combination of values~$\nabh(G')$ on
diagrams~$G'$ of the same
shape as~$G$, 
but where orientations of the edges of the face~$F$ may have changed and where
the colors of these edges have
changed by additive constants~$(\pm\lambda,\pm 1)$. There are infinitely
many possible choices of~$\lambda\in\Q^*$ such that equation~(\ref{e:nabhH}) 
can be applied~$\ell$ times as needed above. For any such choice of~$\lambda$, we can
apply case~2a to all diagrams~$G'$ as above. This completes the proof.
\end{proof}

Let $B=\Q[[h]][h^{-1}]$ be the quotient field of~$\Q[[h]]$, $M=\exp(\Q h)\subset B^*$, and
\mbox{$W=\Q$.} Define~$\beta:M\times W\longrightarrow M$ by~$\beta(m,w)=m^w$.
Let~$G$ be a trivalent graph with admissible coloring.
Define the colored graph~$q(G)$ by replacing all colors~$(\lambda,\mu)$ 
by~$(\exp(\lambda h/4),\mu)\in M\times W$.
Then the colors of~$q(G)$ verify condition~2.8.A of~\cite{Vir} for the
$1$-palette $P=(B,M,W,\beta)$.
Viro's Alexander invariant~$\underline{\Delta}^1$ (see Section~6.3 of~\cite{Vir}), 
considered as a map
$G\mapsto\underline{\Delta}^1(q(G))$, verifies the same equations as~$\nabh$
in Proposition~\ref{p:nabunique} (see~7.2.D and Section~9 of~\cite{Vir}). 
This implies the following theorem.

\begin{theorem}\label{t:nabvir}
For an admissibly colored framed trivalent graph~$G\subset S^3$
we have

$$\nabh(G)=\underline{\Delta}^1(q(G)).$$
\end{theorem}

More generally, the invariants~$\nabh(G)$ for all admissible colorings of~$G$ determine
Viro's Alexander invariant of graphs
with 'universal' colors as in Section~7.4 of~\cite{Vir}.
Theorem~\ref{t:nabvir} and~7.7.G of~\cite{Vir} imply the following relation 
between~$\nabh(L)$ and the multi-variable Alexander polynomial~$\nabla_L$.

\begin{coro}\label{c:natnab}
For a link $L$ with $n$ components colored by~$(\lambda_i,0)$
we have

$$\nabla_L(e^{\lambda_1h/2},\ldots,e^{\lambda_nh/2})=\nabh(L).$$
\end{coro}

Corollary~\ref{c:natnab} can also be proven directly by using the
characterization of~$\nabla$ by axioms from~\cite{Tur}. Except for
some technical details
this direct proof is a
translation of the proof given in~\cite{Vir}. 
I learned about
Corollary~\ref{c:natnab} from A.\ Vaintrob, but he never published his proof.
Since~$\nabh$ was defined using the universal Vassiliev invariant~$Z$, 
Corollary~\ref{c:natnab} implies that
the coefficient of~$h^k$ in~$\nabla_L(e^{\lambda_1h/2},\ldots,e^{\lambda_nh/2})$
is a Vassiliev invariant of degree~$k+1$ 
(this can also be proven directly using~\cite{Har} or~\cite{Tur}).


\begin{thebibliography}{BGRT}

\bibitem[BN1]{BN1}
D.\ Bar-Natan, {\em On the Vassiliev knot invariants}, Topology {\bf 34}
(1995), 423--472.

\bibitem[BN2]{BN2}
D.\ Bar-Natan, {\em Non--associative tangles}, Geometric topology
proceedings of the Georgia International Topology Conference (W.\
H.\ Kazez ed.),  139--183, Amer.\ Math.\ Soc.\ and international
Press, Providence (1997).

\bibitem[BN3]{BN3} D.\ Bar-Natan, {\em Vassiliev and quantum invariants of braids},
Proc.\ of Symp.\ in Appl.\ Math.\ {\bf 51} (1996), 129--144.

\bibitem[BN4]{BN4} D.\ Bar-Natan, {\em On Associators and the Grothendieck-Teichm\"uller
Group I}, Sel. Math., New Ser.\ {\bf 4}, No.\ 2 (1998), 183--212.

\bibitem[BNG]{BNG} D.\ Bar-Natan and S.\ Garoufalidis, {\em On the Melvin-Morton-Rozansky conjecture},
Invent.\ Math.\ {\bf 125} (1996), 103--133.

\bibitem[BeS]{BeS} A.-B.\ Berger and I.\ Stassen, {\em The skein
relation for the $({\mathfrak g}_2, V)$-link invariant}, Comment.\
Math.\ Helv.\ {\bf 75}, No.\ 1 (2000), 134--155.

\bibitem[Dr1]{Dr1} V.\ G.\ Drinfeld, {\em Quasi-Hopf algebras}, Leningrad Math.\ J.\ {\bf 1} 
(1990), 1419--1457.

\bibitem[Dr2]{Dri} V.\ G.\ Drinfeld, {\em On quasitriangular quasi-Hopf algebras and a group closely
connected with $\Gal(\overline{\Q}/\Q)$}, Algebra i Analiz {\bf 2:4}
(1990), 149--181. English transl.: Leningrad Math.\ J.\ {\bf 2} (1991),
829--860.

\bibitem[FKV]{FKV} J.\ M.\ Figueroa-O'Farrill, T.\ Kimura, A.\ Vaintrob,
{\em The universal Vassiliev invariant for the Lie superalgebra $\gloo$}, Commun.\ Math.\ 
Phys.\ {\bf 185} (1997), 93--127.

\bibitem[Har]{Har} R.\ Hartley, {\em The Conway potential function for links}, Comment.\ Math.\ 
Helv.\ {\bf 58} (1983), 365--378.

\bibitem[Kac]{Kac} V.\ C.\ Kac, {\em A sketch of Lie superalgebra theory}, 
Comm.\ Math.\ Phys.\ {\bf 53} (1977), 31--64.

\bibitem[Kas]{Kas} C.\ Kassel, {\em Quantum groups},
GTM 155, Springer-Verlag, New York 1995.

\bibitem[Le]{Le} T.\ Q.\ T.\ Le, {\em On denominators of the
Kontsevich integral and the universal perturbative invariant of
$3$-manifolds}, Invent.\ Math.\ {\bf 135}, No.\ 3 (1999), 689--722.

\bibitem[LM1]{LeM} T.\ Q.\ T.\ Le and J. Murakami, {\em The universal
Vassiliev-Kontsevich invariant for framed oriented links}, Comp.\ 
Math.\ {\bf 102} (1996), 41--64. 

\bibitem[LM2]{LM2} T.\ Q.\ T.\ Le and J.\ Murakami, {\em Parallel version of the
universal Vassiliev-Kontsevich invariant}, J.\ Pure and Appl.\ Algebra {\bf 121}, No.\ 3 (1997),
271--291.

\bibitem[Les]{Les} C.\ Lescop, {\em About the uniqueness and the
denominators of the Kontsevich integral}, CNRS Institut Fourier
preprint (2000), math.GT/0004094.


\bibitem[MiM]{MiM} J.\ W.\ Milnor and J.\ C.\ Moore, {\em On the structure of Hopf algebras},
Ann.\ Math.\ {\bf 81}, No.\ 2 (1965), 211--264.

\bibitem[Mur]{Mur} J.\ Murakami, {\em A state model for the multi-variable Alexander 
polynomial}, Pac.\ J.\ Math.\ {\bf 157}, No.\ 1 (1993), 109--135.

\bibitem[MuO]{MuO} J.\ Murakami and T.\ Ohtsuki, {\em Topological quantum field theory
for the universal quantum invariant}, Commun.\ Math.\ Phys.\ {\bf 188}, 
No.\ 3 (1997), 501--520.

\bibitem[Thu]{Thu} D.\ P.\ Thurston, {\em Wheeling~: A diagrammatic analogue of the Duflo
isomorphism}, Ph.\ D.\ Thesis, U.\ C.\ Berkeley (2000), math.QA/0006083.

\bibitem[Tu1]{Tur} V.\ Turaev, {\em Reidemeister torsion in knot
theory}, Russian Math.\ Surveys {\bf 41}, No.\ 1 (1986), 119--182.

\bibitem[Tu2]{Tu2} V.\ Turaev, {\em Quantum invariants of knots and $3$-manifolds}, 
Walter de Gruyter, Berlin, New York 1994.

\bibitem[Vas]{Vas} V.\ A.\ Vassiliev, {\em Cohomology of knot spaces}, Theory of singularities
and its applications (V.\ I.\ Arnold, ed.), Amer.\ Math.\ Soc., Providence (1990).

\bibitem[Vir]{Vir} O.\ Viro, {\em Quantum relatives of Alexander
polynomial}, in preparation.
\end{thebibliography}
\end{document}